\documentclass{amsart}

\newtheorem{theorem}{Theorem}[section]

\theoremstyle{definition}

\theoremstyle{remark}

\newtheorem{assumption}{Assumption}

\numberwithin{equation}{section}

\usepackage{lipsum}
\usepackage{amsfonts}
\usepackage{graphicx}
\usepackage{epstopdf}
\usepackage{algorithmic}
\ifpdf
  \DeclareGraphicsExtensions{.eps,.pdf,.png,.jpg}
\else
  \DeclareGraphicsExtensions{.eps}
\fi
\usepackage{amsmath,amssymb, amsfonts}
\usepackage{float}
\usepackage{graphicx}
\usepackage{booktabs} 
\usepackage[table]{xcolor}
\usepackage{bm}
\usepackage{setspace}
\usepackage{amsthm}
\usepackage{enumitem}
\usepackage{graphicx}
\usepackage{algorithm}
\usepackage[mathscr]{euscript}
\usepackage{cite}
\usepackage{mathtools}

\renewcommand{\vec}[1] {\ensuremath{\boldsymbol{#1}}}

\usepackage{amsopn}

\def\R{{\mathbb R}}
\def\U{{\mathcal U}}
\def\Z{{\mathcal Z}}
\def\t{{\vec{\theta}}}
\def\u{{\vec{u}}}

\def\z{{\vec{z}}}
\def\zbar{{\overline{\z}}}
\def\ztilde{{\tilde{\z}}}

\def\e{{\vec{e}}}

\def\d{{\delta}}
\def\A{{\vec{A}}}
\def\B{{\vec{B}}}
\def\E{{\vec{E}}}
\def\M{{\vec{M}}}
\def\K{{\vec{K}}}
\def\G{{\vec{G}}}

\def\L{{\vec{L}}}
\def\C{{\vec{C}}}

\def\H{{\vec{H}}}
\def\P{{\vec{P}}}

\def\X{{\vec{X}}}
\def\V{{\vec{V}}}

\def\I{{\vec{I}}}

\def\W{{\vec{W}}}

\def\F{{\vec{F}}}
\DeclareMathAlphabet\mathbfcal{OMS}{cmsy}{b}{n}

\usepackage{array}
\newcolumntype{P}[1]{>{\centering\arraybackslash}p{#1}}
\newcolumntype{M}[1]{>{\centering\arraybackslash}m{#1}}

\usepackage{amsaddr}

\begin{document}

\title[Posterior Optimal Solution Sampling]{Hyper-differential sensitivity analysis with respect to model discrepancy: Posterior Optimal Solution Sampling}


\author{Joseph Hart}
\address{Scientific Machine Learning, Sandia National Laboratories, P.O. Box 5800, Albuquerque, NM 87123-1320}
\email{joshart@sandia.gov} 

\author{Bart van Bloemen Waanders}
\address{Scientific Machine Learning, Sandia National Laboratories, P.O. Box 5800, Albuquerque, NM 87123-1320}
\email{bartv@sandia.gov}



\date{}

\dedicatory{}

\begin{abstract}
Optimization constrained by high-fidelity computational models has potential for transformative impact. However, such optimization is frequently unattainable in practice due to the complexity and computational intensity of the model. An alternative is to optimize a low-fidelity model and use limited evaluations of the high-fidelity model to assess the quality of the solution. This article develops a framework to use limited high-fidelity simulations to update the optimization solution computed using the low-fidelity model. Building off a previous article~\cite{hart_bvw_cmame}, which introduced hyper-differential sensitivity analysis with respect to model discrepancy, this article provides novel extensions of the algorithm to enable uncertainty quantification of the optimal solution update via a Bayesian framework. Specifically, we formulate a Bayesian inverse problem to estimate the model discrepancy and propagate the posterior model discrepancy distribution through the post-optimality sensitivity operator for the low-fidelity optimization problem. We provide a rigorous treatment of the Bayesian formulation, a computationally efficient algorithm to compute posterior samples, a guide to specify and interpret the algorithm hyper-parameters, and a demonstration of the approach on three examples which highlight various types of discrepancy between low and high-fidelity models.  
\end{abstract}

\maketitle

\section{Introduction}
Computational models have made significant advances toward enabling the exploration of optimization and uncertainty quantification problems. However, optimization capabilities remain elusive in many models. Gradient-based optimization is needed to explore high-dimensional spaces and account for uncertainties. Yet accessing derivatives in complex simulation codes is challenging. Furthermore, the number of evaluations required by the most efficient optimization algorithms is still impractical in many cases.

To make optimization practical, a low-fidelity model may be constructed for which optimization algorithms are more amenable. Examples of low-fidelity models include: averaging small length scales, omission of physical processes, projection into lower dimensions, and data-driven model approximation. Although such simplifications help to enable optimization, the resulting solutions are suboptimal as a result of modeling errors.

To take advantage of low-fidelity models and mitigate modeling errors, a novel approach was introduced in~\cite{hart_bvw_cmame} that uses hyper-differential sensitivity analysis (HDSA) to learn an improved optimization solution using a limited number of high-fidelity simulations. The approach in~\cite{hart_bvw_cmame} was restricted to considering a mean optimal solution update without quantifying the certainty (or conversely the uncertainty) of the solution. This article provides uncertainty quantification by developing a Bayesian framework as part of the updating process. Our contributions include:
\begin{enumerate}
\item[$\bullet$] A different formulation of the prior and noise models that address shortcomings of the previous formulation.
\item[$\bullet$] Mathematical developments to enable efficient posterior sampling. 
\item[$\bullet$] A novel formulation of projected optimal solution updates which address theoretical concerns and improves computational efficiency.
\item[$\bullet$] Mathematical analysis to interpret hyper-parameters (parameters defining the prior distribution) and guide their specification. 
\end{enumerate}

Our approach is related to multi-fidelity methods, which seek to perform outer loop analysis such as optimization and uncertainty quantification by using low-fidelity model evaluations to reduce the total number of high-fidelity model evaluations required in the analysis~\cite{multifidelity_review_peherstorfer}. Monte Carlo type analyses for general uncertainty quantification applications can be accelerated using such multi-fidelity methods~\cite{multifidelity_gorodetsky}. Furthermore, optimization under uncertainty may be performed more efficiently by using multi-fidelity estimators to compute objective function statistics~\cite{multi_fidelity_opt_willcox,multifidelity_wind_plant_geraci}. In the context of deterministic optimization,~\cite{multifidelity_quasinewton_bryson,trmm_lewis,biegler_rom_opt,willcox_multifi_opt_2012} have explored the use of multi-fidelity algorithms by leveraging both low and high-fidelity model evaluations within the optimization iterates. Similar to these works, we assume access to both a low and high-fidelity model with the goal of determining the minimizer of the high-fidelity optimization problem. However, our approach differs in that we focus on cases where the high-fidelity model evaluations are computationally expensive and must be computed offline rather than being called from within the optimization algorithm. Furthermore, we assume that derivatives of the high-fidelity model are not available. This is motivated by the myriad of practical challenges associated with running complex computational models on high performance computing platforms (e.g. parallel adjoint computations). Although our approach does not provide theoretical guarantees of convergence to high-fidelity minima, our framework is pragmatic and can provide improvements in optimal solutions for many applications where multi-fidelity optimization is infeasible due to the complexity of the high-fidelity model.

Defined as the difference between model and reality, model discrepancy has been studied extensively in Bayesian frameworks. Pioneering work by Kennedy and O'Hagan~\cite{ohagan2001} and extensions~\cite{ohagan2002,ohagan2014,Higdon_2008,Arendt_2012,Maupin,Ling_2014,Gardner_2021} explore model calibration in the presence of model discrepancy. Our work is conceptually similar in that we introduce an additive discrepancy function which is calibrated in a Bayesian framework. However, our algorithmic strategy differs considerably in that we focus on general deterministic optimization problems and consider the discrepancy to be the difference in the high and low-fidelity models. Furthermore, we assume that the high-fidelity model may be queried for different values of the optimization variables as specified by the user. This is motivated by optimal design and control problems rather than the model calibration context of Kennedy and O'Hagan.

Fundamentally, our approach builds on post-optimality sensitivity analysis~\cite{shapiro_SIAM_review,fiacco}, which seeks to analyze the influence of auxiliary parameters on the solution of optimization problems. This literature originated in linear programing and was extended for more general applications in nonlinear optimization. Subsequent work from~\cite{Griesse_part_1,Griesse_part_2,griesse2}, and references therein, leveraged post-optimality sensitivity analysis for problems constrained by partial differential equations (PDEs). More recently, HDSA, was developed to enable efficient computation of post-optimality sensitivities in high dimensions by coupling tools from PDE-constrained optimization and numerical linear algebra~\cite{HDSA,saibaba_gsvd,sunseri_hdsa,hart_2021_bayes,sunseri_2,hart_hdsa_control_oed}. These references focus on computing sensitivities with respect to model parameters, particularly in cases where the parameters are infinite dimensional. We leverage this body of literature while building on~\cite{hart_bvw_cmame}, which is the first article to consider the computation of sensitivities with respect to model discrepancy.  It is worth noting that our approach can be viewed as an interpretable learning process, as depicted in Figure 1 from ~\cite{hart_bvw_cmame}. HDSA is used to update optimization solutions by learning from high-fidelity data while providing clear mathematical interpretation of the process.

The remainder of this article is organized as follows. Section~\ref{sec:optimization_formulation} introduces notation and the optimization problem under consideration. We discuss the formulation of post-optimality sensitivity analysis with respect to model discrepancy in Section~\ref{sec:post_opt_sen}. The core mathematical advances to enable posterior optimal solution sampling is detailed in Section~\ref{sec:method_and_derivation}. These advances include a novel formulation of the model discrepancy prior, derivations to enable closed form posterior samples, and the introduction of a posterior optimal solution projection. Section~\ref{sec:computing_samples} provides a guide to computing posterior samples and specification of hyper-parameters used to define the prior covariances. Numerical results demonstrate our advances in Section~\ref{sec:numerical_results}. Section~\ref{sec:conclusion} concludes the article with a discussion of extensions and potential impact of this work.

\section{Optimization Formulation} \label{sec:optimization_formulation}

Consider optimization problems of the form
\begin{align}
\label{eqn:true_opt_prob}
& \min_{z \in \Z} J(S(z),z) 
\end{align}
where $z$ denotes an optimization variable in the Hilbert space $\Z$, $S:\Z \to \U$ is the solution operator for a differential equation with state variable $u$ in an infinite dimensional Hilbert space $\U$, and $J:\U \times \Z \to \R$ is the objective function. Although we would like to solve~\eqref{eqn:true_opt_prob}, it is frequently not possible or practical due to the complexity of the high-fidelity model. Rather, we solve
\begin{align}
\label{eqn:approx_opt_prob}
& \min_{z \in \Z} J(\tilde{S}(z),z) 
\end{align}
where $\tilde{S}:\Z \to \U$ is the solution operator for a low-fidelity model. The goal in this article is to update the optimal solution of~\eqref{eqn:approx_opt_prob} using a limited number of high-fidelity simulations, i.e. to approximate the optimal solution of~\eqref{eqn:true_opt_prob} without requiring evaluation of $S(z)$ within the optimization algorithm. 

The Hilbert space $\U$, and possibly $\Z$, are infinite dimensional. Assume they are discretized and let $\{\eta_i\}_{i=1}^m$ be a basis for a finite dimensional subspace $\U_h \subset \U$ and $\{\varphi_j\}_{j=1}^n$ be a basis for a finite dimensional subspace $\Z_h \subset \Z$. If $\Z$ is finite dimensional then $\Z_h=\Z$, but we consider the general case throughout the article with a focus on problems where $\Z_h$ is high dimensional as a result of $\Z$ being infinite dimensional. Let $\u \in \R^m$ and $\z \in \R^n$ denote coordinates of $u$ and $z$, respectively, in these bases. The remainder of the article focuses on the finite dimensional problem, so for simplicity we abuse notation and use $J:\R^m \times \R^n \to \R$, $S:\R^n \to \R^m$, and $\tilde{S}:\R^n \to \R^m$ to denote the discretized objective function and solution operators. We use $\nabla_\u$ and $\nabla_\z$ to denote the gradient of a scalar-valued function (or the Jacobian of a vector-valued function) with respect to $\u$ and $\z$, respectively. Similarly, $\nabla_{\u,\z}$ and $\nabla_{\z,\z}$ will be used to denote the Hessian with respect to $(\u,\z)$ and $(\z,\z)$, respectively. 

Our approach builds on PDE-constrained optimization. We do not review the tools and algorithms required to solve such problems, but mention a few core components which include adjoint-based derivative calculations, matrix-free linear algebra methods, and iterative solvers. We refer the reader to~\cite{Vogel_99, Archer_01,Haber_01,Vogel_02,Biegler_03,Biros_05,Laird_05,Hintermuller_05,Hazra_06,Biegler_07,Borzi_07,Hinze_09,Biegler_11,frontier_in_pdeco} for a sampling of the relevant literature.

\section{Post-optimality sensitivity with respect to model discrepancy} \label{sec:post_opt_sen}
As introduced in~\cite{hart_bvw_cmame}, consider the parameterized optimization problem
\begin{align}
\label{eqn:dis_approx_opt_prob_pert_rs}
 \min_{\z \in \R^n} \hspace{1 mm} \mathcal{J}(\z,\t):=J(\tilde{S}(\z)+\d(\z,\t),\z) ,
\end{align}
where $\d:\R^n \times \R^p \to \R^m$, which we call the model discrepancy, will approximate the difference between the high and low-fidelity solution operators. That is, we parameterize the discrepancy and seek to determine a parameter vector $\t \in \R^p$ such that $\d(\z,\t) \approx S(\z)-\tilde{S}(\z)$. Assume that $\d$ is parameterized such that $\d(\z,\vec{0})=\vec{0}$ $\forall \z$, i.e. $\mathcal{J}(\z,\vec{0})$ is the low-fidelity optimization objective. Hence we may solve the low-fidelity optimization problem~\eqref{eqn:approx_opt_prob} and then vary $\t$ to analyze how the optimal solution changes with respect to variations in the discrepancy. 

Assume that $\mathcal{J}$ is twice continuously differentiable with respect to $(\z,\t)$ and let $\tilde{\z} \in \R^n$ denote a local minimum of the low-fidelity problem,  i.e. a minimizer of the parameterized problem~\eqref{eqn:dis_approx_opt_prob_pert_rs} when $\t=\vec{0}$. We assume that $\ztilde$ satisfies the first and second order optimality conditions:
\begin{align}
\label{eqn:first_order_opt}
\nabla_\z \mathcal{J}(\ztilde,\vec{0}) = \vec{0}
\end{align}
and that $\nabla_{\z,\z} \mathcal{J}(\ztilde,\vec{0})$ is positive definite. Applying the Implicit Function Theorem to~\eqref{eqn:first_order_opt} gives the existence of an operator $\F: \mathscr N(\vec{0},\epsilon) \to \R^n$, defined on a neighborhood $\mathscr N(\vec{0},\epsilon) = \{ \t \in \R^p \vert \vert \vert \t \vert \vert < \epsilon \}$ for some $\epsilon>0$, such that 
\begin{align*}
\nabla_\z \mathcal{J}(\F(\t),\t) = \vec{0}
\end{align*}
for all $\t \in \mathscr N(\vec{0},\epsilon)$. If the Hessian of $\nabla_{\z,\z} \mathcal{J} (\F(\t),\t)$ is positive definite for all  $\t \in \mathscr N(\vec{0},\epsilon)$, then we interpret $\F(\t)$ as a mapping from the model discrepancy to the solution of the parameterized optimization problem. The Jacobian of $\F$ with respect to $\t$, evaluated at $\t=\vec{0}$, is given by
\begin{eqnarray}
\label{eqn:dis_sen_op}
\nabla_\t \F(\vec{0}) = - \H^{-1} \B \in \R^{n \times p}
\end{eqnarray}
where $\H =\nabla_{\z,\z} \mathcal{J}(\tilde{\z},\vec{0}) \in \R^{n \times n}$ is the Hessian of $\mathcal{J}$ and $\B = \nabla_{\z,\t} \mathcal{J}(\tilde{\z},\vec{0}) \in \R^{n \times p}$ is the Jacobian of $\nabla_\z \mathcal{J}$ with respect to $\t$, both are evaluated at the low-fidelity solution $(\tilde{\z},\vec{0})$. We interpret the post-optimality sensitivity operator $\nabla_\t \F(\vec{0})$ in~\eqref{eqn:dis_sen_op} as the change in the optimal solution as a result of model discrepancy perturbations.

To determine a parameterization of the model discrepancy $\d(\z,\t)$, consider that the post-optimality sensitivity operator~\eqref{eqn:dis_sen_op} depends on two matrices, the inverse Hessian $\H^{-1}$ and the mixed second derivative $\B$. The Hessian $\H$ does not depend on the discrepancy since it is the $(\z,\z)$ derivative evaluated at $\t=\vec{0}$ and $\d(\z,\vec{0}) \equiv \vec{0}$, i.e. $\H$ is the Hessian computed in the low-fidelity optimization. To compute $\B$, we differentiate $\mathcal{J}$ to get
\begin{eqnarray*}
\nabla_{\z,\t} \mathcal{J} = \nabla_\z\tilde{S}^T \nabla_{\u,\u} J \nabla_\t\d + \nabla_\z\d^T \nabla_{\u,\u} J \nabla_\t\d + \nabla_\u J \nabla_{\z,\t} \d + \nabla_{\z,\u} J \nabla_\t\d ,
\end{eqnarray*}
where we adopt the convention that gradients are row vectors. Hence~\eqref{eqn:dis_sen_op} depends on $\nabla_\z \d(\tilde{\z},\vec{0})$, $\nabla_\t \d(\tilde{\z},\vec{0})$, and $\nabla_{\z,\t} \d(\tilde{\z},\vec{0})$. Representing $\d(\z,\t)$ as a Taylor series centered at $\tilde{\z}$, we observe that all nonlinear terms (second order and above) will vanish when evaluated at $\tilde{\z}$. Hence it is sufficient to consider discrepancies which are affine functions of $\z$ parameterized by $\t$. To respect the function space structure and discretization of the problem we consider an arbitrary affine mapping from $\Z_h$ to $\U_h$ and apply the Riesz representation theorem to move into the coordinate space. Then we arrive at the model discrepancy $\d:\R^n \times \R^p \to \R^m$ defined as
\begin{eqnarray}
\label{eqn:delta_kron}
\d(\z,\t) =
\left( \begin{array}{cc}
\I_m & \I_m \otimes \z^T \M_z
\end{array} \right) \t
\end{eqnarray}
where the parameter dimension is $p=m(n+1)$, $\I_m \in \R^{m \times m}$ is the identity matrix, and $\M_\z \in \R^{n \times n}$ is the mass matrix whose $(i,j)$ entry is $(\varphi_i,\varphi_j)_\Z$, where the subscript $\Z$ indicates the inner product defined in the Hilbert space $\Z$. We refer the reader to~\cite{hart_bvw_cmame} for a detailed derivation. If $\Z$ is a function space then $p=m(n+1)$ is very large, for instance, a modest problem in two spatial dimension may yield $m=n=10^4 \implies p>10^8$. The Kronecker product structure of~\eqref{eqn:delta_kron} is critical for the analysis which follows where we derive expressions which do not require computation in $\R^p$. This was enabled by ordering the entries of $\t$ appropriately, as discussed in~\cite{hart_bvw_cmame}.

To arrive at an expression for the post-optimality sensitivity operator~\eqref{eqn:dis_sen_op}, we assume that $\nabla_{\z,\u} J = \vec{0}$\footnote{This is common on a wide range of problems when the objective function admits the structure $J(u,z) = J_{mis}(u) + J_{reg}(z)$, where $J_{mis}$ is a state misfit or design criteria and $J_{reg}$ is an optimization variable regularization.} and evaluate $\nabla_{\z,\t} \mathcal{J} (\ztilde,\vec{0})$ to get
\begin{eqnarray}
\label{eqn:B}
\hspace{10 mm}
\B = 
\nabla_\z \tilde S^T \nabla_{\u, \u} J
\left( \begin{array}{cc}
 \I_m & \I_m \otimes \tilde{\z}^T \M_z 
 \end{array} \right)
 +
\left( \begin{array}{cc}
\vec{0} & \nabla_\u J \otimes \M_z
 \end{array} \right)
 \in \R^{n \times p} 
\end{eqnarray}
where $\nabla_\z \tilde S$ is evaluated at $\tilde{\z}$, and $\nabla_\u J$ and $\nabla_{\u,\u} J$ are evaluated at $(\tilde{S}(\tilde{\z}),\tilde{\z})$.

\section{Updating and quantifying uncertainty in the optimal solution} \label{sec:method_and_derivation}
Assume that the high-fidelity model $S$ is evaluated $N$ times. We seek to use this data to update the low-fidelity optimal solution $\ztilde$. Let $\{ \z_\ell \}_{\ell=1}^N$ denote the inputs for the model evaluations and $\vec{d}_\ell=S(\z_\ell)-\tilde{S}(\z_\ell)$, $\ell=1,2,\dots,N,$ denote the corresponding evaluation of the difference in high and low-fidelity solution operators.

Our proposed approach consist of two steps:
\begin{itemize}

\item[] (calibration) determine $\t \in \R^p$ such that 
$$\d(\z_\ell,\t) \approx S(\z_\ell)-\tilde{S}(\z_\ell), \qquad \ell=1,2,\dots,N,$$

\item[] (optimal solution update) approximate the high-fidelity optimal solution as
$$\ztilde + \nabla_\t \F(\vec{0}) \t.$$

\end{itemize}

We are particularly interested in cases where $N$ is small, for instance, $N<5$, as is common in applications where high-fidelity simulations are computationally intensive. Such limited data implies that the calibration step is ill-posed in the sense that there may be many $\t$'s such that $\d(\z_\ell,\t) \approx S(\z_\ell)-\tilde{S}(\z_\ell)$. Hence we pose a Bayesian inverse problem for the calibration step. This gives a posterior probability distribution $\pi_\t$ which weights vectors $\t \in \R^p$ according to how well $\d(\z_\ell,\t)$ approximates $S(\z_\ell)-\tilde{S}(\z_\ell)$. Then in the optimal solution updating step we compute a posterior distribution for the optimal solution by drawing samples from $\pi_\t$ and propagating them through the post-optimality sensitivity operator $\nabla_\t \F(\vec{0})$. Hence we compute two distributions which we will refer to as the posterior discrepancy and the posterior optimal solution. The character $\vec{\Theta}$ will be used to denote a random variable on $\R^p$ and $\t$ will be used to denote a vector corresponding to a realization of the random variable. 

In what follows, we detail the formulation of the Bayesian inverse problem, provide closed form expressions for posterior samples, and demonstrate how these samples may be efficiently propagated through the post-optimality sensitivity operator. Due to the high dimensionality of the model discrepancy parameterization, $\t \in \R^p$, it is critical that the Kronecker product structure of the discrepancy is preserved in our derivations.  We make judicious choices in the problem formulation and manipulate the linear algebra to achieve this. The approach never requires computation in $\R^p$, and hence is efficient for large-scale optimization problems. We subsequently (in Subsection~\ref{ssec:hess_proj}) introduce a projector $\P$ and the corresponding projected update $\ztilde + \P \nabla_\t \F(\vec{0}) \vec{\Theta}$. Theoretical and computational benefits are demonstrated by the use of this projector.

\subsection{Defining the prior} \label{ssec:prior}
The Bayesian inverse problem to estimate the posterior distribution of $\vec{\Theta}$ consists of two probabilistic modeling choices, a prior distribution and a likelihood. The former plays a critical role in ill-posed problems by encoding expert knowledge to augment limited data, and thus inject more information into the posterior. To facilitate efficient computations we use a Gaussian prior for $\vec{\Theta}$ with mean $\vec{0}$ and a covariance matrix designed to enable computational efficiency and aid interpretability of the posterior.

The probability density function of a mean zero Gaussian random vector with covariance $\vec{\Gamma}^{-1}$ is proportional to $\exp\left( - \frac{1}{2} \vert \vert \mathbf{x} \vert \vert_{\vec{\Gamma}}^2 \right)$, where $\vert \vert \mathbf{x} \vert \vert_{\vec{\Gamma}}^2 = \mathbf{x}^T \vec{\Gamma} \vec{x}$ is the inner product weighted by $\vec{\Gamma}$. Hence the inverse of the covariance, called the precision matrix, defines an inner product such that elements with a small norm will have the largest probability. Accordingly, we seek to define an inner product on the space of model discrepancies and use the corresponding weighting matrix as the prior precision matrix. 

Since the discrepancy is a discretization of an operator mapping between infinite dimensional spaces, it is important to do analysis in the appropriate function spaces to derive the inner product~\cite{stuart_inv_prob}. There are many possible inner products in the space of discrepancies. Due to the dimension of the discrepancy, it is critical that we choose an inner product to achieve computationally advantageous expressions for the discretized prior precision. However, interpretability is also important to aid in specifying hyper-parameters in the prior. Accordingly, we define inner products in the state and optimization variable spaces as building blocks to define an inner product on the model discrepancy parameters through the Kronecker product structure of $\d$.

Let $\W_\u^{-1} \in \R^{m \times m}$ be a covariance matrix on the state space whose precision matrix induces a norm appropriate for measuring elements in the range of the discrepancy. Well established results from infinite dimensional Bayesian inverse problems may be leveraged to guide this modeling choice~\cite{ghattas_infinite_dim_bayes_1,ghattas_infinite_dim_bayes_2,stuart_inv_prob}. In Section~\ref{sec:computing_samples}, we define $\W_\u^{-1}$ as the inverse of a Laplacian-like differential operator raised to an appropriate power to ensure smoothness. This common choice for covariance matrix is equivalent to modeling the discrepancy as a Whittle-Mat\'ern Gaussian random field~\cite{lindgren11,whittle54} and offers computational advantages thanks to the maturity of elliptic solvers. 

The post-optimality sensitivity operator is a linearization about the low-fidelity optimal solution $\ztilde$, so we center our analysis around $\ztilde$  and measure the discrepancy in a neighborhood of the optimal solution. Since the inverse problem is ill-posed, a smoothing prior is needed to constrain the space of possible discrepancies. Analogous to a Sobolev norm which penalizes the magnitude of a function and its derivative, we consider the squared norm of $\d(\cdot,\t)$, corresponding to the inner product of $\d$ with itself, as
\begin{align}
\label{eqn:discrepancy_norm}
\vert \vert \d(\ztilde,\t) \vert \vert_{\W_\u}^2 + \vert \vert \nabla_\z \d(\cdot,\t) \vert \vert^2
\end{align}
where $\vert \vert \nabla_\z \d(\cdot,\t) \vert \vert$ is an appropriately chosen norm on the discrepancy Jacobian. Recalling that $\d(\z,\t)$ depends linearly on $\z$, we interpret~\eqref{eqn:discrepancy_norm} as measuring a nominal value (the discrepancy evaluated at $\ztilde$) plus its variation (the discrepancy Jacobian which does not depend on $\z$). 

Although the norm $\vert \vert \d(\ztilde,\t) \vert \vert_{\W_\u}$ is well defined, there are many possible choices for the operator norm $\vert \vert \nabla_\z \d(\cdot,\t) \vert \vert$. Motivated by a desire for computational efficiency and consistency with a well defined norm in infinite dimensions, we use the Schatten 2-norm
\begin{align*}
\vert \vert \nabla_\z \d(\cdot,\t) \vert \vert^2 = \text{Tr}(\nabla_\z \d(\cdot,\t)^* \nabla_\z \d(\cdot,\t) )
\end{align*}
where $\nabla_\z \d(\cdot,\t)^*$ is the adjoint of $\nabla_\z \d(\cdot,\t)$, and $\text{Tr}$ denotes the trace. This norm requires the definition of inner products on the domain and range of $\nabla_\z \d(\cdot,\t)$. The range is equipped with the $\W_\u$ weighted inner product.

The domain inner product corresponds to weighting directions in the optimization variable space $\mathcal{Z}_h$. Let $\W_\z^{-1} \in \R^{n \times n}$ be a covariance matrix on the optimization variable space which assigns high probability to physically reasonable or realistic $\z$'s. This serves to eliminate discrepancy variations in $\z$ directions that are not relevant (for instance, oscillations which are faster than achievable in a control process). A natural approach is to weight the domain inner product using $\W_\z$, but results in nonphysical $\z$'s being given small weight. Such ignorance of nonphysical directions produces a posterior distribution where $\d(\z,\t)$ varies significantly in nonphysical $\z$ directions, and as a result, produces nonphysical samples from the optimal solution posterior. To avoid such scenarios, the domain inner product is weighted by $\W_\z^{-1}$ so that the posterior discrepancy will have small variation in the nonphysical directions. This interpretation and choice is a significant difference between the prior in~\cite{hart_bvw_cmame} and the prior defined in this article. 

Given our inner product choices, the norm~\eqref{eqn:discrepancy_norm} can be written as
$$\vert \vert \d(\ztilde,\t) \vert \vert_{\W_\u}^2 + \vert \vert \nabla_\z \d(\cdot,\t) \vert \vert^2 = \t^T \W_\t \t,$$ 
where
\begin{eqnarray}
\label{eqn:M_theta}
\W_\t =  \left( \begin{array}{cc}
\W_\u & \W_\u \otimes \ztilde^T \M_z \\
\W_\u \otimes \M_z \ztilde & \W_\u \otimes \M_\z (\W_\z+\ztilde  \ztilde^T) \M_\z 
\end{array} \right)  \in \R^{p \times p}.
\end{eqnarray} 
We define the prior covariance for $\vec{\Theta}$ as $\W_\t^{-1}$ to favor discrepancies that vary in the physically relevant directions defined by $\W_\z^{-1}$ with a range that is weighted according to $\W_\u^{-1}$. Appendix~A provides a detailed derivation of~\eqref{eqn:M_theta} and Section~\ref{sec:computing_samples} provides additional details about the prior covariances.
 
 \subsection{Defining the likelihood}
To facilitate the analysis we assume:
\begin{assumption}
$\{\z_\ell \}_{\ell=1}^N$ is a linearly independent set of vectors in $\R^n$.
\end{assumption}
That is, we require that the inputs to the high-fidelity model evaluations are linearly independent. This is usually satisfied in practice as one explores the input space in different directions since the number of model evaluations $N$ is small relative to the dimension of $\z$. For notational simplicity, we define
\begin{eqnarray}
\label{eqn:Aell}
\A_\ell = \begin{pmatrix} \I_{m} & \I_{m} \otimes \z_\ell^T \M_z \end{pmatrix} \in \R^{m \times p}, \qquad \ell=1,2,\dots,N,
\end{eqnarray}
so that $\d(\z_\ell,\t)=\A_\ell \t$. By concatenating these matrices
\begin{eqnarray*}
\A = \left(
\begin{array}{c}
\A_1 \\
\A_2 \\
\vdots \\
\A_N
\end{array}
\right) \in \R^{mN \times p} ,
\end{eqnarray*}
$\A \t \in \R^{mN}$ corresponds to the evaluation of $\d(\z,\t)$ for all of the input data. In an analogous fashion, let $\vec{d} \in \R^{mN}$ be defined by stacking $\vec{d}_\ell \in \R^m$, $\ell=1,2,\dots,N$, into a vector. Then we seek $\t$ such that $\A \t \approx \vec{d}$.

To enable a closed form expression for the posterior we consider an additive Gaussian noise model with mean $\vec{0}$ and covariance $\alpha_{\vec{d}} \M_\u^{-1}$, where $\M_\u \in \R^{m \times m}$ is the state space mass matrix whose $(i,j)$ entry  is $(\phi_i,\phi_j)_\U$ and $\alpha_{\vec{d}}>0$ scales the noise covariance. Defining the noise covariance as a scalar multiple of the inverse mass matrix ensures that we achieve mesh independence by computing data misfit norms with the function space norm from $\mathcal U$. The noise covariance has a different interpretation when calibrating to high-fidelity simulation data. This is discussed in Subsection~\ref{ssec:hyper-parameter_selection}.

Thanks to our formulation with Gaussian prior and noise models, and the linearity of $\d(\z,\t)$ as a function of $\z$, the posterior is Gaussian with mean and covariance~\cite{gentle_intro_bayes}
\begin{eqnarray}
\label{eqn:post_mean_product}
\hspace{8 mm} \overline{\t} = \alpha_{\vec{d}}^{-1}  \vec{\Sigma} \A^T (\I_N \otimes \M_\u) \vec{d} \quad \text{and} \quad \vec{\Sigma} = \left( \W_\t + \alpha_{\vec{d}}^{-1} \A^T (\I_N \otimes \M_\u) \A \right)^{-1}, \hspace{-6 mm}
\end{eqnarray}
where $\I_N$ denotes the $N\times N$ identity matrix. In principle, other (non-Gaussian) prior and noise models may be considered. However, care must be taken to define them properly in function spaces. Furthermore, the derivations below leverage the form of the posterior mean and covariance~\eqref{eqn:post_mean_product} to provide computationally efficient algorithms. Analogous derivations are possible with other prior and noise models, but extensive effort may be required. We focus on the Gaussian case for this article. 

%

\subsection{Expressions for posterior sampling}\label{subsec:posterior_sampling_expressions}
Our modeling choices enable closed-form expressions for the posterior mean and covariance~\eqref{eqn:post_mean_product}. Nonetheless, these expressions are not sufficient to enable efficient computation of posterior samples due to the dimension $p$. Rather, leveraging the Kronecker product structure of the model discrepancy discretization and prior covariance matrix, we seek to express the posterior quantities in terms of Kronecker products of solves and factorizations in the state and optimization variable spaces. In what follows, we introduce linear algebra constructs needed to state the resulting posterior sample expressions. These results are stated as Theorems and their proofs are given in the Appendix. 

Let $\vec{Z} = \begin{pmatrix} \z_1 & \z_2 & \dots & \z_N \end{pmatrix} \in R^{n \times N}$ be the matrix of optimization variable data, i.e. the inputs to the high-fidelity model evaluations, and define 
\begin{align*}
& \G= \e \e^T + (\vec{Z} - \ztilde \e^T)^T \W_\z^{-1} (\vec{Z}-\ztilde \e^T) \in \R^{N \times N},
\end{align*}
where $\e \in \R^N$ is the vector of ones. Let $(\vec{g}_i,\mu_i)$ denote the eigenvectors and eigenvalues of $\vec{G}$ and define
\begin{eqnarray}
\vec{y}_i = \vec{Z} \vec{g}_i - (\e^T \vec{g}_i) \ztilde \qquad \text{and} \qquad s_i = (\e^T \vec{g}_i) - \vec{y}_i^T \W_\z^{-1} \ztilde 
\end{eqnarray}
for $i=1,2,\dots,N$. These quantities arise naturally from the Appendix proofs which decomposes the $p=m(n+1)$ dimensional posterior into a $mN$ dimensional data informed component and its $m(n-N+1)$ dimensional complement. All of these quantities may be computed at a negligible cost when $N$ is small.

Theorem~\ref{thm:post_mean} expresses the posterior mean in a form amenable for computation. 
\begin{theorem}
\label{thm:post_mean}
The posterior mean for the model discrepancy parameters $\vec{\Theta}$ is 
\begin{eqnarray*}
\overline{\t} = \frac{1}{\alpha_{\vec{d}}} \sum\limits_{\ell=1}^N 
\left[
\left(
\begin{array}{c}
a_\ell \u_\ell \\
 \u_\ell  \otimes  \M_\z^{-1} \W_\z^{-1} (\z_\ell-\ztilde) 
 \end{array}
\right)
-
\sum\limits_{i=1}^N b_{i,\ell}
\left(
\begin{array}{c}
s_i \vec{u}_{i,\ell} \\
 \vec{u}_{i,\ell} \otimes   \M_\z^{-1} \W_\z^{-1} \vec{y}_i 
\end{array}
\right)
\right]
\end{eqnarray*}
where
\begin{align*}
& a_{\ell} = 1-\ztilde^T \W_\z^{-1} (\z_\ell-\ztilde) \qquad  & b_{i,\ell} = (\z_\ell-\ztilde)^T \W_\z^{-1} \vec{Z} \vec{g}_i + (\e^T \vec{g}_i) a_\ell \\
& \vec{u}_\ell = \W_\u^{-1} \M_\u \vec{d}_\ell & \vec{u}_{i,\ell} = \left( \alpha_{\vec{d}} \W_\u + \mu_i \M_\u \right)^{-1} \M_\u \u_\ell .
\end{align*}
\end{theorem}
It follows from Theorem~\ref{thm:post_mean} that the posterior mean may be computed with a cost of $\mathcal O(N)$ linear solves involving priors precision operators on the state and optimization variables spaces.

To facilitate an expression for posterior samples, we make the additional assumption.
 \begin{assumption}
 $\z_1=\ztilde$.
 \end{assumption}
 That is, we require that data is collected at the nominal optimal solution $\ztilde$; a reasonable assumption given that our goal is to use high-fidelity forward solves to improve it. To arrive at an expression for posterior samples, we introduce $\{ \breve{\z}_1,\breve{\z}_2,\dots,\breve{\z}_{n-N+1} \} \subset \R^n$, which we assume to be an orthonormal set of vectors which are orthogonal to $\{ \W_\z^{-\frac{1}{2}} (\z_\ell-\ztilde) \}_{\ell=2}^N$. As demonstrated later in the article, we will never need to explicitly compute these vectors. However, their existence and orthogonality properties play an important role in the analysis needed to express posterior samples in a convenient form. This set of vectors arises to model uncertainty in directions which are not informed by the data $\{\z_\ell\}_{\ell=1}^N$.
 
Samples from the posterior distribution are given by Theorem~\ref{thm:post_samples}.
\begin{theorem}
\label{thm:post_samples}
Posterior samples of $\vec{\Theta}$ take the form
\begin{eqnarray*}
\overline{\t} + \hat{\t} + \breve{\t} 
\end{eqnarray*}
where $\overline{\t}$ is given in Theorem~\ref{thm:post_mean}, 
\begin{eqnarray*}
\hat{\t} =  \sqrt{\alpha_{\vec{d}}} \sum\limits_{i=1}^N \frac{1}{\sqrt{\mu_i}} 
\left(
\begin{array}{c}
s_i \hat{\u}_i \\
\hat{\u}_i \otimes \M_\z^{-1} \W_\z^{-1} \vec{y}_i
\end{array}
\right)
\qquad
\text{and} 
\qquad
\breve{\t} = 
\sum\limits_{k=1}^{n-N+1}
\left(
\begin{array}{c}
\breve{s}_k \breve{\u}_k\\
\breve{\u}_k \otimes \breve{\vec{y}}_k
\end{array}
\right),
\end{eqnarray*}
with
\begin{align*}
& \hat{\u}_i \sim \mathcal N(\vec{0},\left( \alpha_{\vec{d}} \W_\u + \mu_i \M_\u \right)^{-1} ) \qquad & \breve{s}_k= -\ztilde^T \W_\z^{-\frac{1}{2}} \breve{\z}_k \\
& \breve{\u}_k \sim \mathcal N(\vec{0}, \W_\u^{-1} )  & \breve{\vec{y}}_k = \M_\z^{-1} \W_\z^{-\frac{1}{2}} \breve{\z}_k .
\end{align*}
\end{theorem}
We use the notation $\sim \mathcal N(\vec{0},\vec{\Gamma}^{-1})$ to denote a Gaussian random vector with mean $\vec{0}$ and covariance $\vec{\Gamma}^{-1}$. Theorem~\ref{thm:post_samples} provides a computationally efficient expression for posterior samples of the parameters $\vec{\Theta}$. However, $\vec{\Theta}$ does not have a physical interpretation. Rather, our focus is on the model discrepancy $\d$ which it parameterizes, and the effect of the model discrepancy on the solution of the optimization problem. The subsections which follow propagate the posterior samples from Theorem~\ref{thm:post_samples} through the discrepancy and post-optimality sensitivity operator.

\subsection{Posterior model discrepancy}
Substituting $\overline{\t}$, $\hat{\t}$, and $\breve{\t}$ from Theorems~\ref{thm:post_mean} and~\ref{thm:post_samples} into~\eqref{eqn:delta_kron}, we have a decomposition of the posterior model discrepancy
\begin{eqnarray}
\label{eqn:delta_post}
\d(\z, \overline{\t} + \hat{\t} + \breve{\t}) = \overline{\d}(\z) + \hat{\d}(\z) + \breve{\d}(\z)
\end{eqnarray}
where
\begin{align}
\label{eqn:delta_bar}
 \overline{\d}(\z) &=  \frac{1}{\alpha_{\vec{d}}} \sum\limits_{\ell=1}^N (1+(\z_\ell-\ztilde)^T \W_\z^{-1} (\z-\ztilde)) \u_\ell  \\
 & -\frac{1}{\alpha_{\vec{d}}} \sum\limits_{\ell=1}^N \sum\limits_{i=1}^N b_{i,\ell} ( \vec{e}^T \vec{g}_i + \vec{y}_i^T \W_\z^{-1} (\z-\ztilde)) \u_{i,\ell}  \nonumber 
\end{align}
is the mean model discrepancy defined by linear combinations of the discrepancy data preconditioned by the prior,
\begin{align}
\label{eqn:delta_hat}
  \hat{\d}(\z) &= \sqrt{\alpha_{\vec{d}}} \sum\limits_{i=1}^N \frac{1}{\sqrt{\mu_i}} ( \e^T \vec{g}_i + \vec{y}_i^T \W_\z^{-1} (\z-\ztilde)) \hat{\u}_i 
\end{align}
represents discrepancy uncertainty in the directions defined by the data $\{\z_\ell\}_{\ell=1}^N$, and
\begin{align}
\label{eqn:delta_breve}
   \breve{\d}(\z) &= \sum\limits_{k=1}^{n-N+1} \left( \breve{\z}_k^T \W_\z^{-\frac{1}{2}} (\z-\ztilde) \right) \breve{\u}_k
\end{align}
represents discrepancy uncertainty in directions informed exclusively by the prior.  $\overline{\d}: \R^n \to \R^m$ is a deterministic function whereas $\hat{\d}$ and $\breve{\d}$ are random functions mapping from $\R^n$ to $\R^m$. We have written~\eqref{eqn:delta_hat} and~\eqref{eqn:delta_breve} as deterministic functions of the independent random vectors $\{ \hat{\u}_i \}_{i=1}^N$ and $\{ \breve{\u}_k \}_{k=1}^{n-N+1}$, which are distributed as $\mathcal N(\vec{0},\left( \alpha_{\vec{d}} \W_\u + \mu_i \M_\u \right)^{-1} )$ and $\mathcal N(\vec{0},\W_\u^{-1})$, respectively. 

This decomposition into data informed and prior informed directions is apparent when we evaluate $\overline{\d}$, $\hat{\d}$, and $\breve{\d}$ for $\z$'s in the subspaces defined by the data $\{\z_\ell\}_{\ell=1}^N$ and its orthogonal complement. Observe that $\breve{\d}$ vanishes in data informed directions, $\breve{\d}(\z_\ell)=0$, $\ell=1,2,\dots,N$, since $\{\breve{\z}_k\}_{k=1}^{n-N+1}$ are orthogonal to $\{\W_\z^{-\frac{1}{2}} (\z_\ell-\ztilde)\}_{\ell=2}^N$. Furthermore, observe that 
\begin{eqnarray*}
\overline{\d}(\ztilde + \W_\z^{\frac{1}{2}} \breve{\z}_k)=\overline{\d}(\ztilde) \qquad \text{and} \qquad \hat{\d}(\ztilde + \W_\z^{\frac{1}{2}} \breve{\z}_k)=\hat{\d}(\ztilde), \qquad k=1,2,\dots,n-N+1,
\end{eqnarray*}
do not vary in directions orthogonal to the data. Additionally, observe that $\hat{\d}(\z)  \to 0$ as $\alpha_{\vec{d}} \to 0$. In other words, making the noise variance arbitrarily small implies that there will be no uncertainty in the data informed directions. On the other hand, $\breve{\d}(\z)$ has no dependence on $\alpha_{\vec{d}}$ since it models uncertainty in directions orthogonal to the data. This decomposition of the discrepancy is unsurprising since the Bayesian inverse problem is linear with Gaussian noise and prior models. Nonetheless, we emphasize the benefit of the linear algebra manipulations in the Appendix which enables these interpretable and computationally advantageous expressions which would otherwise be intractable due to high dimensionality. 

The mean $\overline{\d}$ and samples of the data informed direction uncertainty $\hat{\d}$ can be computed with $\mathcal O(N)$ linear solves and samples in $\R^m$ and $\R^n$. However, computing a realization of $\breve{\d}$ is more difficult since it involves the basis vectors $\{ \breve{\z}_k \}_{k=1}^{n-N+1}$ which are impractical to compute or store when $n$ is large. This challenge may be avoided by introducing the orthogonal projector onto the low dimensional subspace $\{\W_\z^{-\frac{1}{2}} (\z_\ell-\ztilde)\}_{\ell=2}^N$ and then representing $\{ \breve{\z}_k \}_{k=1}^{n-N+1}$ using the projector's orthogonal complement. Theorem~\ref{thm:delta_breve_samples} provides an efficient expression to compute samples.

\begin{theorem}
\label{thm:delta_breve_samples}
A sample of $\breve{\d}(\z)$ may be computed as $\gamma(\z) \vec{\nu}_\u$ where $\vec{\nu}_\u \sim \mathcal N(\vec{0},\W_\u^{-1})$ and 
\begin{align*}
\gamma(\z) = \sqrt{ (\z-\ztilde)^T \left( \W_\z^{-1} - \W_\z^{-1} \vec{Z}_c \left(  \vec{Z}_c^T \W_\z^{-1}  \vec{Z}_c \right)^{-1}  \vec{Z}_c^T \W_\z^{-1} \right) (\z-\ztilde) }
\end{align*}
where $ \vec{Z}_c = \begin{pmatrix}  \z_2-\ztilde & \cdots & \z_N - \ztilde \end{pmatrix} \in \R^{n \times (N-1)}$.
\end{theorem}
 Realizations of $\breve{\d}(\z)$ take the form of a sample from the mean zero state prior with covariance $\W_\u^{-1}$, multiplied by a $\z$ dependent coefficient which equals $0$ when evaluated at $\z=\ztilde$ and is scaled by the optimization variable prior covariance $\W_\z^{-1}$. This highlights the importance of $\W_\z^{-1}$ which determines the posterior model discrepancy rate of change with respect to $\z$ in the uninformed directions. 

\subsection{Posterior optimal solution}
To characterize uncertainty in the solution of the optimization problem, we propagate posterior samples $\overline{\t}+\hat{\t}+\breve{\t}$ through the post-optimality sensitivity operator $\nabla_\t \F(\vec{0})$. Given the Kronecker product structure, we compute the action of $\B$ on $\overline{\t}$, $\hat{\t}$, and $\breve{\t}$ separately, and then apply $-\H^{-1}$ to the sum of the resulting vectors. Observe that
\begin{align}
\label{eqn:B_thetabar}
\B \overline{\t} &=   \frac{1}{\alpha_{\vec{d}}} \nabla_\z \tilde{S}^T \nabla_{\u, \u} J \left[ \sum\limits_{\ell=1}^N \left(  \u_\ell-  \sum\limits_{i=1}^N b_{i,\ell} (\e^T \vec{g}_i) \vec{u}_{i,\ell} \right)\right] \\
& + \frac{1}{\alpha_{\vec{d}}} \sum\limits_{\ell=1}^N ( \nabla_\u J \u_\ell) \W_\z^{-1} ( \z_\ell-\ztilde) \nonumber \\
&- \frac{1}{\alpha_{\vec{d}}} \sum\limits_{\ell=1}^N \sum\limits_{i=1}^N b_{i,\ell} (\nabla_\u J \vec{u}_{i,\ell}) \W_\z^{-1} \vec{y}_i \nonumber
\end{align}
\begin{align}
\label{eqn:B_thetahat}
\B \hat{\t} & =   \sqrt{\alpha_{\vec{d}}}  \nabla_\z \tilde{S}^T \nabla_{\u, \u} J \left( \sum\limits_{i=1}^N \frac{\e^T \vec{g}_i}{\sqrt{\mu_i}} \hat{\u}_i \right) +  \sqrt{\alpha_{\vec{d}}}   \sum\limits_{i=1}^N \frac{\nabla_\u J \hat{\u}_i}{\sqrt{\mu_i}} \W_\z^{-1} \vec{y}_i 
\end{align}
and
\begin{align}
\label{eqn:B_thetatilde}
\B \breve{\t}= \sum\limits_{k=1}^{n-N+1} \left( \nabla_\u J \breve{\u}_k \right) \W_\z^{-\frac{1}{2}} \breve{\z}_k .
\end{align}
Hence, we compute samples of the optimal solution as $\ztilde - \H^{-1} \B (\overline{\t} + \hat{\t} + \breve{\t})$ without requiring computation in $\R^p$. Equations \eqref{eqn:B_thetabar}-\eqref{eqn:B_thetatilde} provides a systematic and interpretable combination of high-fidelity data in $\vec{u}_\ell$ and $\vec{u}_{i,\ell}$, the low-fidelity model in $ \nabla_\z \tilde{\vec{S}}^T$, prior information via $\W_\u^{-1}$ and $\W_\z^{-1}$, and the optimization objective in $\nabla_\u J$ and $\nabla_{\u,\u} J$. 

Analogous to the posterior discrepancy, we compute $\B \overline{\t}$ and samples of $\B \hat{\t}$ efficiently with $\mathcal O(N)$ linear solves and samples. However, computing samples of $\B \breve{\t}$ involves a sum over the computationally prohibitive basis vectors $\{ \breve{\z}_k \}_{k=1}^{n-N+1}$. Analogous to Theorem~\ref{thm:delta_breve_samples}, we leverage an orthogonal projector to compute samples of $\B \breve{\t}$ which only requires $\mathcal O(N)$ linear solves and samples. Theorem~\ref{thm:B_breve_samples} gives this efficient expression.
\begin{theorem}
\label{thm:B_breve_samples}
A sample of $\B \breve{\t}$ can be computed as $\Gamma(\z) \vec{\nu}_\z$ where $\vec{\nu}_\z \sim \mathcal N(\vec{0},\W_\z^{-1})$ and 
\begin{eqnarray}
\label{eqn:B_theta_tilde}
\Gamma(\z) = \sqrt{\nabla_\u J \W_\u^{-1} \nabla_\u J^T}  \left( \I_n - \W_\z^{-1} \vec{Z}_c  \left( \vec{Z}_c^T \W_\z^{-1} \vec{Z}_c \right)^{-1} \vec{Z}_c^T \right).
\end{eqnarray}
\end{theorem}
We stress the importance of our prior covariances. Samples of $\B \breve{\t}$ involve samples from a Gaussian with covariance $\W_\z^{-1}$, which defines characteristic variances and length scales of the optimization variable.

\subsection{Hessian Projection} \label{ssec:hess_proj}


The ultimate goal is to learn an update of the optimal solution by applying the post-optimality sensitivity operator $\nabla_\t F(\vec{0})=-\H^{-1} \B$ in directions defined by the high-fidelity data, and characterize uncertainty in this update due to the data sparsity (having a limited number of high-fidelity simulations). In many applications, there is a low dimensional subspace of the optimization variable space for which the objective function exhibits sensitivity, while it has small variations in directions orthogonal to this subspace. As a result, small errors estimating the model discrepancy may cause large changes in the optimal solution that have a negligible effect on the objective function. In such cases, it is advantageous to project the posterior optimal solution onto the subspace of greatest objective function sensitivity. This will reduce the complexity of the posterior optimal solution by implicitly regularizing the effect of the model discrepancy calibration on it. Additionally, projecting the posterior optimal solution can reduce the computational cost to generate samples by exposing a low dimensional subspace on which the samples may be computed.

To this end, observe that the optimization objective function can be approximated in a neighborhood of the minimizer by its second order Taylor polynomial, 
\begin{align*}
\mathcal{J}(\z,\vec{0}) & \approx \mathcal{J}(\ztilde,\vec{0}) + \nabla_{\z} \mathcal{J}(\ztilde,\vec{0}) (\z-\ztilde) + \frac{1}{2} (\z-\ztilde)^T \H (\z-\ztilde) \\
& = \mathcal{J}(\ztilde,\vec{0}) + \frac{1}{2} (\z-\ztilde)^T \H (\z-\ztilde) .
\end{align*}
The post-optimality sensitivity update consists of perturbing the gradient by varying the model discrepancy, and then taking a Newton like step, applying $-\H^{-1}$ to the perturbed gradient, to update the optimal solution. The greatest sensitivities occur in directions of the leading eigenvalues of $\H^{-1}$, or equivalently the smallest eigenvalues of $\H$. These are directions for which the objective function has less curvature, i.e. modifying the optimal solution in those directions will yield the least benefit in terms of reducing the value of the objective function. Hence, we compute the projection of the sensitivity operator onto a subspace for which the objective function has large curvature, i.e. the span of the leading eigenvectors of $\H$. 

Samples from the projected posterior optimal solution, which we will refer to as the posterior optimal solution for simplicity whenever appropriate, take the form $\ztilde - \P \H^{-1} \B (\overline{\t} + \hat{\t} + \breve{\t})$, where $\P$ is a projector onto the subspace of objective function sensitivity. To define this projector, we compute the leading eigenvectors of $\H$ in the $\W_\z$ weighted inner product. This inner product choice ensures mesh independence and encourages directions in the optimization variable space which are consistent with the prior information encoded in $\W_\z$. Assume that 
\begin{align}
\label{eqn:hess_gen_eig}
\H \vec{v}_j = \rho_j \W_\z \vec{v}_j, \qquad j=1,2,\dots,r,
\end{align}
where $r$ is the truncation rank (the eigenpairs are ordered with $\rho_1 \ge \rho_2 \ge \cdots$), and let $\vec{V} = \begin{pmatrix} \vec{v}_1 & \vec{v}_2 & \cdots \vec{v}_r \end{pmatrix} \in \R^{m \times r}$. We define the projector as $\P = \vec{V} \vec{V}^T \W_\z$. $\vec{V}$ is an orthogonal matrix\footnote{This is guaranteed by the Spectral theorem if the eigenvalues are distinct. If the eigenvalues are repeated we can orthonormalize the eigenvectors to define $\vec{V}$.} in the $\W_\z$ weighted inner product since $\H$ is symmetric positive definite. Then we observe that for any $\vec{x} \in \R^m$,
\begin{align}
\label{eqn:proj_hess_inv}
\P \H^{-1} \vec{x} = \sum\limits_{j=1}^r \frac{1}{\rho_j} \left( \vec{v}_j^T \vec{x} \right) \vec{v}_j.
\end{align}

Computing samples $\ztilde - \P \H^{-1} \B (\overline{\t} + \hat{\t} + \breve{\t})$ has the theoretical advantage that it avoids consideration of large variations which yield negligible improvements in the objective function. Furthermore, there is a computational advantage when computing projected samples. Computing $N_\text{post}$ posterior optimal solution samples without projection requires $N_\text{post}$ inverse Hessian-vector products. Since the Hessian is generally only available through matrix-vector products, which involve two linear PDE solves (the incremental state and adjoint equations), the action of the inverse Hessian must be computed through an iterative linear solver. Hence the cost to compute $N_\text{post}$ inverse Hessian-vector products is $2 N_\text{post} L_\text{iter}$ PDE solves, where $L_\text{iter}$ is the average number of iterations required by the linear solver. In contrast, computing the projection of the inverse Hessian via~\eqref{eqn:proj_hess_inv} requires the leading $r$ eigenvalues and eigenvectors of the Hessian, which can be precomputed with $\mathcal O(r)$ PDE solves.

\section{Computing samples} \label{sec:computing_samples}
The previous section provided expressions to compute samples which require $\mathcal O(N)$ linear solves which involve the prior covariances $\W_\u^{-1}$ and $\W_\z^{-1}$. In this section we propose to use Laplacian-like differential operators to define these prior covariances. This is motivated by the prevalence of this prior covariance choice in the infinite dimensional Bayesian inverse problems~\cite{stuart_inv_prob} and Gaussian random field~\cite{lindgren11} literature, and the interpretability of the hyper-parameters which define these covariances. Furthermore, we propose to leverage a generalized Singular Value Decomposition (GSVD) of the inverse Laplacian-like operator which enables efficient computation of the samples $\hat{\u}_i \sim \mathcal N(\vec{0},\left( \alpha_{\vec{d}} \W_\u + \mu_i \M_\u \right)^{-1} )$ required to sample the discrepancy and optimal solution posteriors.

\subsection{Laplacian-like operators}
\label{subsec:laplace_like_operators}
Let $\K_{\u}$ and $\K_\z$ denote the stiffness matrices\footnote{The stiffness matrix $(i,j)$ entry is the inner product of the $i^{th}$ and $j^{th}$ basis function's gradients.} in $\mathcal U_h$ and $\mathcal Z_h$, respectively.  We define Laplacian-like operators $\E_\u = \beta_\u \K_\u + \M_\u$ and $\E_\z = \beta_\z \K_\z + \M_\z$, where $\beta_\u,\beta_\z \ge 0$ are constants used to control smoothness. We have assumed zero Neumann boundary conditions, although other boundary conditions are used as appropriate~\cite{stadler_cov_bc,khristenko19}. Define the precision matrices as 
\begin{align*}
\W_\u = \frac{1}{\alpha_\u} \E_\u^T \M_\u^{-1} \E_\u \qquad \text{and} \qquad \W_\z = \frac{1}{\alpha_\z} \E_\z^T \M_\z^{-1} \E_\z,
\end{align*}
where $\alpha_\u, \alpha_\z>0$ are constants which scale the variance in the prior\footnote{This definition of the prior covariances ensures that they correspond to the discretization of trace class covariance operators as long as the domains of $u$ and $z$ has dimension $ \le 3$, which covers many problems of practical relevance since these domains are typically spatial coordinates.}. 

The computation of samples requires factorizations involving the matrices $\W_\u^{-1}$ and $(\alpha_{\vec{d}} \W_\u + \mu_i \M_\u)^{-1}$, of which the latter is particularly challenging. However, factorizations can be efficiently approximated using a GSVD of $\E_\u^{-1}$ with inner products weighted by $\M_\u^{-1}$ and $\M_\u$. Specifically, the GSVD is given by 
\begin{align*}
\E_\u^{-1} =  \V \vec{\Pi} \vec{\daleth}^T \M_\u^{-1}
\end{align*}
where $\V^T \M_\u \V = \I_m$ are the $\M_\u$ orthonormal left singular vectors, $\vec{\Pi}$ is the diagonal matrix of singular values, and $\vec{\daleth}^T \M_\u^{-1} \vec{\daleth}=\I_m$ are the $\M_\u^{-1}$ orthonormal right singular vectors. It follows that
\begin{align}
\label{eqn:W_u_inv}
\W_\u = \frac{1}{\alpha_\u} \M_\u \V \vec{\Pi}^{-2} \V^T \M_\u \qquad \text{and} \qquad \W_\u^{-1} = \alpha_\u \V \vec{\Pi}^{2} \V^T .
\end{align}
The Sherman-Morrison-Woodbury formula gives 
\begin{align}
\label{eqn:W_u_shift_inv}
(\alpha_{\vec{d}} \W_\u+\mu_i \M_\u)^{-1} = \alpha_\u \V \vec{\aleph} \V^T
\end{align}
where $\vec{\aleph}$ is a diagonal matrix whose diagonal entries are given by 
\begin{align*}
\vec{\aleph}_{j,j} = \frac{ \pi_j^2}{\alpha_{\vec{d}}+\alpha_\u \mu_i \pi_j^2},
\end{align*}
where $\pi_j$ is the $j^{th}$ diagonal entry of $\vec{\Pi}$.

We generally take $\beta_\u$ sufficiently large so that $\E_\u^{-1}$ imposes smoothness on the prior discrepancy. Hence the GSVD is truncated at a rank $q$. Then we can approximate $\W_\u^{-1}$ and $(\alpha_{\vec{d}} \W_\u+\mu_i \M_\u)^{-1}$ via rank $q$ truncations of~\eqref{eqn:W_u_inv} and~\eqref{eqn:W_u_shift_inv}, respectively, with a truncation error that is $\mathcal O(\pi_{q+1}^2)$. Similarly, approximate sample from $\mathcal N(\vec{0},\W_\u^{-1})$ and $\mathcal N(\vec{0},(\alpha_{\vec{d}} \W_\u+\mu_i \M_\u)^{-1})$ can be computed by factorizing the rank $q$ truncations of~\eqref{eqn:W_u_inv} and~\eqref{eqn:W_u_shift_inv}, respectively, with a truncation error that is $\mathcal O(\pi_{q+1})$. Since the relative truncation error in samples is $\mathcal O( \frac{\pi_{q+1}}{\pi_1})$, and truncation underestimating uncertainty (particularly for smaller correlation lengths), we seek a rank $q$ such that $\mathcal O( \frac{\pi_{q+1}}{\pi_1}) << 1$. The truncation rank may also be assessed by plotting the singular vectors to observe their correlation lengths.

Algorithm~\ref{alg:post_discrepancy_samples} summarizes the computation of projected posterior optimal solution samples using the rank $q$ prior truncation and rank $r$ projector of the optimal solution updates. The right column of Algorithm~\ref{alg:post_discrepancy_samples} describes the dominant cost at each step, which is typically large linear system solves. Samples from $\mathcal N(\vec{0},\W_\u^{-1})$ and $\mathcal N(\vec{0},(\alpha_{\vec{d}} \W_\u+\mu_i \M_\u)^{-1})$ are not included in the computational cost analysis since they are enabled by the GSVD of $\E_\u^{-1}$. A detailed computational cost analysis is given in Appendix~C.
\begin{algorithm}[ht!!]
	\caption{Projected Optimal Solution Posterior Sampling}
	\begin{algorithmic} [1] \label{alg:post_discrepancy_samples}
		\STATE Input: Number of samples $s$, state prior rank $q$, projector rank $r$
		\STATE Compute the truncated GSVD of $\E_\u^{-1}$  \hspace{30 mm} \COMMENT{\textcolor{red}{$\mathcal O(q)$ $\E_\u^{-1}$ multiplies}} 
		\STATE Compute the posterior mean $\overline{\t}$  \hspace{44 mm} \COMMENT{\textcolor{red}{$N$ $\W_\z^{-1}$ multiplies}}
		\STATE Compute $s$ samples $\hat{\t}$ and $\breve{\t}$ \hspace{44 mm} \textcolor{red}{$s$ $\mathcal N(\vec{0},\W_\z^{-1})$ samples}
		\STATE Compute $\B (\overline{\t} + \hat{\t} + \breve{\t})$  \hspace{54 mm} \textcolor{red}{$s+1$ $\nabla_\z S^T$ multiplies}
		\STATE Compute the truncated GEVD of $\H$ \hspace{20 mm} \textcolor{red}{$\mathcal O(r)$ $\H$ and $\W_\z^{-1}$ multiplies}
		\STATE Compute $\ztilde + \P \H^{-1} \B(\overline{\t} + \hat{\t} + \breve{\t})$
	\end{algorithmic}
\end{algorithm}

\subsection{Guidance for hyper-parameter specification} \label{ssec:hyper-parameter_selection}
The proposed approach has five hyper-parameters defining the prior and noise covariances, $\alpha_\u,\beta_\u,\alpha_z,\beta_\z,\alpha_{\vec{d}}$. The interpretation and influence of the hyper-parameters on the posterior can be understood thanks to the closed form expressions derived for sampling. We discuss the interpretation of the hyper-parameters and provide general insight to guide their specification. 

\subsubsection*{State variance and correlation length}
The state variance $\alpha_\u$ and correlation length parameter $\beta_\u$, which define $\W_\u^{-1}$, determine the magnitude and smoothness of samples from the state prior $\mathcal N(\vec{0},\W_\u^{-1})$. One can show that for samples $\t \sim \mathcal N(\vec{0},\W_\t^{-1})$, the prior discrepancy evaluated at $\z=\ztilde$ follows the state prior distribution, i.e. $\d(\ztilde,\t) \sim \mathcal N(\vec{0},\W_\u^{-1})$. Furthermore, recalling Theorem~\ref{thm:delta_breve_samples}, the posterior discrepancy in data uniformed directions, $\breve{\d}(\z)$, can be sampled as $\vec{\nu}_\u \sim \mathcal N(\vec{0},\W_\u^{-1})$ scaled by a scalar coefficient $\gamma(\z)$. Hence we must specify $\alpha_\u$ and $\beta_\u$ so that samples from $\mathcal N(\vec{0},\W_\u^{-1})$ have magnitudes and smoothness characteristics commensurate to $S-\tilde{S}$.

\subsubsection*{Optimization variable correlation length}
Assuming that $\z$ corresponds to the discretization of a function \footnote{If $\z$ does not correspond to the discretization of a function, then $\W_\z$ would be specified differently and would not involve a correlation length hyper-parameter.}, the hyper-parameter $\beta_\z$ defines the smoothness of the optimal solution perturbations. In particular, Theorem~\ref{thm:B_breve_samples} indicates that samples of $\B \breve{\t}$ are proportional to samples from $\mathcal N(\vec{0},\W_\z^{-1})$ and hence the smoothness of such samples is determined by $\beta_\z$. The hyper-parameter $\beta_\z$ should be chosen to reflect a correlation length that is commensurate to the optimization variable. 

\subsubsection*{Optimization variable variance}
The optimization variable variance $\alpha_\z$ determines the magnitude of samples $\mathcal N(\vec{0},\W_\z^{-1})$ and, as a result of Theorem~\ref{thm:B_breve_samples}, the magnitude of uncertainty in the optimal solution posterior. However, unlike $\alpha_\u$ which could be chosen based on characteristics of $S-\tilde{S}$, the specification of $\alpha_\z$ is less clear. Since post-optimality sensitivity analysis is local about the low-fidelity optimal solution $\ztilde$, $\alpha_\z$ serves to restrict the analysis within a reasonable neighborhood of $\ztilde$. Rather than specifying $\alpha_\z$ based on characteristics of the optimization variable or the optimal solution posterior, we focus on how $\alpha_\z$ affects the model discrepancy.

Recall again from Theorem~\ref{thm:delta_breve_samples} that samples of the posterior discrepancy in data uniformed directions are scaled by $\gamma(\z)$. Furthermore, for variations of $\z$ in a $n-N+1$ dimensional subspace we have
\begin{align*}
\gamma(\z) = \sqrt{ (\z-\ztilde)^T  \W_\z^{-1} (\z-\ztilde) } = \sqrt{\alpha_\z} \vert \vert \E_\z^{-1} (\z-\ztilde) \vert \vert_{\M_\z}.
\end{align*}
Hence the magnitude of the posterior model discrepancy uncertainty in data uniformed directions is directly proportional to $\sqrt{\alpha_\z}$. Using the posterior to guide specification of the prior hyper-parameters is, in general, questionable. However, since $\breve{\d}(\z,\t)$ is only informed by the input data $\{\z_\ell\}_{\ell=1}^N$ and not the high-fidelity model evaluations, it is justifiable to analyze it when determining prior hyper-parameters. A general approach is to choose $\z \notin \{\z_\ell \}_{\ell=1}^N$ whose magnitude and smoothness characteristics are representative of a possible optimal solution, and evaluate $\breve{\d}(\z)$ to assess uncertainty in the model discrepancy posterior in the uninformed direction defined by $\z$. The hyper-parameter $\alpha_\z$ is chosen so that the magnitude of uncertainty in $\breve{\d}(\z)$ is commensurate to the magnitude of $S-\tilde{S}$ that may be expected if additional high fidelity evaluations are available. The user may also consider a range of values for $\alpha_\z$ and analyze how the optimal solution posterior changes under these varying scenarios of prior assumptions on the magnitude of $S-\tilde{S}$.

\subsubsection*{Noise variance}
 In many Bayesian inverse problems, the noise variance $\alpha_{\vec{d}}$ is determined based on knowledge of the data collection process and the noise it generates. In our context the data is not noisy (assuming we work with high-fidelity simulation data). Nonetheless, we do not necessarily want $\d(\z,\t)$ to fit the data to high precision as $\d$ is a linear approximation of a nonlinear operator. Taking large values for $\alpha_{\vec{d}}$ will result in a poor fit to the data and hence a failure to utilize the high-fidelity simulations. Taking small values for $\alpha_{\vec{d}}$ may encourage overfitting and error due to the linear approximation. Furthermore, it follows from~\eqref{eqn:delta_hat} that the posterior model discrepancy uncertainty is proportional to $\sqrt{\alpha_{\vec{d}}}$. A guiding principle is to tune $\alpha_{\vec{d}}$ so that the discrepancy posterior evaluated at $\{\z_\ell\}_{\ell=1}^N$ fits the data and has a reasonable magnitude of uncertainty. 

\section{Numerical results} \label{sec:numerical_results}
This section presents three examples to highlight various types of model discrepancy and aspects of the proposed optimal solution updating approach. The first example, a boundary value problem in one spatial dimension, illustrates our approach for a model discrepancy that is due to using a spatially homogenous low-fidelity model when the high-fidelity model has spatial heterogeneity. The second example, a system of ordinary differential equations modeling a mass-spring system, demonstrates how our optimal solution updating approach can utilize a low-fidelity model that neglects model coupling, to estimate the optimal solution for the high-fidelity coupled system. The third example, a steady state PDE in two spatial dimensions, shows how a low-fidelity model based on linearization can be augmented with our approach to improve the optimal solution using only one evaluation of the high-fidelity nonlinear PDE model.
 
\subsection{Diffusion reaction example}
\label{subsec:illustrative_example}
Consider the boundary value diffusion reaction equation 
\begin{align*}
& -\kappa \Delta u + R(u) = z \qquad & \text{on } (0,1) \\
& \kappa u' = 0 & \text{on } \{0,1\}
\end{align*}
where $\kappa>0$ is the diffusion coefficient and $R(u)$ is the reaction function. We consider two models for the reaction function with a low-fidelity model corresponding to $R(u) = u^2$ and its high-fidelity counterpart having $R(u) = (1+.7 \sin(2 \pi x)) u^2$. We seek the design of a source injection to achieve a target state, i.e.
\begin{align}
\label{eqn:illustrative_example_opt_prob}
\min_{z} J(\tilde{S}(z),z) \coloneqq \frac{1}{2} \int_0^1 (\tilde{S}(z)(x) - T(x))^2 dx + \frac{\gamma}{2} \int_0^1 z(x)^2 dx
\end{align}
where $z:[0,1] \to \R$ is the optimization variable, $\tilde{S}(z)$ is the solution operator for the low-fidelity model, $T(x) = 20(x+.5)(1.3-x)$ is the target state, and $\gamma=10^{-4}$ is the regularization coefficient. 

\begin{figure}[h]
\centering
  \includegraphics[width=0.49\textwidth]{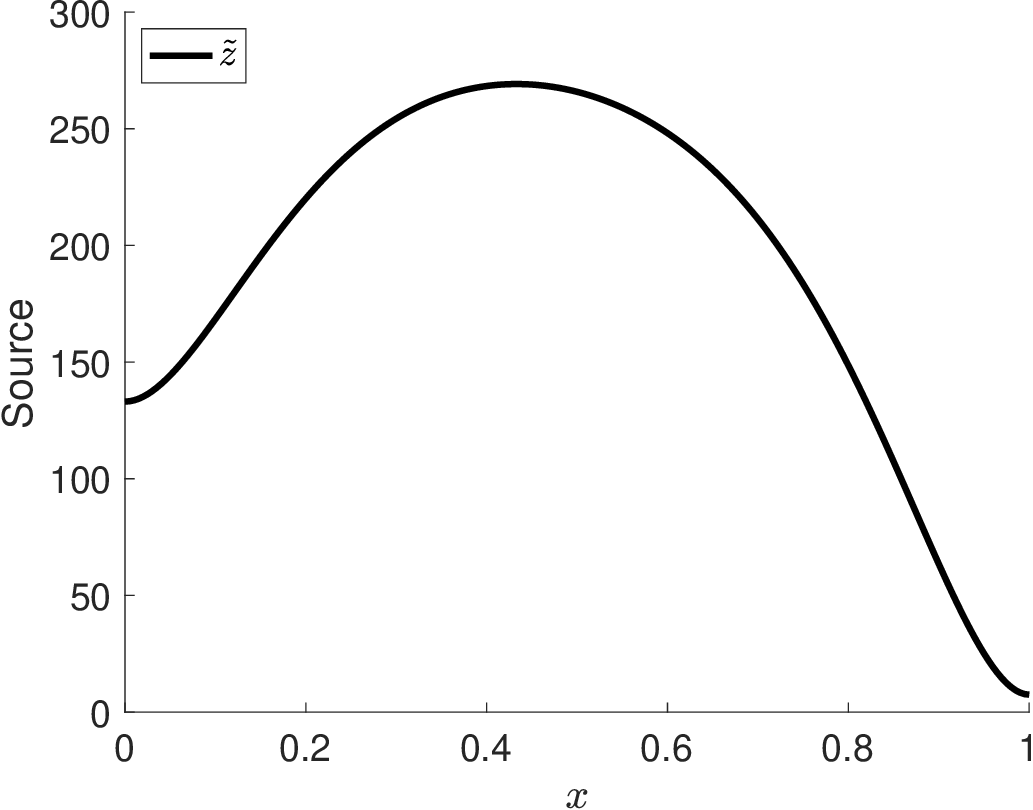}
    \includegraphics[width=0.49\textwidth]{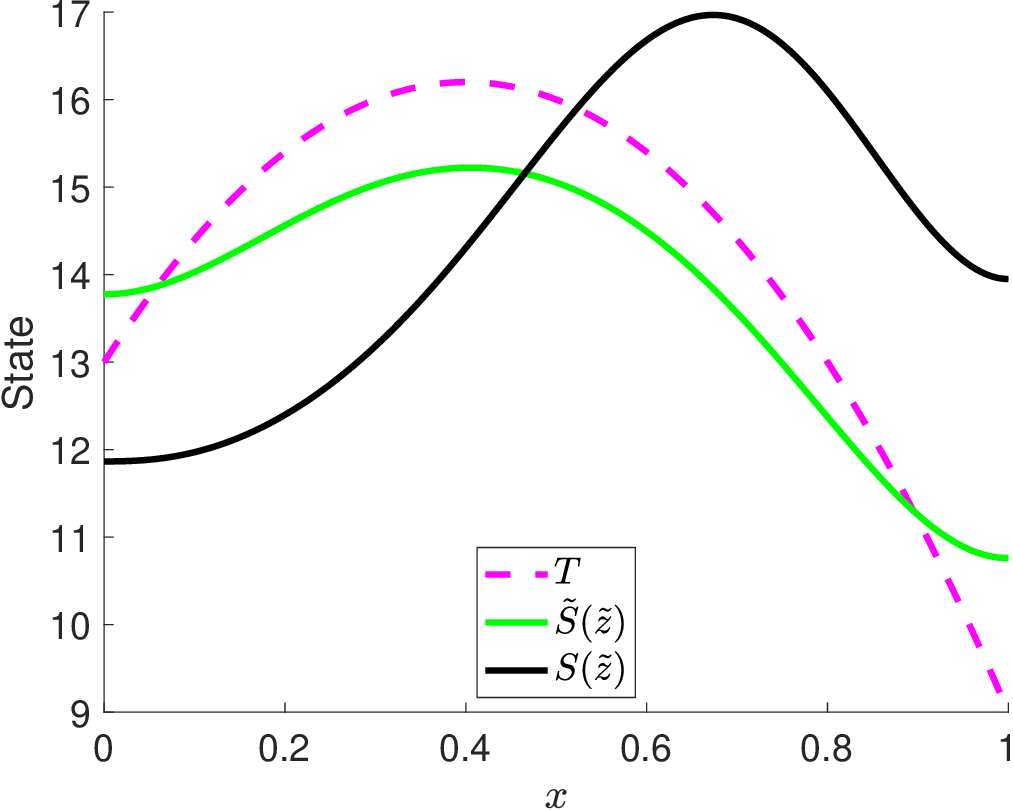}
    \caption{Low-fidelity optimization solution for the reaction diffusion example. Left: optimal source $\tilde{\z}$; right: state target $T$ alongside the low and high-fidelity state solutions evaluated at $\tilde{\z}$, i.e. $\tilde{S}(\tilde{\z})$ and $S(\tilde{\z})$, respectively.}
  \label{fig:illustrative_example_lofi_solution}
\end{figure}

The left panel of Figure~\ref{fig:illustrative_example_lofi_solution} displays the solution of~\eqref{eqn:illustrative_example_opt_prob}, i.e. the optimal source term $\tilde{\z}$. The right panel of Figure~\ref{fig:illustrative_example_lofi_solution} displays the low and high-fidelity solution operators $\tilde{S}(\tilde{\z})$ and $S(\tilde{\z})$ evaluated at the low-fidelity optimal solution, as well as the target state $T$. As expected, $\tilde{S}(\tilde{\z})$ is similar to the target $T$ since $\tilde{\z}$ was designed to minimize their difference, but the high-fidelity model $S(\tilde{\z})$ is significantly different. This highlights the limitation of optimization using a low-fidelity model. Our goal is to use $N$ evaluations of the high-fidelity model to update $\tilde{\z}$, for this example we restrict ourselves to $N=2$ evaluations, i.e. $S(\z_1)$ and $S(\z_2)$, where $\z_1=\tilde{\z}$ and $\z_2$ is randomly sampled with a correlation length commensurate to $\tilde{\z}$. 

We follow the guidance of Subsection~\ref{ssec:hyper-parameter_selection} to determine the prior hyper-parameters and report their values in Table~\ref{tab:diff_react_hyper-parameters}. A detailed description of the hyper-parameter selection process is given in Appendix~D.

\begin{figure}[h]
\centering
  \includegraphics[width=0.49\textwidth]{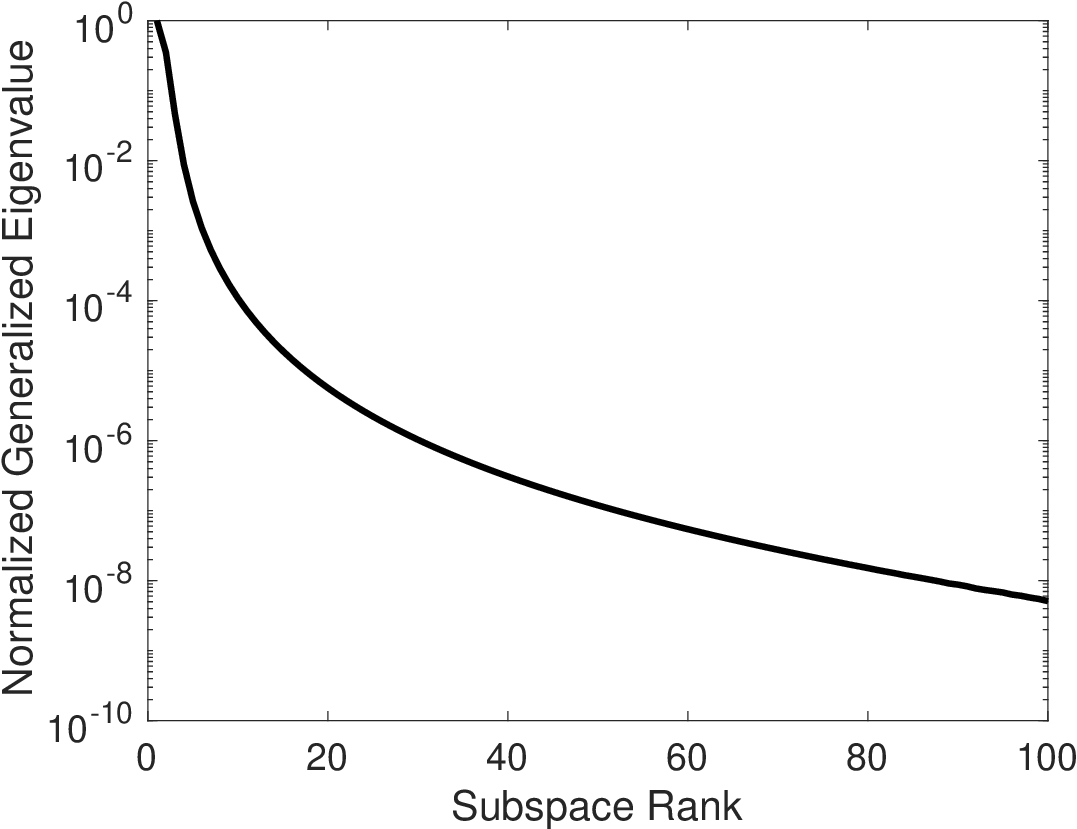}
    \includegraphics[width=0.49\textwidth]{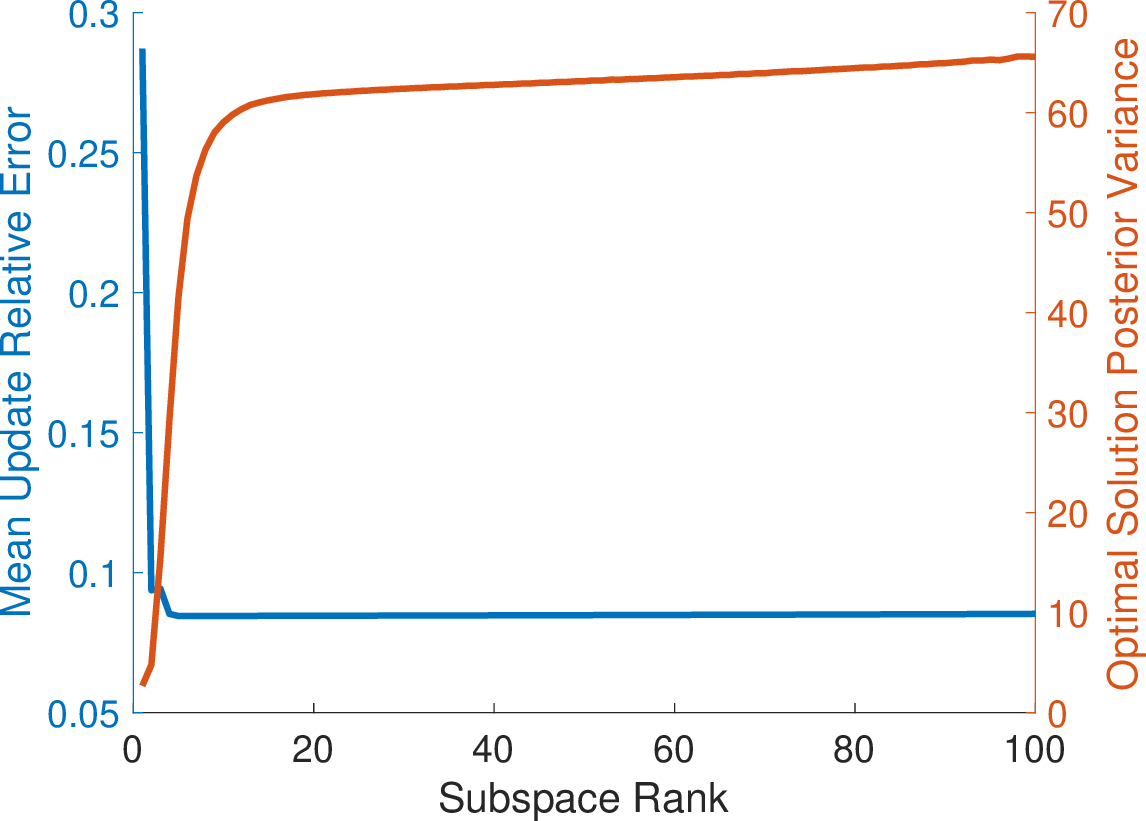}
    \caption{Left: generalized eigenvalues~\eqref{eqn:hess_gen_eig} for the reaction diffusion example; right: the posterior optimal solution mean relative error (in the left axis) and posterior optimal solution variance (in right axis) for various subspace ranks.}
  \label{fig:illustrative_example_vary_rank}
\end{figure}

To analyze the effect of the projector, we consider the posterior optimal solution for various subspace ranks ranging from $r=1,2,\dots,100$. The left panel of Figure~\ref{fig:illustrative_example_vary_rank} displays that the generalized eigenvalues decay rapidly. For comparison, we compute the optimal solution using the high-fidelity model, which is generally not practical to compute but is used here to verify the quality of the posterior optimal solution. For each of the potential ranks, $r=1,2,\dots,100$, we compute the posterior optimal solution projected on the subspace defined by the leading $r$ generalized eigenvectors and compare it with the high-fidelity optimal solution. In particular, we compute the relative error of the posterior mean approximation to the high-fidelity optimal solution. These errors are shown in the right panel of Figure~\ref{fig:illustrative_example_vary_rank}. To understand how the posterior variance is affected by the subspace rank, the right panel of Figure~\ref{fig:illustrative_example_vary_rank} also shows the integrated (over the interval $[0,1]$) variance of the posterior optimal solution. We observe that, in this example, the minimum error in the posterior optimal solution mean is achieved with rank $r=5$. A sharp decrease in the mean error and increase in the variance is observed for small ranks, but both quantities exhibit small changes after a modest rank.

\begin{figure}[h]
\centering
  \includegraphics[width=0.49\textwidth]{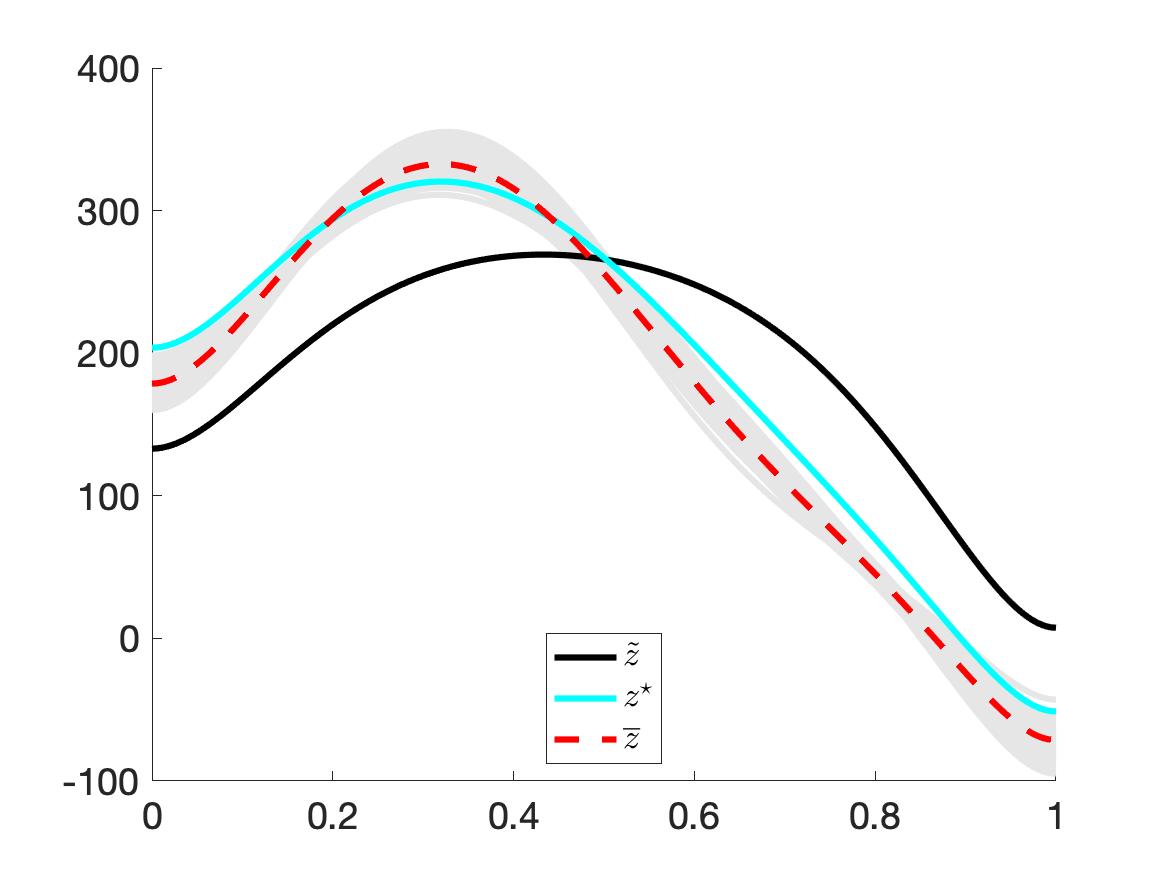}
    \includegraphics[width=0.49\textwidth]{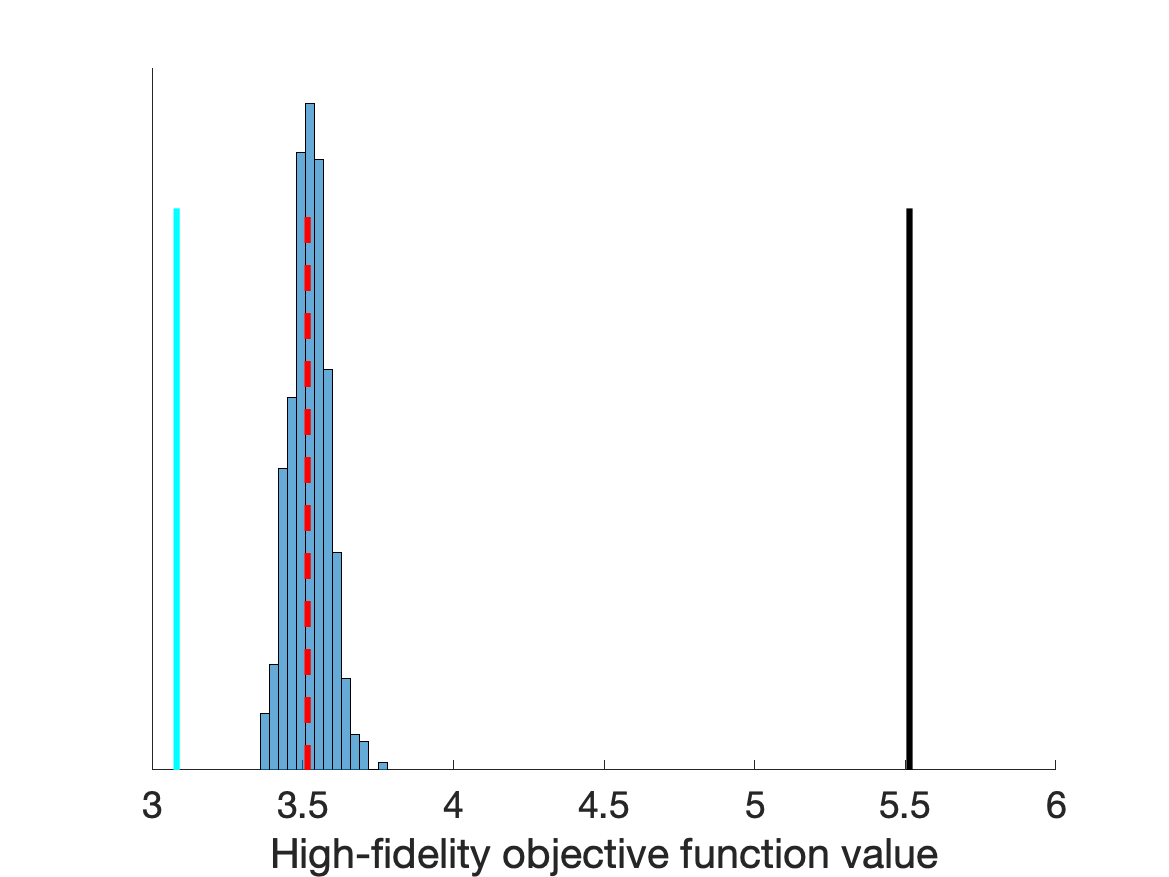} \\
      \includegraphics[width=0.49\textwidth]{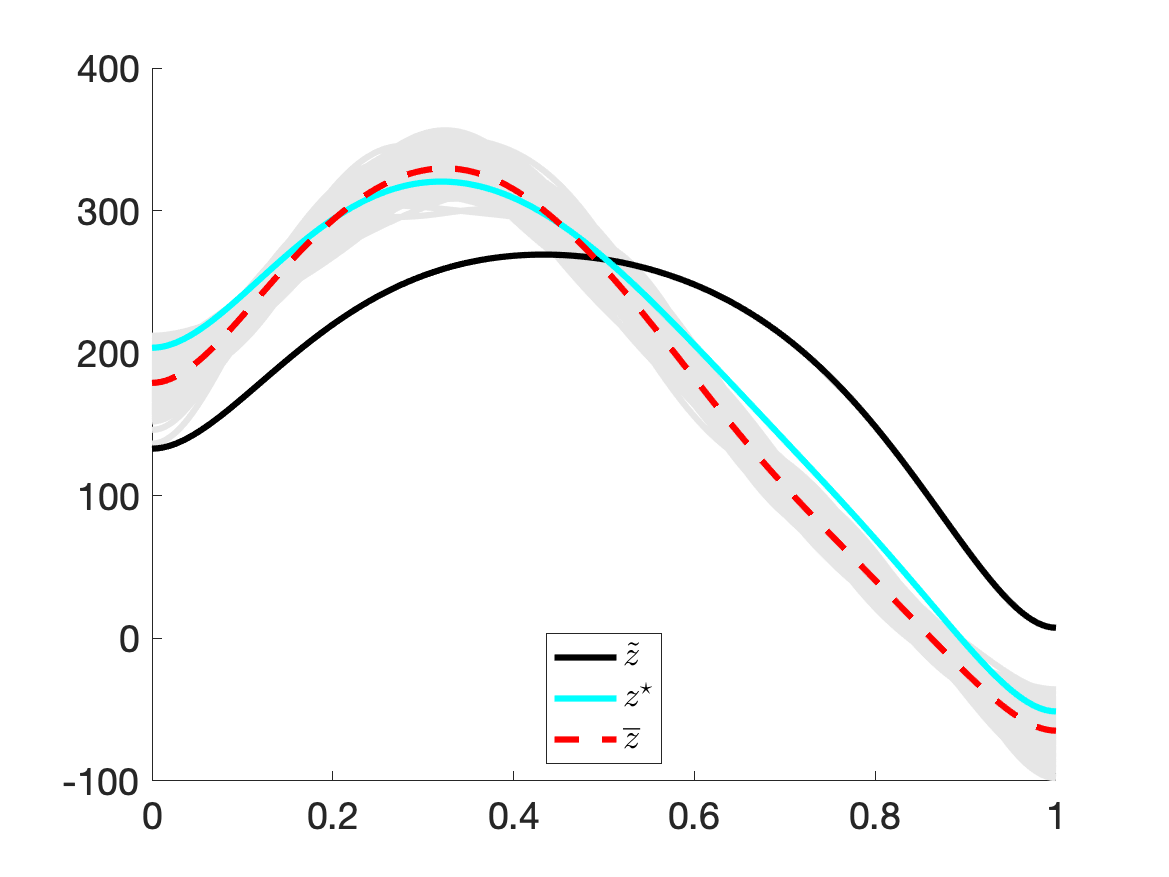}
    \includegraphics[width=0.49\textwidth]{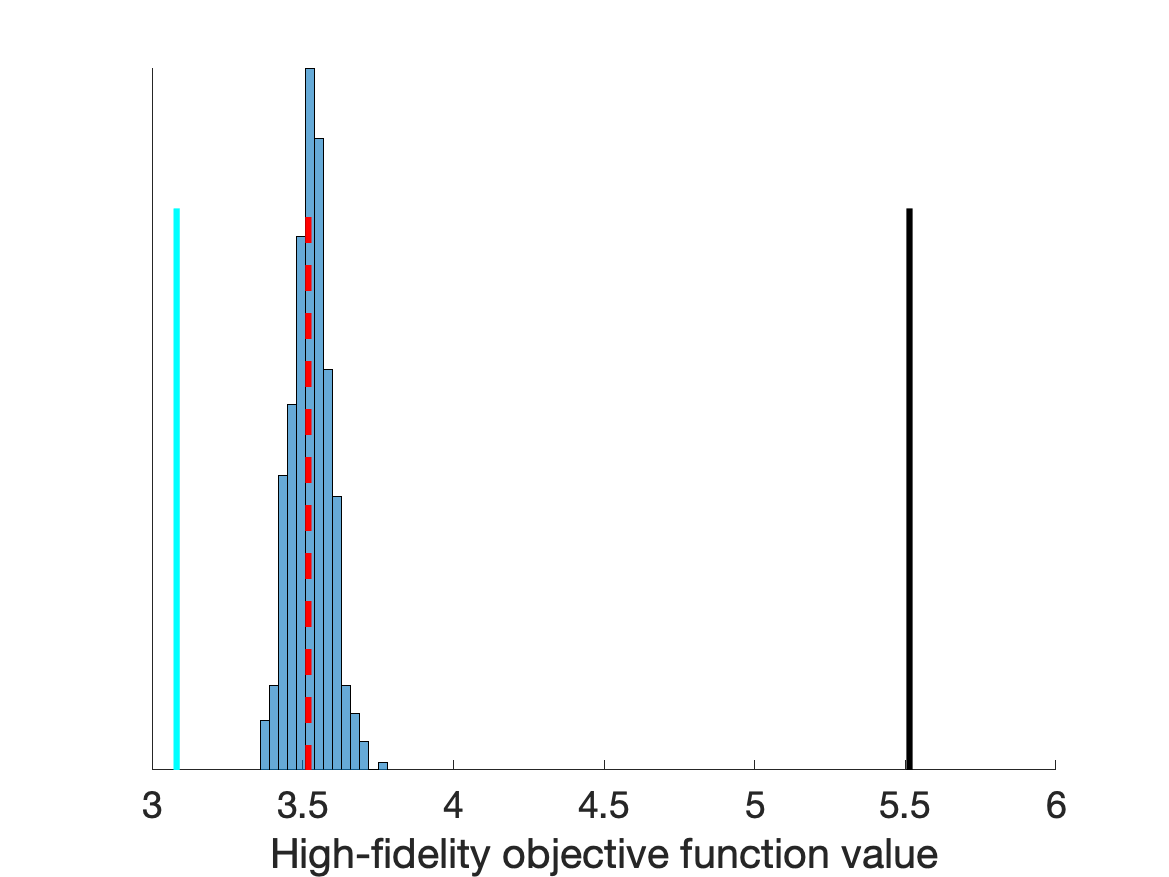}
    \caption{Left: the low and high-fidelity optimal sources $\tilde{\z}$ and $\z^\star$, respectively, alongside the posterior optimal solution mean $\overline{\z}$ with posterior optimal solution samples (in grey); right: the high-fidelity objective function values for 500 optimal source posterior samples, the vertical lines from left to right indicates the value the high-fidelity objective $J(S(\z),\z)$ evaluated at the high-fidelity optimal source, the posterior optimal source mean, and the low-fidelity optimal source. The top and bottom rows correspond to a subspace projection ranks of 4 and 11, respectively.}
  \label{fig:illustrative_example_opt_solution_posterior}
\end{figure}

To illustrate the posterior optimal solution we consider two cases: a rank $r=4$ subspace and a rank $r=11$ subspace, which are displayed in the top and bottom rows of Figure~\ref{fig:illustrative_example_opt_solution_posterior}, respectively. For each, we compute $500$ posterior optimal solution samples. The left panels of Figure~\ref{fig:illustrative_example_opt_solution_posterior} display the low-fidelity optimal source, the optimal source posterior (mean and samples), and the optimal source corresponding to solving the optimization problem~\eqref{eqn:illustrative_example_opt_prob} using the high-fidelity model $S(\z)$ in place of the low-fidelity model $\tilde{S}(\z)$. The right panels of Figure~\ref{fig:illustrative_example_opt_solution_posterior} display the high-fidelity objective function values. The histograms are generated using the 500 optimal source posterior samples, i.e. the data $\{ J(S(\z^k),\z^k) \}_{k=1}^{500} \subset \R$, where $\z^k$ are the posterior optimal solution samples. The leftmost vertical line indicates the value of the minimum of the high-fidelity objective $J(S(\z),\z)$, the rightmost vertical line indicates the high-fidelity objective evaluated at the low-fidelity optimal source, i.e. $J(S(\tilde{\z}),\ztilde)$, and the middle vertical line defines the high-fidelity objective evaluated at the posterior optimal solution mean. This demonstrative figure, which is impractical to compute if high-fidelity evaluations are limited, shows the improvement in the high-fidelity objective function value.

Comparing the rank $r=4$ and $r=11$ cases, we do not observe a noticeable difference in the posterior optimal solution means. Both posterior means provide a measurable improvement relative to the low-fidelity optimal solution. Similarly, comparing the histograms of high-fidelity objective function values are very similar highlighting how changes due to the larger rank yields small changes in the objective function value. This is unsurprising given the rapid eigenvalue decay. The noticeable difference when comparing the top and bottom rows is the spread of posterior samples. Increasing the rank translates to having a greater variance in the posterior optimal solution which is apparent from comparing the posterior samples in the left panel top and bottom rows.

\begin{table}[!ht]
\centering
\begin{tabular}{|c|c|c|c|c|}
\hline
$\alpha_\u$ & $\beta_\u$ & $\alpha_\z$ & $\beta_\z$ & $\alpha_{\vec{d}}$  \\
$4$ & $2 \times 10^{-2}$ & $10^{-10}$ & $3 \times 10^{-2}$ & $10^{-4}$  \\
\hline
\end{tabular}
\caption{Hyper-parameters used in the diffusion-reaction example.}
\label{tab:diff_react_hyper-parameters}
\end{table}

\subsection{Mass-spring system example}
\begin{figure}[h]
\centering
  \includegraphics[width=0.7\textwidth]{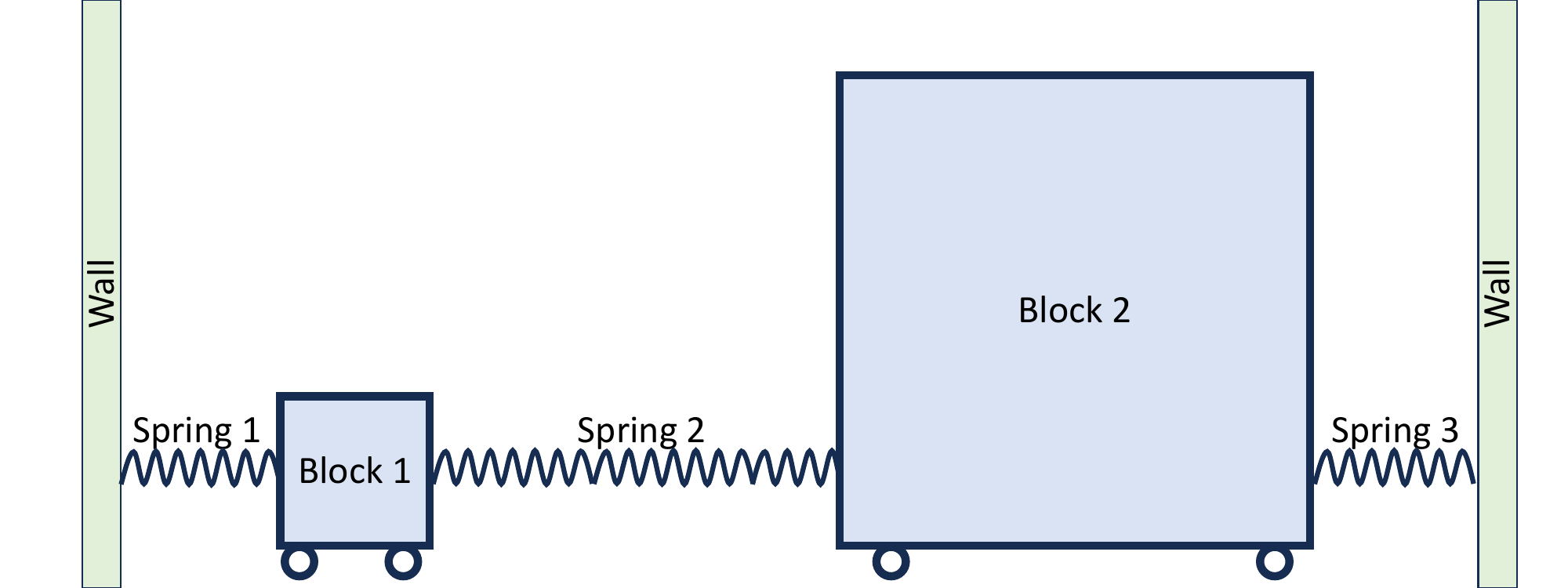}
    \caption{Diagram depicting the mass-spring system. We consider the coupled system for both blocks 1 and 2 as the high-fidelity model. A low-fidelity model is derived by assuming that block 2 is stationary.}
  \label{fig:mass_spring_diagram}
\end{figure}

Consider a simple example of the common modeling practice where a low-fidelity model is derived by neglecting the variability of other processes in the physical system. Specifically, a mass-spring system, depicted in Figure~\ref{fig:mass_spring_diagram}, with two masses connected to one another by a spring, and connected by springs to walls on each end of the domain. Under idealized conditions, the high-fidelity model is 
\begin{align}
\label{eqn:mass_spring_coupled}
& m_1 x_1''(t) = k_2 x_2(t) - (k_1+k_2)x_1(t)+ z(t) \\
& m_2 x_2''(t) = k_2 x_1(t) - (k_2+k_3)x_2(t) \nonumber
\end{align} 
where $x_1(t),x_2(t)$ is the displacements of the two blocks, $m_1,m_2$ is defined as their masses, $k_1, k_2,k_3$ denote the spring constants for each spring, and $z(t)$ is an external forcing acting to move the first block. We assume that the blocks are initially at rest with $x_1(0)=x_2(0)=0$, that $k_1=k_2=k_3=1$, and that $m_1=1$ and $m_2=10$. We consider the optimal control problem
\begin{align}
\label{eqn:mass_spring_opt_prob}
\min_{z} J(S(z),z) \coloneqq \frac{1}{2} \int_0^{10} (S_1(z)(t) - T(t))^2 dt + \frac{\gamma}{2} \int_0^{10} z(t)^2 dt
\end{align}
where $S_1(z)(t)=x_1(t)$ is the solution operator for the displacement of block 1, i.e. evaluating $S_1(z)(t)$ requires solving~\eqref{eqn:mass_spring_coupled} for a given $z$ and extracting the state $x_1$, $z:[0,10] \to \R$ is the controller which seeks to move block 1 along the trajectory $T(t) = 5t^2$, and $\gamma=10^{-6}$ is the regularization coefficient. 

A low-fidelity model is derived by assuming that block 2 is stationary, i.e. $x_2(t)=0$ for all $t$, and modeling block 1 via the simpler differential equation
\begin{align}
\label{eqn:mass_spring_simplified}
m_1 x_1''(t) =  - (k_1+k_2)x_1(t)+ z(t) .
\end{align} 
In this scenario, $\tilde{S}_1(z)$ is the solution operator for~\eqref{eqn:mass_spring_simplified} which will be compared with the high-fidelity model solution operator $S_1(z)$. We use the subscript to emphasize that the solution operator for the high-fidelity model takes values in the product of the function spaces containing $x_1(t)$ and $x_2(t)$, whereas the low-fidelity model takes values in only the function space containing $x_1(t)$. We compare the solution operators for their predictions of $x_1(t)$, and its time derivative $x_1'(t)$, a state variable when the second order system is rewritten as a first order system. The high-fidelity prediction of $x_2(t)$ is not used in our analysis.

\begin{figure}[h]
\centering
  \includegraphics[width=0.49\textwidth]{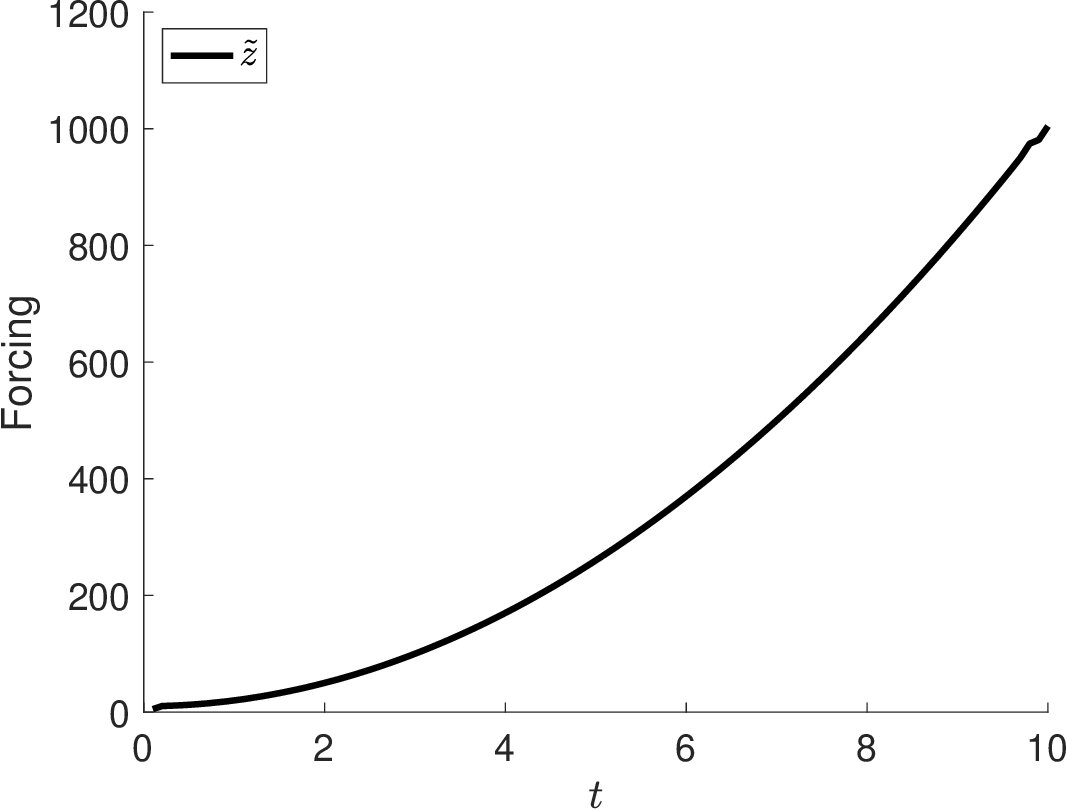}
    \includegraphics[width=0.49\textwidth]{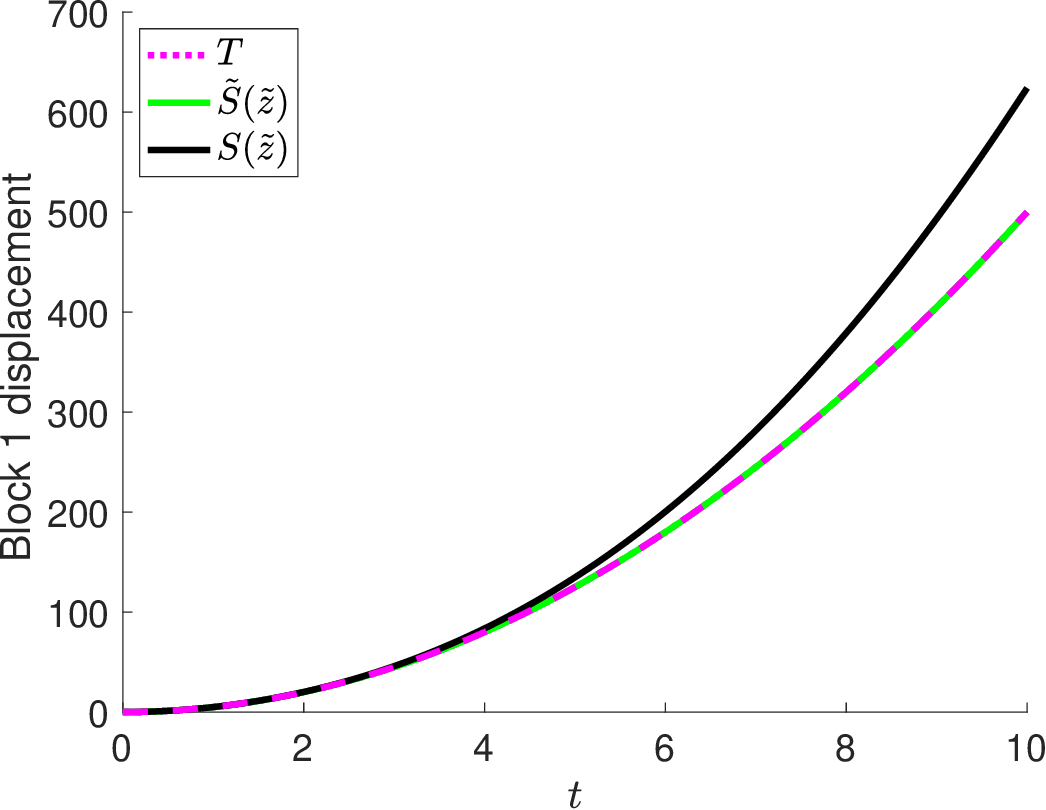}
    \caption{Low-fidelity optimization solution for the mass-spring system. Left: low-fidelity optimal forcing $\tilde{\z}$; right: state trajectory $T$ alongside the low and high-fidelity block 1 displacement solutions evaluated at $\tilde{\z}$, i.e. $\tilde{S}(\tilde{\z})$ and $S(\tilde{\z})$, respectively.}
  \label{fig:mass_spring_lofi_solution}
\end{figure}

Figure~\ref{fig:mass_spring_lofi_solution} shows the low-fidelity optimal forcing and the corresponding block 1 displacement. We observe that the low-fidelity model is nearly identical to the target trajectory as a result of the successful optimization. However, due to the faulty modeling assumption that block 2 is stationary, the trajectory of block 1 according to the high-fidelity model is noticeably different. 

\begin{figure}[h]
\centering
  \includegraphics[width=0.49\textwidth]{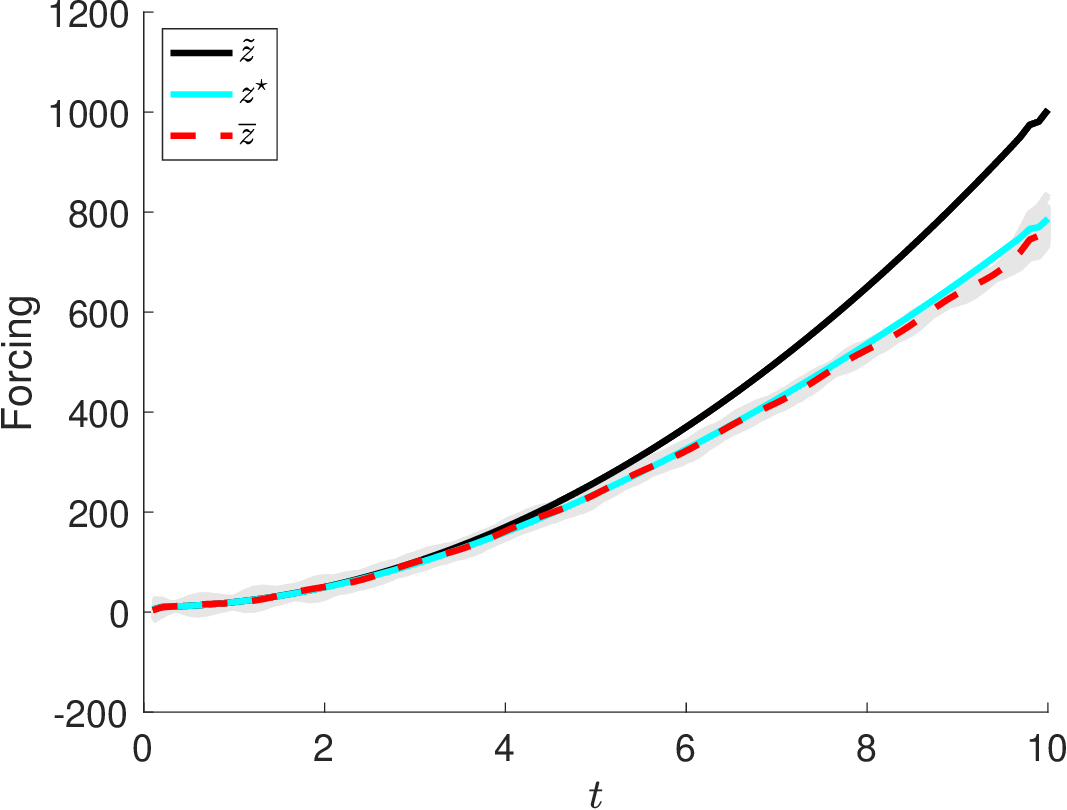}
  \includegraphics[width=0.49\textwidth]{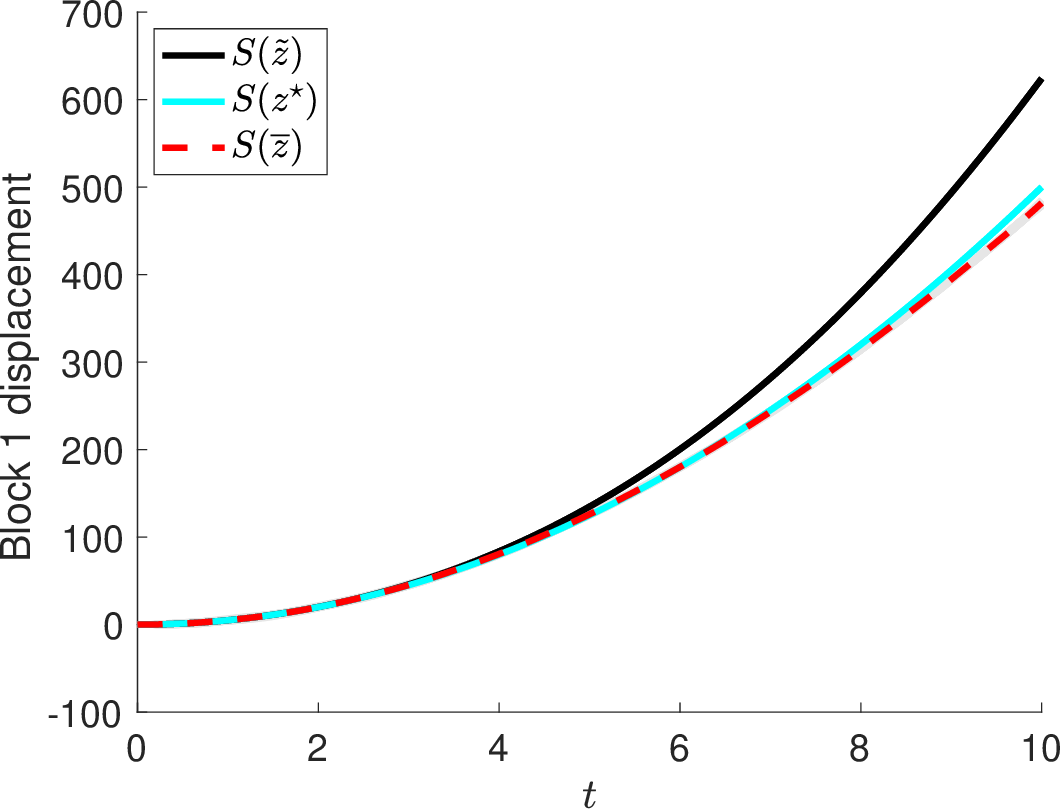}
    \caption{Left: the low and high-fidelity optimal forcings $\tilde{\z}$ and $\z^\star$, respectively, alongside the posterior optimal solution mean $\overline{\z}$ with posterior optimal solution samples (in grey); right: block 1 trajectory under the forcings from the left panel.}
  \label{fig:mass_spring_update_solution}
\end{figure}

We evaluate the high-fidelity model for $N=2$ distinct forcing terms ($\z_1=\tilde{\z}$ and $\z_2$ being randomly sampled) and seek to use the resulting block 1 displacement (and velocity) data to improve the low-fidelity optimal forcing. The prior hyper-parameter values are given in Table~\ref{tab:mass_spring_hyper-parameters}. The Hessian projection~\eqref{eqn:hess_gen_eig} uses a rank $r=17$ corresponding to four orders of magnitude decay in the eigenvalues. Figure~\ref{fig:mass_spring_update_solution} displays the optimal forcing posterior (mean and samples) and the corresponding block 1 displacement. The left panel shows the considerable improvement attained using two high-fidelity evaluations where its mean is close to the high-fidelity optimal solution with posterior samples bound tightly around it. The right panel shows the improvement in the trajectory of block 1 where we see that the posterior optimal forcing achieves trajectories which are close to the trajectory simulated using the high-fidelity optimal forcing. 

\begin{table}[!ht]
\centering
\begin{tabular}{|c|c|c|c|c|}
\hline
$\alpha_\u$ & $\beta_\u$ & $\alpha_\z$ & $\beta_\z$ & $\alpha_{\vec{d}}$  \\
$10^4$ & $5 \times 10^{-2}$ & $10^{-10}$ & $10^{-1}$ & $10^{-1}$  \\
\hline
\end{tabular}
\caption{Hyper-parameters used in the mass-spring system.}
\label{tab:mass_spring_hyper-parameters}
\end{table}

\subsection{Advection diffusion example}
\begin{figure}[h]
\centering
    \hspace{5 mm}    \includegraphics[width=0.4\textwidth]{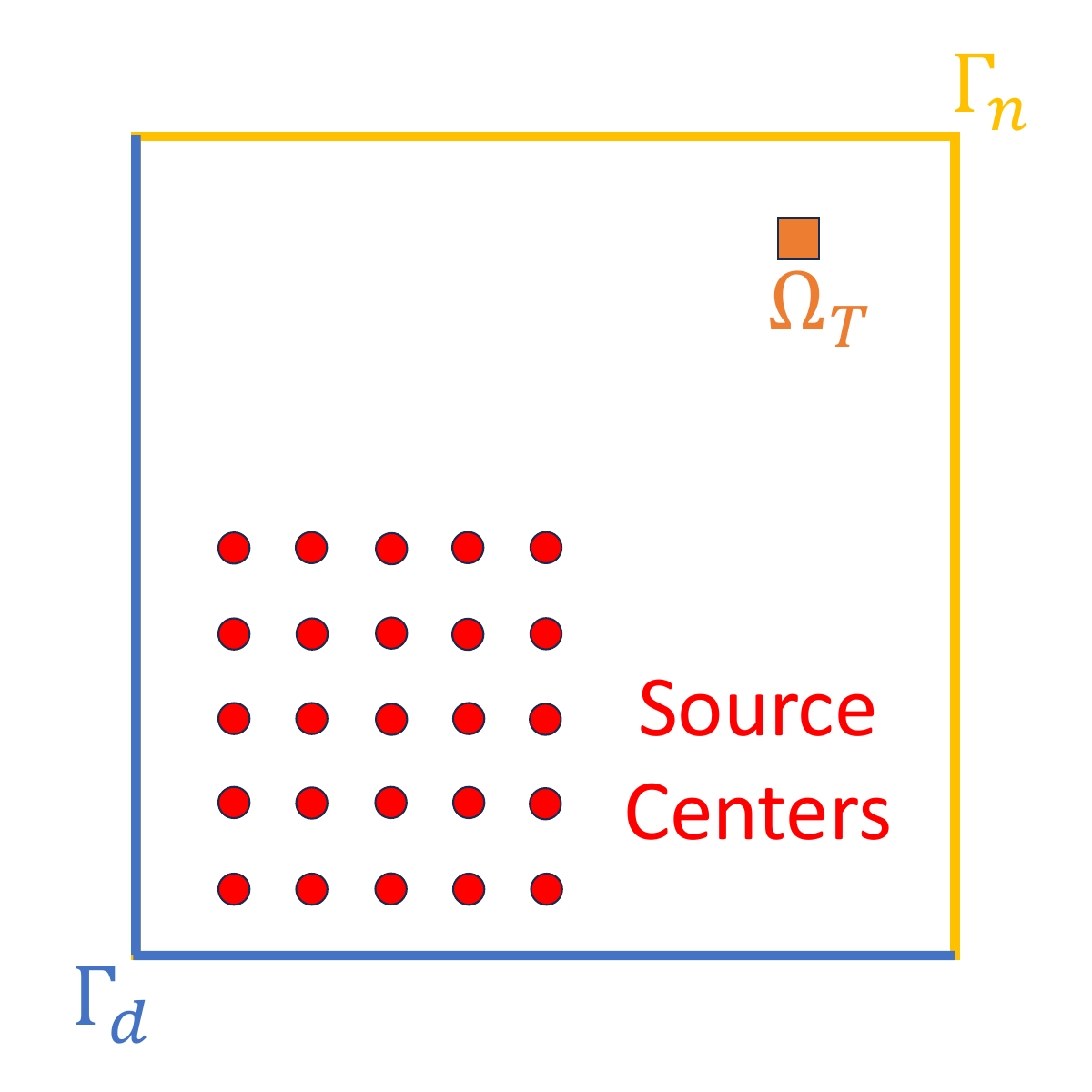} \hspace{5 mm}
  \includegraphics[width=0.49\textwidth]{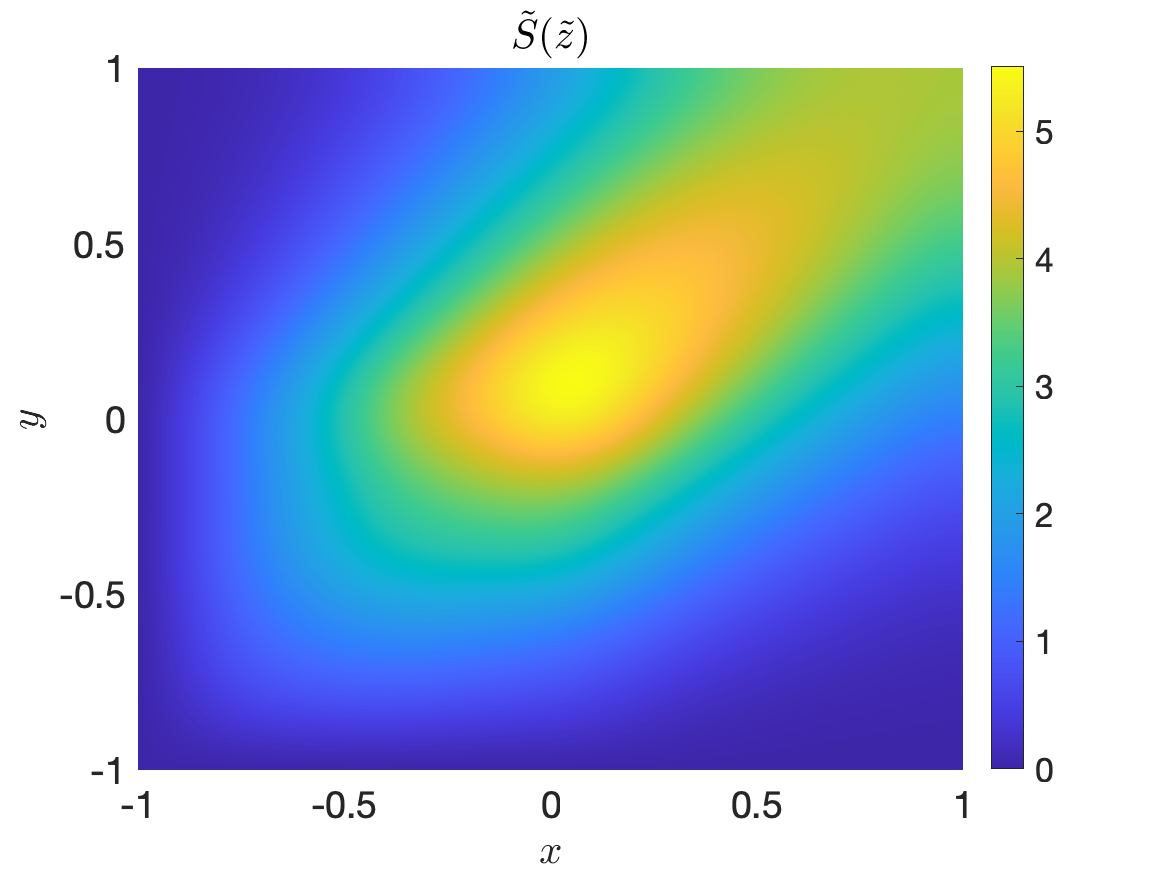}
      \includegraphics[width=0.49\textwidth]{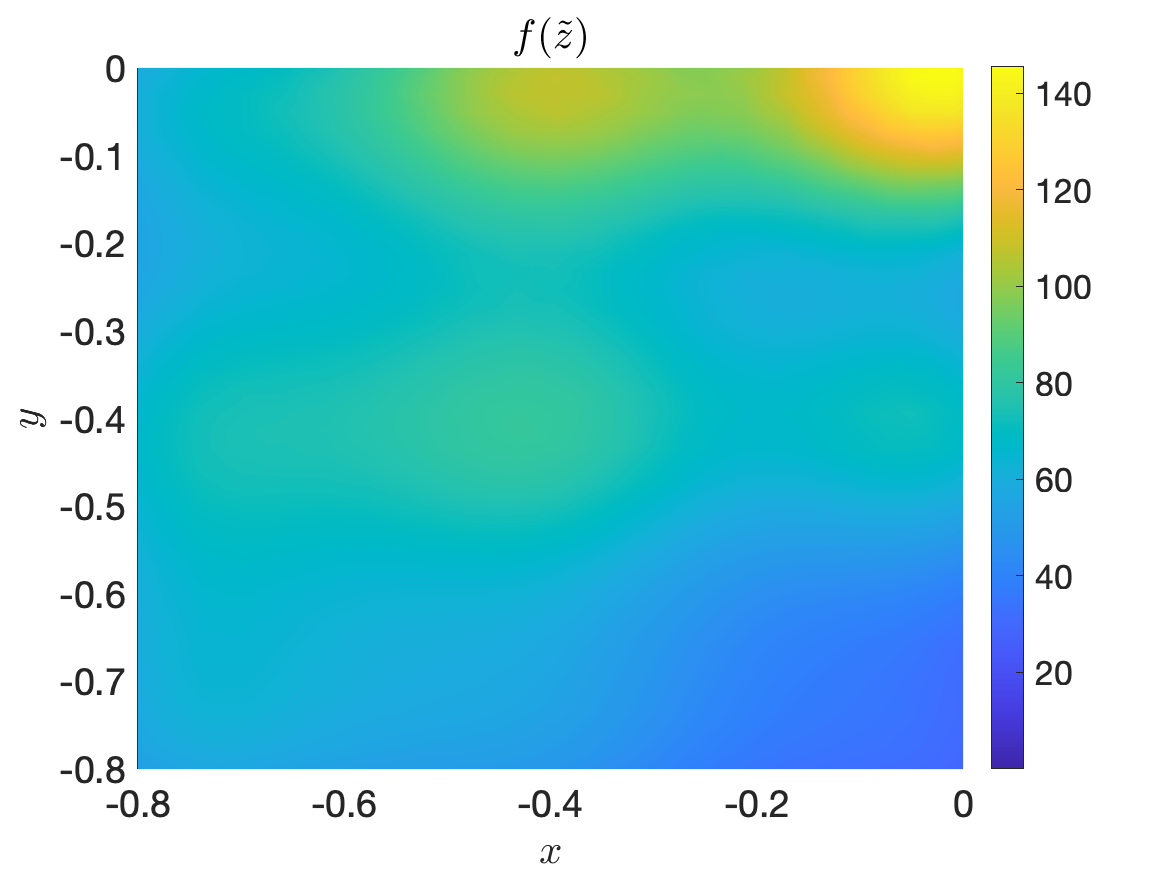}
    \includegraphics[width=0.49\textwidth]{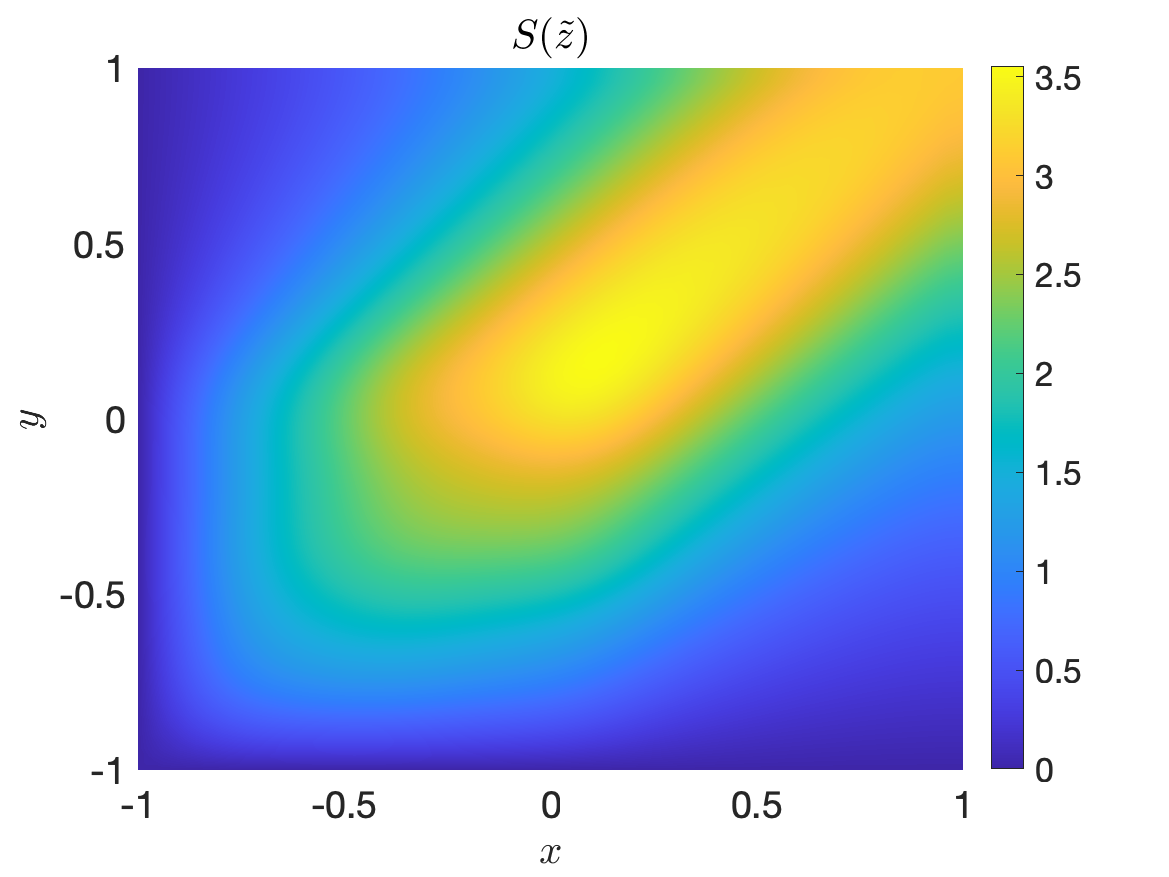}
    \caption{Left: schematic of the computational domain, source, and target locations (top), and low-fidelity optimal source $\tilde{\z}$ (bottom); right: low-fidelity model prediction (top) and and high-fidelity model prediction (bottom) evaluated at $\z=\tilde{\z}$.}
  \label{fig:adv_diff_lofi_solution}
\end{figure}

For a final demonstration, we consider a PDE in two spatial dimensions,
\begin{align*}
& -\kappa \Delta u + \vec{v}(u) \cdot \nabla u = f(\z) \qquad & \text{on } \Omega \\
& - \kappa \nabla u \cdot \vec{n} = 0 & \text{on } \Gamma_n \\
&  u = 0 & \text{on } \Gamma_d
\end{align*}
where $\kappa>0$ is the diffusion coefficient, $\vec{v}(u)$ is the velocity field, $f(\z)$ is a parameterized source control, $\Omega = (-1,1)^2$ is the domain, $\Gamma_n = \{1 \} \times (-1,1) \cup (-1,1) \times \{1\}$ is the Neumann boundary, and $\Gamma_d = \{-1 \} \times (-1,1) \cup (-1,1) \times \{-1\}$ is the Dirichlet boundary. We consider different models for the velocity field with a low-fidelity model corresponding to $\vec{v}(u) = (1,1)$ and its high-fidelity counterpart being $\vec{v}(u) = (u,u)$. In other words, the low-fidelity model is a linearization of the advection operator around the nominal state $u=1$. We seek to design a parameterized source control of the form
\begin{align*}
f(\z) = \sum\limits_{j=1}^{25} z_j \varphi_j
\end{align*}
where $\varphi_j:[-1,1]^2 \to \R$ is given by $\varphi_j(x,y) = \exp \left( -30 \left( (x-x_j)^2 + (y-y_j)^2 \right) \right)$, where $\{ (x_j,y_j) \}_{j=1}^{25}$ are uniformly placed centers on a $5 \times 5$ grid over the region $[-0.8,0.0] \times [-0.8,0.0]$. We seek to determine $\z$ such that the state is as close as possible to the target value $4.0$ over the subdomain $\Omega_T = [0.6,0.7] \times [0.8,0.9]$. The problem formulation is depicted in the top left panel of Figure~\ref{fig:adv_diff_lofi_solution}. We solve the optimization problem
\begin{align}
\label{eqn:adv_diff_example_opt_prob}
\min_{\z} J(\tilde{S}(\z),\z) \coloneqq \frac{1}{2} \int_{\Omega_T} (\tilde{S}(\z)(x) - 4.0)^2 dx + \frac{\gamma}{2} \int_{\Omega} f(\z)(x)^2 dx
\end{align}
where $\z \in \R^{25}$ is the optimization variable, $\tilde{S}(\z)$ is the solution operator for the low-fidelity model, and $\gamma=10^{-7}$ is the regularization coefficient. 

Figure~\ref{fig:adv_diff_lofi_solution} displays the solution of the low-fidelity optimization problem~\eqref{eqn:adv_diff_example_opt_prob} with the optimal source $f(\ztilde)$ in the bottom left panel and the corresponding state solution of the low-fidelity PDE in the top right panel. The solution of the high-fidelity PDE evaluated using the source $f(\ztilde)$ is shown in the bottom right panel to highlight the discrepancy between the models. Since the state $u>1$ in the high concentration regions, the low-fidelity model underestimates the advection speed. As a result, the high-fidelity model predicts that some concentration has advected out of the domain and as a result the peak concentration has a smaller magnitude.

We evaluate the high-fidelity model at $\ztilde$ and seek to use this single observation (i.e. $N=1$) of the model discrepancy to enhance the optimal source. Prior hyper-parameters are specified as discussed in Subsections~\ref{ssec:hyper-parameter_selection} with the exception that an optimization variable correlation length parameter is not needed due to the controller parameterization. Rather, we define $\W_\z=\alpha_\z^{-1} \vec{\Phi}_{basis}^T \M_\z \vec{\Phi}_{basis}$, where $\vec{\Phi}_{basis}$ is the matrix whose columns correspond to evaluating $\{\varphi_j\}_{j=1}^{25}$ at the spatial nodes. This $\W_\z$ corresponds to the $L^2$ inner product of $f(\z)$. The hyper-parameter values are reported in Table~\ref{tab:adv_diff_hyper-parameters}. 

\begin{table}[!ht]
\centering
\begin{tabular}{|c|c|c|c|}
\hline
$\alpha_\u$ & $\beta_\u$ & $\alpha_\z$ & $\alpha_{\vec{d}}$   \\
$4$ & $5 \times 10^{-1}$ & $10^{-8}$ & $10^{-2}$ \\
\hline
\end{tabular}
\caption{Hyper-parameters used in the advection diffusion example.}
\label{tab:adv_diff_hyper-parameters}
\end{table}

\begin{figure}[h]
\centering
  \includegraphics[width=0.32\textwidth]{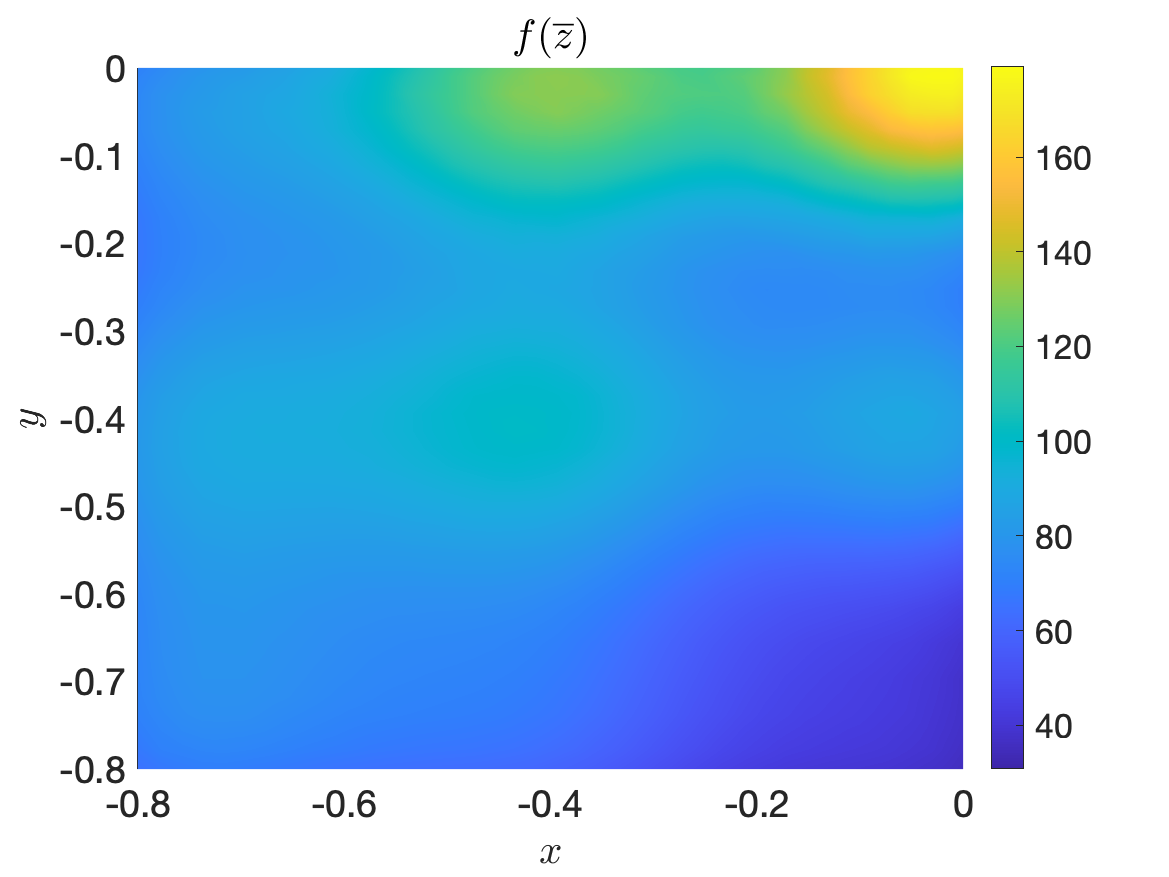}
      \includegraphics[width=0.32\textwidth]{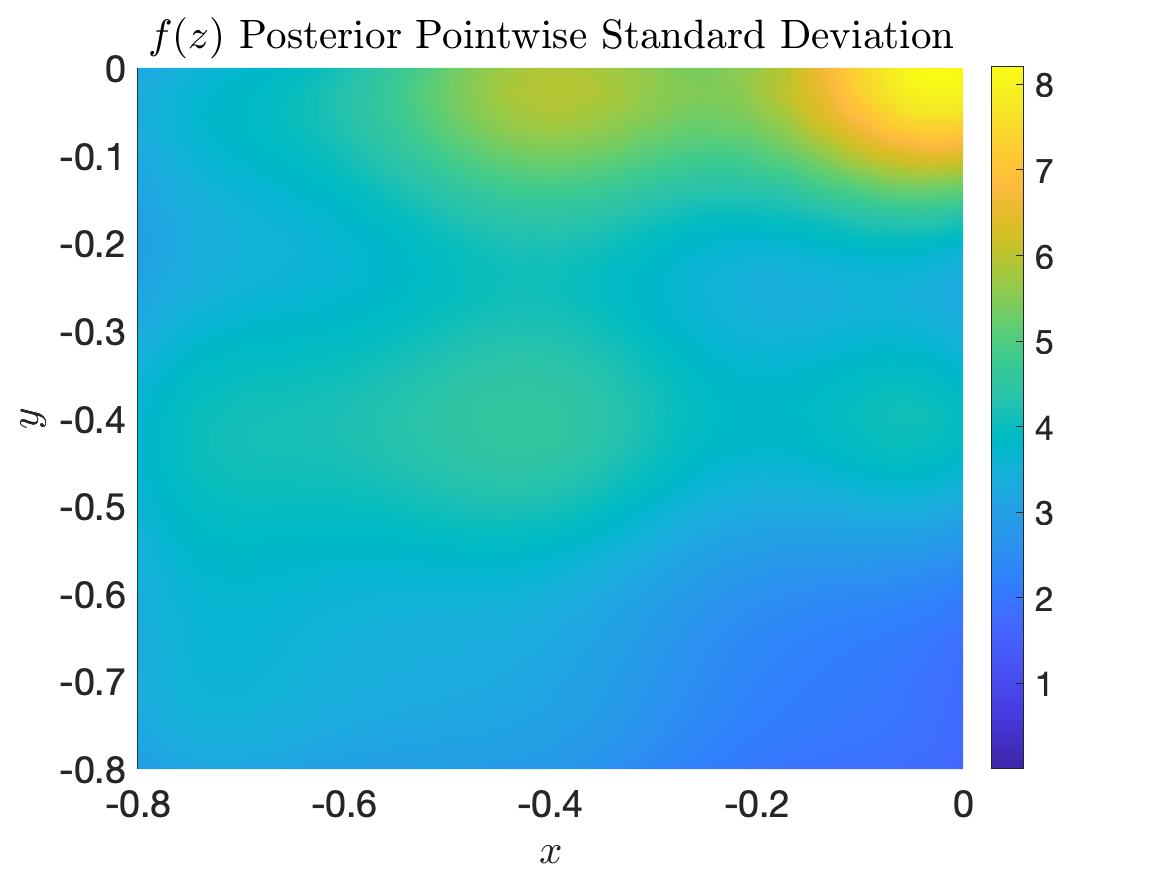}
            \includegraphics[width=0.32\textwidth]{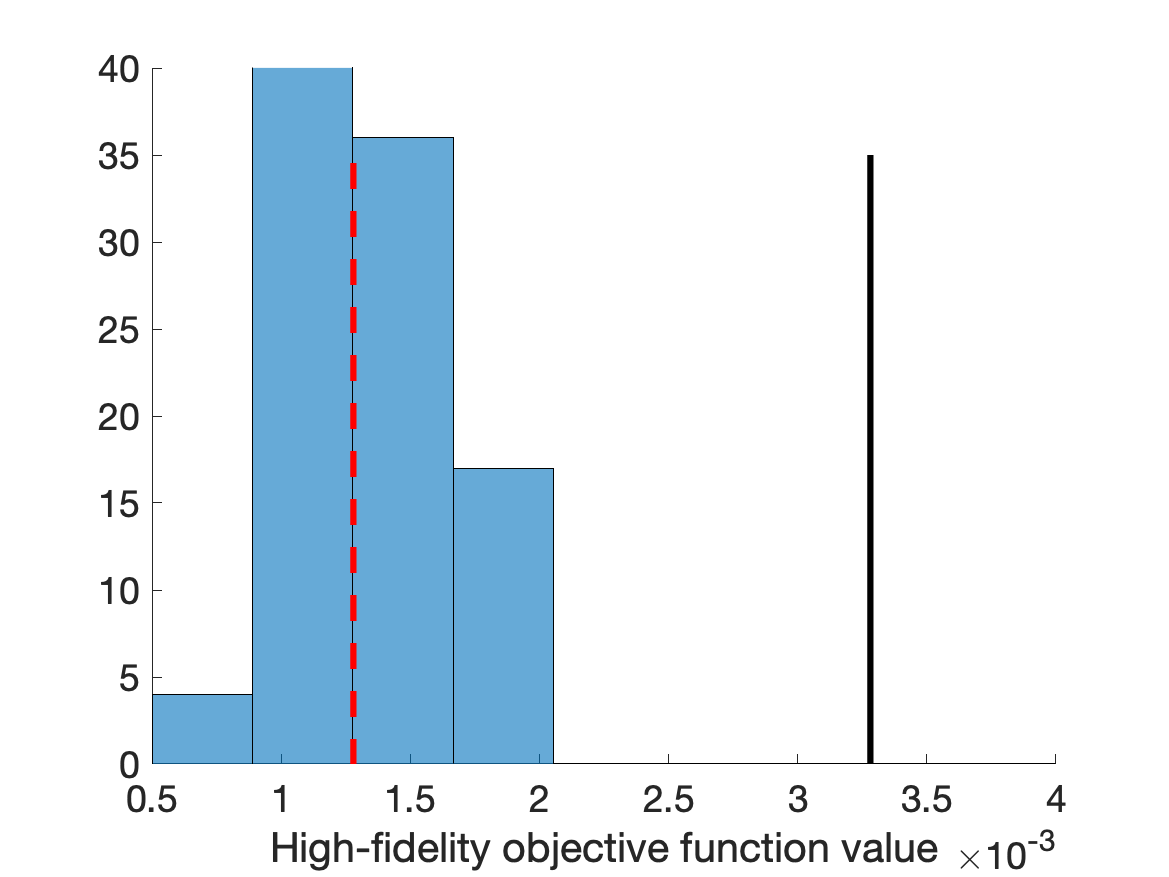}
  \includegraphics[width=0.32\textwidth]{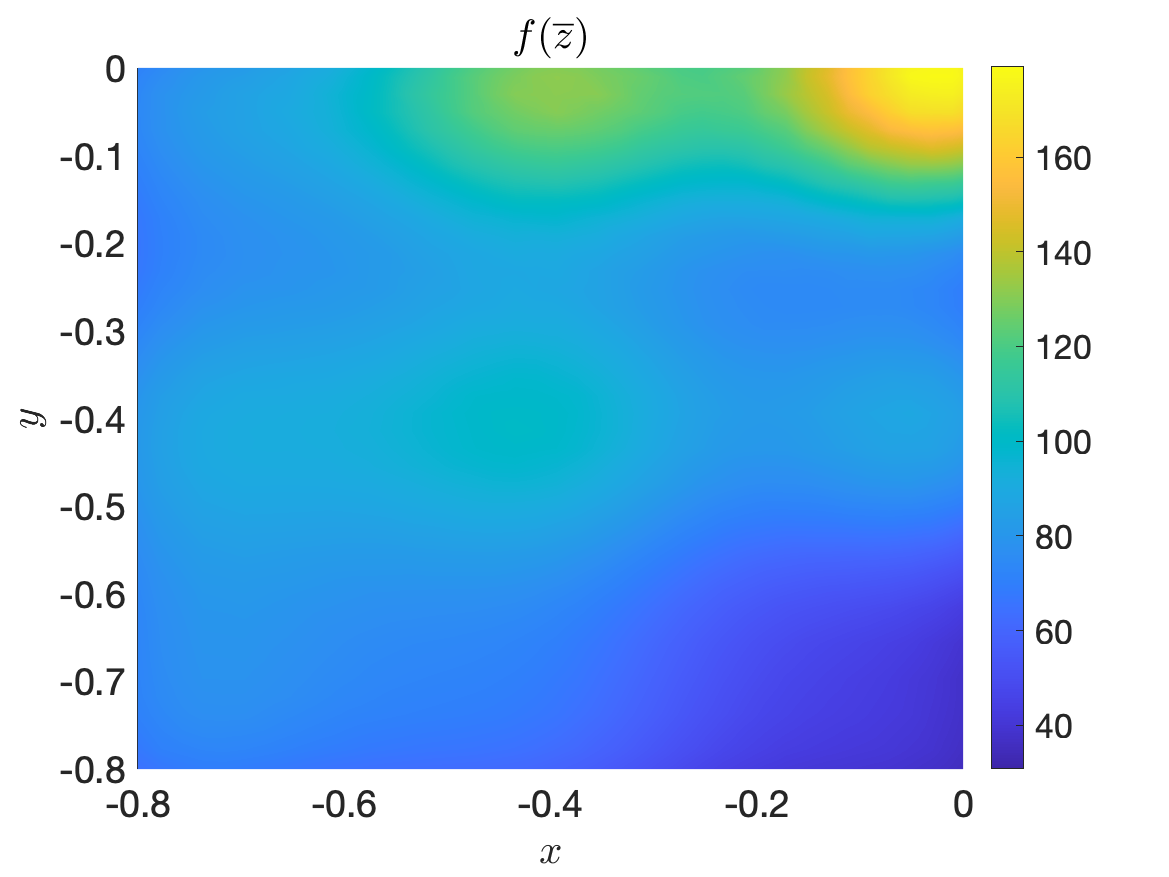}
      \includegraphics[width=0.32\textwidth]{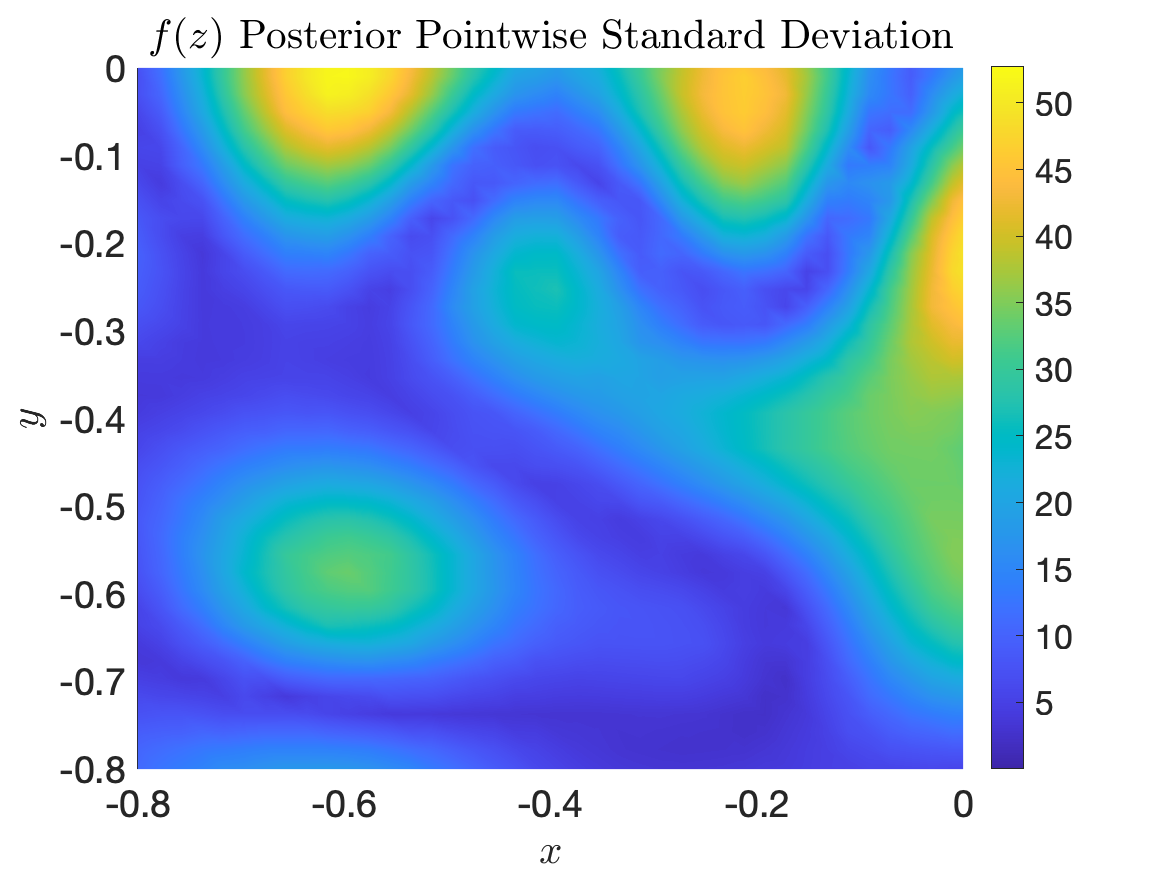}
                  \includegraphics[width=0.32\textwidth]{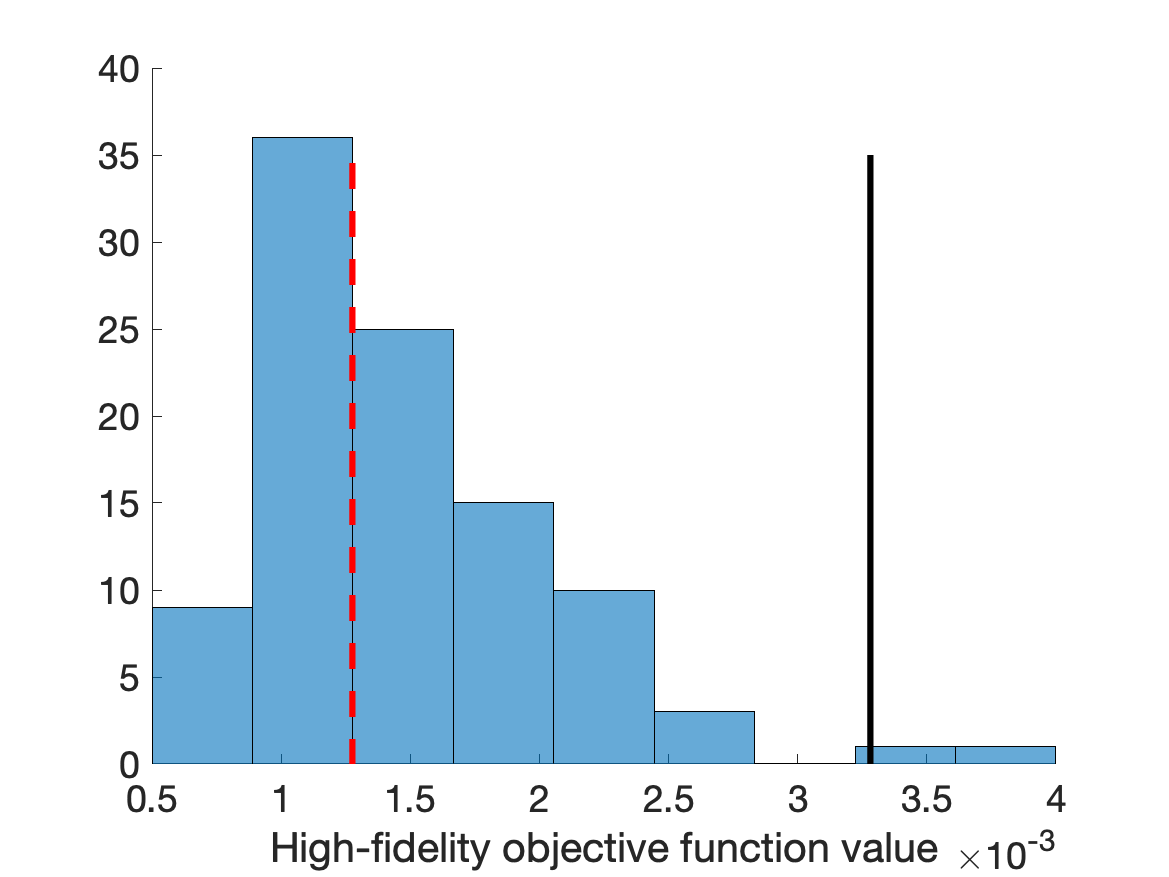}
    \caption{Comparison of the posterior optimal solution using projection subspace ranks $r=1$ (top) and $r=2$ (bottom). Left: posterior optimal solution mean $\overline{\z}$; center: posterior optimal solution pointwise standard deviation; right: histogram of the high-fidelity objective function values corresponding to the posterior optimal solution samples with the vertical lines indicating $J(S(\overline{\z}),\overline{\z})$ (red broken line) and $J(S(\tilde{\z}),\tilde{\z})$ (black solid line).}
  \label{fig:adv_diff_update}
\end{figure}

Computing the generalized eigenvalue decomposition of the Hessian, which is a $25 \times 25$ matrix as a result of the controller parameterization, reveals that the leading eigenvalue $\lambda_1=.125$ is three order of magnitude greater than the remaining 24 eigenvalues which lie within the interval $[5.2 \times 10^{-4},5.6 \times 10^{-4}]$. This corresponds to there being one dominant direction in the controller space for which the misfit $\frac{1}{2} \int_{\Omega_T} (\tilde{S}(\z)(x) - 4.0)^2 dx$ is sensitive, whereas the remaining 24 eigenvectors have little impact on the misfit but rather have their eigenvalue's magnitude dominated by the regularization term $\frac{\gamma}{2} \int_{\Omega} f(\z)(x)^2 dx$. We compute the posterior optimal solution using the Hessian subspace dimensions $r=1$ and $r=2$. Figure~\ref{fig:adv_diff_update} shows the resulting posterior optimal solution means, posterior optimal solution pointwise standard deviations, and the high-fidelity objective function values corresponding to the posterior optimal solution samples, computed with $r=1$ in the top row and $r=2$ in the bottom row. The value of the high-fidelity objective evaluated at the low-fidelity controller is $J(S(\ztilde),\ztilde)=0.0033$ and the corresponding value evaluated at the posterior optimal solution mean $\zbar$ (for both $r=1$ and $r=2$ cases) is $J(S(\overline{\z}),\overline{\z})=0.0013$. Comparing the two rows we observe that the posterior optimal solution mean is nearly identical for both cases. However, the  posterior optimal solution pointwise standard deviation increases by an order of magnitude when increasing the rank from $1$ to $2$. This corresponds to uncertainty in directions for which the misfit is insensitive. This is confirmed by comparing the high-fidelity objective function values in the rightmost column where we see a small difference in the objective function values. The $r=2$ case has a slightly higher concentration of samples which have smaller objective function value relative to the $r=1$ case. However, this difference is small relative to the order of magnitude change in the posterior optimal solution pointwise standard deviation. This highlights the fact that many optimal solutions can be considered which achieve comparable misfit values, and hence there remains considerable uncertainty in the optimal solution; however, this uncertainty has a negligible effect on its performance. 

\section{Conclusion} \label{sec:conclusion}
The proposed optimal solution updating approach is a learning algorithm constrained by the low-fidelity model that enables nontrivial improvements using a small number of high-fidelity simulations. Thanks to the closed form expressions for posterior samples, the model discrepancy and optimal solution posteriors are interpretable. These expressions aid in specification of the prior hyper-parameters and understanding their effect on the posterior. The optimal solution posterior variance may be large as a result of the high-dimensional subspace not informed by the data. However, the introduction of an optimal solution projector makes the algorithm robust by eliminating spurious variations which have minimal influence on the optimization objective. Furthermore, the projector ensures that the computation cost scales with the projector rank. By integrating low-fidelity models, domain expertise (via the prior), high-fidelity data, and the optimality system, our proposed approach is efficient for learning optimization solutions from limited data. 

Future research is needed to bring the proposed approach to its full fruition. One overarching limitation is that the post-optimality solution operator is a linear approximation of the model discrepancy to optimal solution mapping. As a result, the quality of the optimal solution update is dependent upon the low-fidelity model being ``close enough" to the high-fidelity model so that the linear approximation is valid. Ongoing work will explore the use of sequential linearization techniques to broaden the applicability of the proposed approach for complex applications where large discrepancies exists between the models. Another important problem is the design of the high-fidelity data collection. Given that the proposed approach is a learning algorithm which inputs data and predicts an optimization problem's solution, the quality of the data being input will determine the quality of what can be learned. Ongoing work is exploring the use of Bayesian optimal experimental design to direct the high-fidelity model evaluations for maximum information content. This article focused on problems where the number of high-fidelity simulations $N$ is severely limited, for instance, $N < 5$. Preliminary explorations have indicated that taking larger values for $N$ yield diminishing benefit. We observed significant gain from the first few high-fidelity simulations, but the linear approximation limits our ability to learn from additional data that is further away from the low-fidelity solution. Our ongoing research is exploring the relationship between sequential linear approximation and optimal experimental design for larger $N$.

We focused on applications for which the high-fidelity data comes from computational models. However, it is possible to use experimental data in place of high-fidelity model data in some cases. We require that the optimization variables corresponding to the high-fidelity data are known. In some applications, such as design or control of an engineering, experiments can be used to collect high-fidelity data which is then used to improve the design or control strategy. Iterating between optimization of a computational model and physical experiments with our proposed framework has potential to discover novel solutions which are only learned by fusing computational optimization with experimental data.

\section*{Acknowledgements}
This work was supported by the Laboratory Directed Research and Development program at Sandia National Laboratories, a multimission laboratory managed and operated by National Technology and Engineering Solutions of Sandia LLC, a wholly owned subsidiary of Honeywell International Inc. for the U.S. Department of Energy’s National Nuclear Security Administration under contract DE-NA0003525. SAND2024-01492O.

\section*{Appendix A}
This appendix provides a detailed derivation of the model discrepancy prior presented in Section~\ref{ssec:prior}. The squared Schatten 2-norm of $\nabla_\z \d(\cdot,\t)$ is
\begin{align*}
\vert \vert \nabla_\z \d(\cdot,\t) \vert \vert^2  & = \text{Tr}( \nabla_\z \d(\cdot,\t) ^* \nabla_\z \d(\cdot,\t)  ) \\
& = \text{Tr}(  \W_\z \nabla_\z \d(\cdot,\t)^T \W_\u \nabla_\z \d(\cdot,\t)) \\
&= \sum\limits_{k=1}^n (\W_\z^{\frac{1}{2}} \e_k,\W_\z \nabla_\z \d(\cdot,\t)^T \W_\u \nabla_\z \d(\cdot,\t) \W_\z^{\frac{1}{2}} \e_k)_{\W_\z^{-1}}
\end{align*}
where $\{\e_k\}_{k=1}^n$ is the canonical basis for $\R^n$ ($\{\W_\z^{\frac{1}{2}} \e_k\}_{k=1}^n$ is an orthonormal basis in the $\W_\z^{-1}$ weighted inner product). Note that $\W_\u$ and $\W_\z$ arise from the adjoint and trace since the operations are computed in the function spaces (rather than Euclidean space). Hence,
\begin{align*}
\vert \vert \nabla_\z \d(\cdot,\t) \vert \vert^2  = \sum\limits_{k=1}^n \e_k^T \W_\z^\frac{1}{2} \nabla_\z \d(\cdot,\t)^T \W_\u \nabla_\z \d(\cdot,\t) \W_\z^\frac{1}{2} \e_k 
\end{align*}
and it follows that
\begin{align*}
\vert \vert \nabla_\z \d(\cdot,\t) \vert \vert^2 = \t^T
\left( \begin{array}{cc}
\vec{0} & \vec{0}  \\
\vec{0} & \W_\u \otimes \M_\z \W_\z \M_\z 
\end{array} \right) 
\t .
\end{align*}
By observing that $\vert \vert \d(\ztilde,\t) \vert \vert_{\W_\u}^2$ can be represented as 
\begin{align*}
\vert \vert \d(\ztilde,\t) \vert \vert_{\W_\u}^2 = \t^T
\left( \begin{array}{cc}
\W_\u & \W_\u \otimes \ztilde^T \M_z \\
\W_\u \otimes \M_z \ztilde & \W_\u \otimes \M_\z \ztilde \ztilde^T \M_\z 
\end{array} \right) 
\t ,
\end{align*}
the norm~\eqref{eqn:discrepancy_norm} can be written in a quadratic form 
$$\vert \vert \d(\ztilde,\t) \vert \vert_{\W_\u}^2 + \vert \vert \nabla_\z \d(\cdot,\t) \vert \vert^2 = \t^T \W_\t \t,$$ 
where $\W_\t$ is defined as in~\eqref{eqn:M_theta}.

\section*{Appendix B}
This appendix provides detailed derivations to prove Theorems~\ref{thm:post_mean},~\ref{thm:post_samples},~\ref{thm:delta_breve_samples}, and~\ref{thm:B_breve_samples}. The guiding principles for our analysis is that we seek to avoid computation in $\R^p$, that is matrix or vector operations involving the parameter vector $\t$, but rather we seek to express the posterior samples in terms of Kronecker products of computation in the state and optimization variable spaces. 

\subsection*{Proof of Theorem~\ref{thm:post_mean}}
The proof consists of three parts:
\begin{enumerate}
\item factorizing $\W_\u$ and $\W_\t$ to facilitate derivations,
\item factorizing $\A$ through the use of a $\W_\t$ orthogonal basis to rewrite $\vec{\Sigma}^{-1}$,
\item inverting $\vec{\Sigma}^{-1}$ by exploiting orthogonality to invert a sum.
\end{enumerate}

\subsubsection*{Factorizing $\W_\u$ and $\W_\t$}
Let $(\vec{x}_j,\lambda_j)$ be the generalized eigenvalues and eigenvectors of $\W_\u$ in the $\M_\u$ weighted inner product. Let $\vec{X} = \begin{pmatrix} \vec{x}_1 & \vec{x}_2 & \cdots & \vec{x}_m \end{pmatrix}$ and $\vec{\Lambda}$ be the diagonal matrix with entries $\lambda_1,\lambda_2,\dots,\lambda_m$, then the generalized eigenvalue decomposition is given by $\W_\u = \M_\u \vec{X} \vec{\Lambda} \vec{X}^T \M_\u$. We may decompose $\W_\t = \L \L^T$ where
\begin{align}
\label{eqn:L}
\L =
\left(
\begin{array}{cc}
\M_\u \vec{X} \vec{\Lambda}^\frac{1}{2} & \vec{0} \\
\M_\u \vec{X} \vec{\Lambda}^\frac{1}{2}  \otimes \M_\z \ztilde & \M_\u \vec{X} \vec{\Lambda}^\frac{1}{2}  \otimes \M_\z \W_\z^{\frac{1}{2}} 
\end{array}
\right) .
\end{align}
Furthermore, thanks to its block lower triangular structure, the inverse of $\L$ is
\begin{align}
\label{eqn:Linv}
\L^{-1} =
\left(
\begin{array}{cc}
\vec{\Lambda}^{-\frac{1}{2}} \vec{X}^T & \vec{0} \\
\vec{\Lambda}^{-\frac{1}{2}} \vec{X}^T   \otimes \left( -\W_\z^{-\frac{1}{2}} \ztilde \right) & \vec{\Lambda}^{-\frac{1}{2}} \vec{X}^T  \otimes  \W_\z^{-\frac{1}{2}} \M_\z^{-1} 
\end{array}
\right) 
\end{align}
and the inverse of $\W_\t$ may be computed by multiplying $\W_\t^{-1} = \L^{-T} \L^{-1}$.

\subsubsection*{Factorizing $\A$}
To write the posterior precision matrix $\vec{\Sigma}^{-1}$ in a form amenable for inversion, we decompose $\A$ in a chosen inner product. In particular, the generalized singular value decomposition (GSVD) of $\A$ is used with the $\I_N \otimes \M_\u$ weighted inner product on its row space and $\W_\t$ weighted inner product on its column space. This gives $\A = \vec{\Xi} \vec{\Phi} \vec{\Psi}^T \W_\t$, where  $\vec{\Phi}$ is the diagonal matrix of singular values, and $\vec{\Xi}$ and $\vec{\Psi}$ are matrices containing the singular vectors which satisfy $\vec{\Xi}^T (\I_N \otimes \M_\u) \vec{\Xi} = \I_{mN}$ and $\vec{\Psi}^T \W_\t \vec{\Psi} = \I_p$. To determine the singular vectors, note that 
\begin{align*}
( \I_N \otimes \M_\u) \A \W_\t^{-1} \A^T ( \I_N \otimes \M_\u)= ( \I_N \otimes \M_\u) \vec{\Xi} \vec{\Phi}^2 \vec{\Xi}^T( \I_N \otimes \M_\u).
\end{align*}
 Hence determine $\vec{\Xi}$ and $\vec{\Phi}$ from the generalized eigenvalue decomposition of
$$( \I_N \otimes \M_\u) \A \W_\t^{-1} \A^T ( \I_N \otimes \M_\u)$$
 in the $(\I_N \otimes \M_\u)$ weighted inner product, or equivalently the eigenvalue decomposition of $\A \W_\t^{-1} \A^T ( \I_N \otimes \M_\u)$. Using the expressions~\eqref{eqn:Aell} for $\A$ and~\eqref{eqn:Linv} for $\W_\t^{-1}=\L^{-T} \L^{-1}$, we can write
\begin{eqnarray*}
 \A \W_\t^{-1} \A^T ( \I_N \otimes \M_\u) =\vec{G} \otimes  \W_\u^{-1} \M_\u
\end{eqnarray*}
where
\begin{align*}
& \G= \e \e^T + (\vec{Z} - \ztilde \e^T)^T \W_\z^{-1} (\vec{Z}-\ztilde \e^T) \in \R^{N \times N},
\end{align*}
$\vec{Z} = \begin{pmatrix} \z_1 & \z_2 & \dots & \z_N \end{pmatrix}$ is the matrix of optimization variable data, and $\e \in \R^N$ is the vector of ones, as introduced in Subsection~\ref{subsec:posterior_sampling_expressions}. Hence the left singular vectors $\vec{\Xi}$ correspond to the eigenvectors of $\vec{G} \otimes \W_\u^{-1} \M_\u$ and the squared singular values $\vec{\Phi}^2$ correspond to the eigenvalues.

Let $(\vec{g}_i,\mu_i)$ be the eigenvectors and eigenvalues of $\vec{G}$, and note that $(\vec{x}_j,\frac{1}{\lambda_j})$ are the eigenvectors and eigenvalues of $\W_\u^{-1} \M_\u$. Recalling properties of the eigenvalue decomposition of a Kronecker product, we observe that the squared generalized singular values of $\vec{A}$ are given by $\frac{\mu_i}{\lambda_j}$ and are associated with the left singular vectors $\vec{\xi}_{i,j}=\vec{g}_i \otimes \vec{x}_j$. Rewriting the GSVD to solve for $\vec{\Psi}$, the right singular vector associated with $\vec{\xi}_{i,j}$ is given by
 \begin{eqnarray}
 \label{eqn:psi_ij}
\vec{\psi}_{i,j}= \frac{1}{\sqrt{\mu_i \lambda_j}} 
\left(
 \begin{array}{cc}
s_i  \vec{x}_j \\
 \vec{x}_j \otimes \M_\z^{-1} \W_\z^{-1} \vec{y}_i
 \end{array}
 \right),
 \end{eqnarray}
$ i=1,2,\dots,N \ \ j=1,2,\dots,m,$ where
\begin{eqnarray}
\label{eqn:w_vecs}
\vec{y}_i = \vec{Z} \vec{g}_i - (\e^T \vec{g}_i) \ztilde \qquad \text{and} \qquad s_i = (\e^T \vec{g}_i) - \vec{y}_i^T \W_\z^{-1} \ztilde .
\end{eqnarray}

Given these decompositions, we express $\vec{\Sigma}^{-1}$ as
\begin{eqnarray}
\label{eqn:sigma_inv_factor}
\vec{\Sigma}^{-1} =  \frac{1}{\alpha_{\vec{d}}} \L \C \L^T 
\end{eqnarray}
where
\begin{eqnarray}
\label{eqn:X}
\C =  \alpha_{\vec{d}} \vec{I} +  \L^T \vec{\Psi} \vec{\Phi}^2 \vec{\Psi}^T \L .
\end{eqnarray}

\subsubsection*{Inverting $\vec{\Sigma}^{-1}$}
This factorization of $\vec{\Sigma}^{-1}$ implies that 
\begin{eqnarray}
\label{eqn:sigma_factors}
\vec{\Sigma}= \alpha_{\vec{d}} \L^{-T} \C^{-1} \L^{-1}
\end{eqnarray}
 and facilitates manipulation because $\C$ is composed of diagonal and orthogonal matrices. Applying the Sherman-Morrison-Woodbury formula to~\eqref{eqn:X}, we have
\begin{eqnarray*}
\C^{-1} = \frac{1}{\alpha_{\vec{d}}} \left( \vec{I} - \L^T \vec{\Psi} \vec{D}\vec{\Psi}^T \L \right)
\end{eqnarray*}
where $\vec{D}\in \R^{mN \times mN}$ is a diagonal matrix whose entries are given by $\frac{\mu_i}{\mu_i + \alpha_{\vec{d}} \lambda_j}.$
Algebraic simplifications yields
\begin{eqnarray}
\label{eqn:Sigma}
\vec{\Sigma} = \W_\t^{-1} - \vec{\Psi} \vec{D}\vec{\Psi}^T .
\end{eqnarray}

Theorem~\ref{thm:post_mean} follows from~\eqref{eqn:Sigma},~\eqref{eqn:post_mean_product}, and manipulations to rewrite sums over eigenpairs in terms of linear system solves.

\subsection*{Proof of Theorem~\ref{thm:post_samples}}
We build off the previous proof and extend it in two parts:
\begin{enumerate}
\item factorizing $\vec{\Sigma}$ by determining its eigenvalue decomposition, 
\item sampling the posterior through matrix-vector products in $\R^m$ and $\R^n$.
\end{enumerate}

\subsubsection*{Factorizing $\vec{\Sigma}$}
We seek a factorization of $\vec{\Sigma}$ so that samples are computed via a matrix-vector product of the covariance factor with a standard normal random vector. Such a factorization of $\vec{\Sigma}$ is not easily attained from~\eqref{eqn:Sigma}. Rather, we go back to~\eqref{eqn:sigma_factors} and seek to compute the eigenvalue decomposition of $\C$. 

$(\L^T \vec{\Psi})^T ( \L^T \vec{\Psi}) = \vec{\Psi}^T \W_\t \vec{\Psi} = \vec{I}_p$, and hence the columns of $\L^T \vec{\Psi}$ are orthonormal in the Euclidean inner product. This implies that $\L^T \vec{\psi}_{i,j}$ is an eigenvector of $\C$ with eigenvalue $\alpha_{\vec{d}} +\frac{\mu_i}{\lambda_j}$. By multiplying~\eqref{eqn:L} and~\eqref{eqn:psi_ij}, we express these eigenvectors as
\begin{eqnarray}
\label{eqn:Lt_psi_ij}
\L^T \vec{\psi}_{i,j} = 
\frac{1}{\sqrt{\mu_i}}
\left(
\begin{array}{c}
 (\vec{e}^T \vec{g}_i) \vec{e}_j \\
\vec{e}_j \otimes  \W_\z^{-\frac{1}{2}} \vec{y}_i
\end{array}
\right), 
\end{eqnarray}
$j=1,2,\dots,m, \qquad i = 1,2,\dots,N .$

This gives a total of $mN$ eigenpairs of $\C \in \R^{p \times p}$. There are $p-mN$ additional eigenpairs of $\C$.  These have repeated eigenvalue $\alpha_{\vec{d}}$ and eigenvectors from a set of orthonormal vectors that are orthogonal to the columns of $\{ \L^T \vec{\psi}_{i,j} \}$. The determination of these eigenvalues motivates the assumption that $\z_1=\ztilde$ since the assumption implies that $\text{span}\{ \vec{y}_i\}_{i=1}^N = \text{span}\{ \z_\ell - \ztilde \}_{\ell=2}^N$ defines a $N-1$ dimensional subspace and hence ensures existence of a set of orthonormal vectors $\{ \breve{\z}_1,\breve{\z}_2,\dots,\breve{\z}_{n-N+1} \} \subset \R^n$ which are orthogonal to $\{ \W_\z^{-\frac{1}{2}} (\z_\ell-\ztilde) \}_{\ell=2}^N$. Then remaining eigenvectors of $\C$ are given by
 \begin{align}
 \label{eqn:X_evec_ortho}
 \left(
 \begin{array}{c}
 \vec{0} \\
 \vec{e}_j \otimes \breve{\z}_k
 \end{array}
 \right), \qquad j=1,2,\dots,m, \qquad k=1,2,\dots,n-N+1.
 \end{align}
 
Collecting this set of $p$ eigenvectors, through the union of~\eqref{eqn:Lt_psi_ij} and~\eqref{eqn:X_evec_ortho}, in an orthogonal matrix $\vec{Q}$ and the corresponding eigenvalues into a diagonal matrix $\vec{\Upsilon}$, the eigenvalue decomposition is given by
\begin{eqnarray*}
\C = \vec{Q} \vec{\Upsilon} \vec{Q}^T.
\end{eqnarray*}
Hence, recalling~\eqref{eqn:sigma_factors}, we decompose the covariance matrix as
\begin{eqnarray}
\label{eqn:Sigma_factor}
\vec{\Sigma} = \vec{T} \vec{T}^T
\end{eqnarray}
where 
\begin{eqnarray}
\label{eqn:T}
\vec{T} =  \sqrt{\alpha_{\vec{d}}} \L^{-T} \vec{Q} \vec{\Upsilon}^{-\frac{1}{2}}.
\end{eqnarray}

\subsubsection*{Computing posterior samples}
The expression~\eqref{eqn:T} will facilitate efficient sampling as posterior samples take the form $\overline{\t} + \vec{T} \vec{\omega}$, where $\vec{\omega} \sim \mathcal N(\vec{0},\vec{I}_p)$ follows a standard normal distribution in $\R^p$. To prove Theorem~\ref{thm:post_samples}, we make the additional observation that if $\vec{\nu} \sim \mathcal N(\vec{0},\vec{I}_p)$, then we can express $\vec{\omega}$ as
\begin{align*}
\vec{\omega} = \left( \vec{I}_{n+1} \otimes \M_\u^\frac{1}{2} \vec{X} \right)^T \vec{\nu},
\end{align*}
 since $(\vec{I}_{n+1} \otimes \M_\u^\frac{1}{2} \vec{X})^T \in \R^{p \times p}$ is an orthogonal matrix and Gaussian distributions are invariant under orthogonal transformations. Matrix multiplication to compute $\vec{T} \vec{\omega} = \vec{T} \left( \vec{I}_{n+1} \otimes \M_\u^\frac{1}{2} \vec{X} \right)^T \vec{\nu}$ completes the proof of Theorem~\ref{thm:post_samples}.

\subsection*{Proof of Theorem~\ref{thm:delta_breve_samples}}
Let $\breve{\P}$ be the orthogonal projector onto $\text{span}\{ \breve{\z}_k \}_{k=1}^{n-N+1}$. Let $\{ \breve{\z}_k \}_{k=n-N+2}^n$ be an orthonormal basis for $\text{span} \{  \W_\z^{-\frac{1}{2}} (\z_\ell-\ztilde) \}_{\ell=2}^N$ so that $\{ \breve{\z}_k\}_{k=1}^n$ is an orthonormal basis for $\R^n$. Observe that $\breve{\P} \breve{\z}_k=\breve{\z}_k$ for $k=1,2,\dots,n-N+1$ and $\breve{\P} \breve{\z}_k=\vec{0}$ for $k=n-N+2,\dots,n$. Furthermore, recall that $\breve{\u}_k$ is sampled as $\breve{\u}_k =  \X \vec{\Lambda}^{-\frac{1}{2}} \vec{\omega}_k$, where $\vec{\omega}_k \sim \mathcal N(\vec{0},\I_m)$. Then we have
\begin{align*}
   \breve{\d}(\z) &= \sum\limits_{k=1}^{n-N+1} \left( \breve{\z}_k^T \W_\z^{-\frac{1}{2}} (\z-\ztilde) \right) \breve{\u}_k \\
   & =  \sum\limits_{k=1}^{n} \left( \breve{\z}_k^T \breve{\P} \W_\z^{-\frac{1}{2}} (\z-\ztilde) \right) \X \vec{\Lambda}^{-\frac{1}{2}} \vec{\omega}_k \\
   & = \X \vec{\Lambda}^{-\frac{1}{2}}  \sum\limits_{k=1}^{n} \left( \breve{\z}_k^T \breve{\P} \W_\z^{-\frac{1}{2}} (\z-\ztilde) \right) \vec{\omega}_k 
\end{align*}
where $\vec{\omega}_k \sim \mathcal N(\vec{0},\I_m)$, $k=1,2,\dots,n$, are independent identically distributed (i.i.d.) Gaussian random vectors. Notice that $\sum_{k=1}^{n} \left( \breve{\z}_k^T \breve{\P} \W_\z^{-\frac{1}{2}} (\z-\ztilde) \right) \vec{\omega}_k $ is a linear combination of Gaussian random vectors, which is itself a Gaussian random vector with mean $\vec{0}$ and covariance $( \sum_{k=1}^{n} ( \breve{\z}_k^T \breve{\P} \W_\z^{-\frac{1}{2}} (\z-\ztilde) )^2 ) \I_m$. Hence, we can sample $\breve{\d}(\z)$ by computing
\begin{align}
\label{eqn:delta_breve_post_sample}
\sqrt{ \sum\limits_{k=1}^{n} \left( \breve{\z}_k^T \breve{\P} \W_\z^{-\frac{1}{2}} (\z-\ztilde) \right)^2 } \X \vec{\Lambda}^{-\frac{1}{2}}  \vec{\omega} ,
\end{align}
where $\vec{\omega} \sim \mathcal N(\vec{0},\I_m)$. To compute the coefficient, which depends the computationally intractable basis functions $\{ \breve{\z}_k \}_{k=1}^n$, we employ Parseval's identity to write
\begin{align*}
\sum\limits_{k=1}^{n} \left( \breve{\z}_k^T \breve{\P} \W_\z^{-\frac{1}{2}} (\z-\ztilde) \right)^2  & = \vert \vert \breve{\P} \W_\z^{-\frac{1}{2}} (\z-\ztilde) \vert \vert_2^2 .
\end{align*}

Constructing $\breve{\P}$ explicitly as the identity minus the orthogonal projector onto $\text{span} \{  \W_\z^{-\frac{1}{2}} (\z_\ell-\ztilde) \}_{\ell=2}^N$, we arrive at the final computable expression
\begin{align*}
\vert \vert \breve{\P} \W_\z^{-\frac{1}{2}} (\z-\ztilde) \vert \vert_2^2 = (\z-\ztilde)^T \left( \W_\z^{-1} - \W_\z^{-1} \vec{Z}_c \left(  \vec{Z}_c^T \W_\z^{-1}  \vec{Z}_c \right)^{-1}  \vec{Z}_c^T \W_\z^{-1} \right) (\z-\ztilde)
\end{align*}
where $ \vec{Z}_c = \begin{pmatrix}  \z_2-\ztilde & \cdots & \z_N - \ztilde \end{pmatrix} \in \R^{n \times (N-1)}$.

\subsection*{Proof of Theorem~\ref{thm:B_breve_samples}}
Observe that $\B \breve{\t}$ is a linear combination of the vectors $ \W_\z^{-\frac{1}{2}} \breve{\z}_k$ with coefficients
\begin{eqnarray*}
\nabla_\u J \breve{\u}_k  = \nabla_\u J \W_\u^{-\frac{1}{2}}  \vec{\nu}_{N+k}.
\end{eqnarray*}
These coefficients are i.i.d. normally distributed (scalar) random variables with mean 0 and variance $ \nabla_\u J \W_\u^{-1} \nabla_\u J^T$. Accordingly, we consider the equivalent sample
\begin{eqnarray*}
\sqrt{\nabla_\u J \W_\u^{-1} \nabla_\u J^T} \W_\z^{-\frac{1}{2}}  \sum\limits_{k=1}^{n-N+1} c_k \breve{\z}_k = \sqrt{\nabla_\u J \W_\u^{-1} \nabla_\u J^T} \W_\z^{-\frac{1}{2}}  \breve{\P} \sum\limits_{k=1}^{n} c_k  \breve{\z}_k
\end{eqnarray*}
where $c_k$, $k=1,2,\dots,n$, are i.i.d. samples from a standard normal distribution. Since $\{ \breve{\z}_k\}_{k=1}^n$ is an orthonormal basis for $\R^n$, the sum over an orthonormal basis with i.i.d standard normal coefficients corresponds to rotating a Gaussian and does not change its distribution. Hence we equivalently compute samples as
\begin{eqnarray*}
\sqrt{\nabla_\u J \W_\u^{-1} \nabla_\u J^T} \W_\z^{-\frac{1}{2}}  \breve{\P} \vec{\omega}_n
\end{eqnarray*}
where $\vec{\omega}_n \sim \mathcal N(\vec{0},\I_n)$. Writing $\breve{\P}$ as the identity minus the orthogonal projector onto $\text{span} \{  \W_\z^{-\frac{1}{2}} (\z_\ell-\ztilde) \}_{\ell=2}^N$ and manipulating matrix multiplication completes the proof.

\section*{Appendix C}
This appendix provides a detailed computational cost analysis. We introduce computational cost parameters and provide expressions for the total cost. Cost is measured in the number of PDE solves as we assume all other computation is negligible. Let $f$ and $\tilde{f}$ denote the computational cost of a high and low-fidelity forward solve, respectively. Let $\tilde{a}$ denote the cost of a low-fidelity adjoint solve, that is, the cost of a linear system solve where the coefficient matrix is the state Jacobian of the PDE constraint. Let $e_\u$ and $e_\z$ denote the computational cost of the elliptic solves involved in the prior covariances $\W_\u^{-1}$ and $\W_\z^{-1}$.

We assume that the low-fidelity model is optimized using a trust region Newton-CG algorithm, let $\tilde{n}_\text{iter}$ denote the number of iterations required by it, and let $\tilde{n}_\text{adjoint}$ denote the average number of Hessian-vector products required in each iteration. Then the low-fidelity optimization has a total computational cost of 
\begin{align*}
\tilde{O} = \tilde{n}_\text{iter} \tilde{f} + \tilde{n}_\text{iter}(1+2\tilde{n}_\text{adjoint}) \tilde{a}.
\end{align*}
 The ``1" arises from the adjoint solve needed for gradient computation and the ``2" corresponds to the incremental state and incremental adjoint solves required for a Hessian-vector product (which we assume have the cost $\tilde{a}$).

We compare with the computational cost of Algorithm~\ref{alg:post_discrepancy_samples} with $s$ samples, a rank $r$ projector, and rank $q$ approximation of the state prior elliptic operator. Line~2 costs $2(q+\ell_E)e_\u $, where $\ell_E$ is the oversampling factor required to compute the GSVD of $\E_\u^{-1}$ using a randomized algorithm~\cite{saibaba_gsvd}. Typically, $10 \le \ell_E \le 20$ is sufficient. If $\W_z^{-1}$ is defined as a squared inverse elliptic operator, line~3 costs $2N e_\z$ and line~4 costs $s e_\z$. Line~5 requires $s+1$ matrix-vector products with $\nabla_\z \tilde{S}$, whose cost is $(s+1)\tilde{a}$ since $\nabla_\z \tilde{S}$ takes the form of an adjoint solve with a different right hand side. The Hessian-vector products in line~6 cost $4(r+\ell_H)\tilde{a}$, where $\ell_H$ is the oversampling factor required to compute the GEVD using a randomized algorithm. Additionally, $2(r+\ell_H)$ matrix-vector products with $\W_\z^{-1}$ are required, which costs $4(r+\ell_H)e_\z$. Including the high-fidelity forward solves, the total computational cost of our post-optimality analysis is
\begin{align*}
\tilde{P} = N f + (s+1+4r+4\ell_H) \tilde{a} + 2(q+\ell_E)e_\u + (s+4r+4 \ell_H + 2N) e_\z .
\end{align*}

To make the comparison concrete, we assume cost values given in Table~\ref{tab:cost_comparison}, algorithmic parameters given in Table~\ref{tab:cost_comparison_alg_params}, and iteration counts $\tilde{n}_\text{iter}=50$ and $\tilde{n}_\text{adjoint}=50$. Then we have the total computational costs
\begin{align*}
 \tilde{O} = 1.59 \times 10^4 \qquad \tilde{P} = 2.587 \times 10^3 .
\end{align*}
In this illustration, the post-optimality computation costs an order of magnitude less than the low-fidelity optimization. 

\begin{table}[!ht]
\centering
\begin{tabular}{|c|c|c|c|c|}
\hline
$f$ & $\tilde{f}$ & $\tilde{a}$ & $e_\u$ & $e_\z$ \\
$100$ & $15$ & $3$ & $1$ & $1$ \\
\hline
\end{tabular}
\caption{Cost values for illustrative computational cost comparison.}
\label{tab:cost_comparison}
\end{table}

\begin{table}[!ht]
\centering
\begin{tabular}{|c|c|c|c|c|c|}
\hline
$N$ & $s$ & $q$ & $r$ & $\ell_E$ & $\ell_H$ \\
$2$ & $100$ & $500$ & $50$ & $10$ & $10$ \\
\hline
\end{tabular}
\caption{Algorithm parameters for illustrative computational cost comparison.}
\label{tab:cost_comparison_alg_params}
\end{table}

\section*{Appendix D}
\begin{figure}[h]
\centering
  \includegraphics[width=0.49\textwidth]{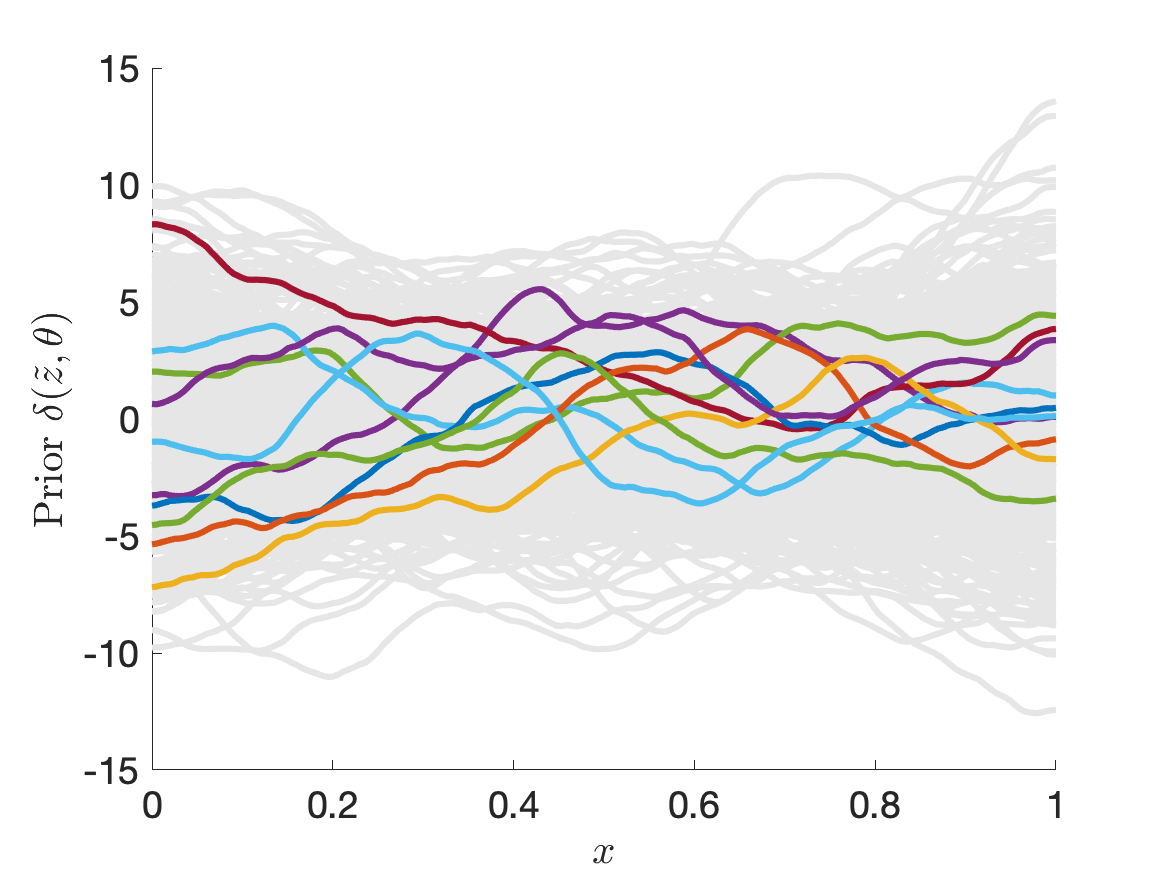}
    \includegraphics[width=0.44\textwidth]{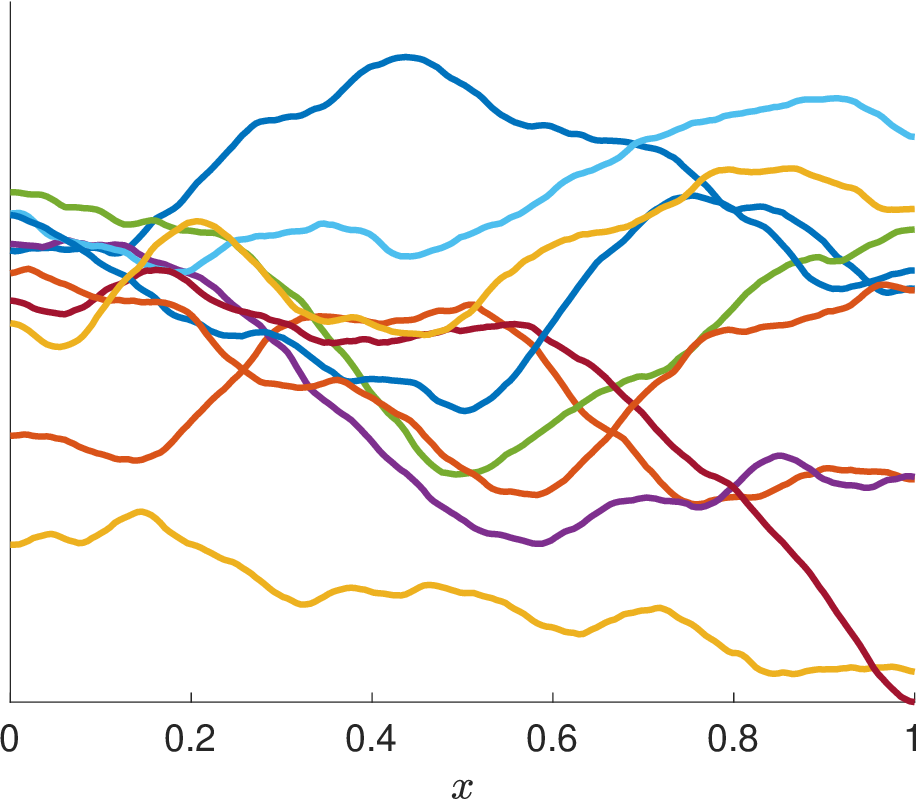}
    \caption{Left: samples of the prior discrepancy $\d(\z,\t)$ evaluated at $\z=\tilde{\z}$, i.e. a mean zero Gaussian distribution with covariance $\W_\u^{-1}$; right: samples from a mean zero Gaussian distribution with covariance $\W_\z^{-1}$.}
  \label{fig:illustrative_example_state_prior}
\end{figure}

This appendix provides a detailed illustration prior hyper-parameter selection for the diffusion reaction example presented in Section~\ref{subsec:illustrative_example}. We begin with specification of the correlation length hyper-parameter and variance for the state prior covariance $\W_\u^{-1}$. Recall that the prior discrepancy, evaluated at $\z=\ztilde$, follows a mean zero Gaussian distribution with covariance $\W_\u^{-1}.$ Using samples from $\mathcal N(\vec{0},\W_\u^{-1})$, we tune $\beta_\u$ so that the samples have a correlation length commensurate to that of the high-fidelity model. A similar process is followed tuning $\alpha_\u$ so that the prior discrepancy samples have a magnitude commensurate to the scales characteristic of $S(\z)-\tilde{S}(\z)$. This may be done by computing the correlation length and variance of $S(\z_\ell)-\tilde{S}(\z_\ell)$, but ideally should be based on some prior information about the high-fidelity model. Such prior information is typically available if the model discrepancy is related to unresolved processes whose correlation length and variance can be estimated from knowledge of the physical process. The left panel of Figure~\ref{fig:illustrative_example_state_prior} displays 500 prior samples of the discrepancy evaluated at $\ztilde$. We highlight the range of magnitude variability showing 490 samples with grey curves while demonstrating the correlation length by plotting 10 samples in contrasting colors. 

To define the correlation length hyper-parameter for the optimization variable prior covariance $\W_\z^{-1}$, we follow a similar approach to tune $\beta_\z$ such that the samples are representative of the correlation lengths expected for the source term. The right panel of Figure~\ref{fig:illustrative_example_state_prior} displays 10 samples from a mean zero Gaussian distribution with covariance $\W_\z^{-1}$. Note that the vertical axis is intentionally not labeled as the magnitude of the samples is determined by the variance hyper-parameter which is specified later. 

\begin{figure}[h]
\centering
  \includegraphics[width=0.49\textwidth]{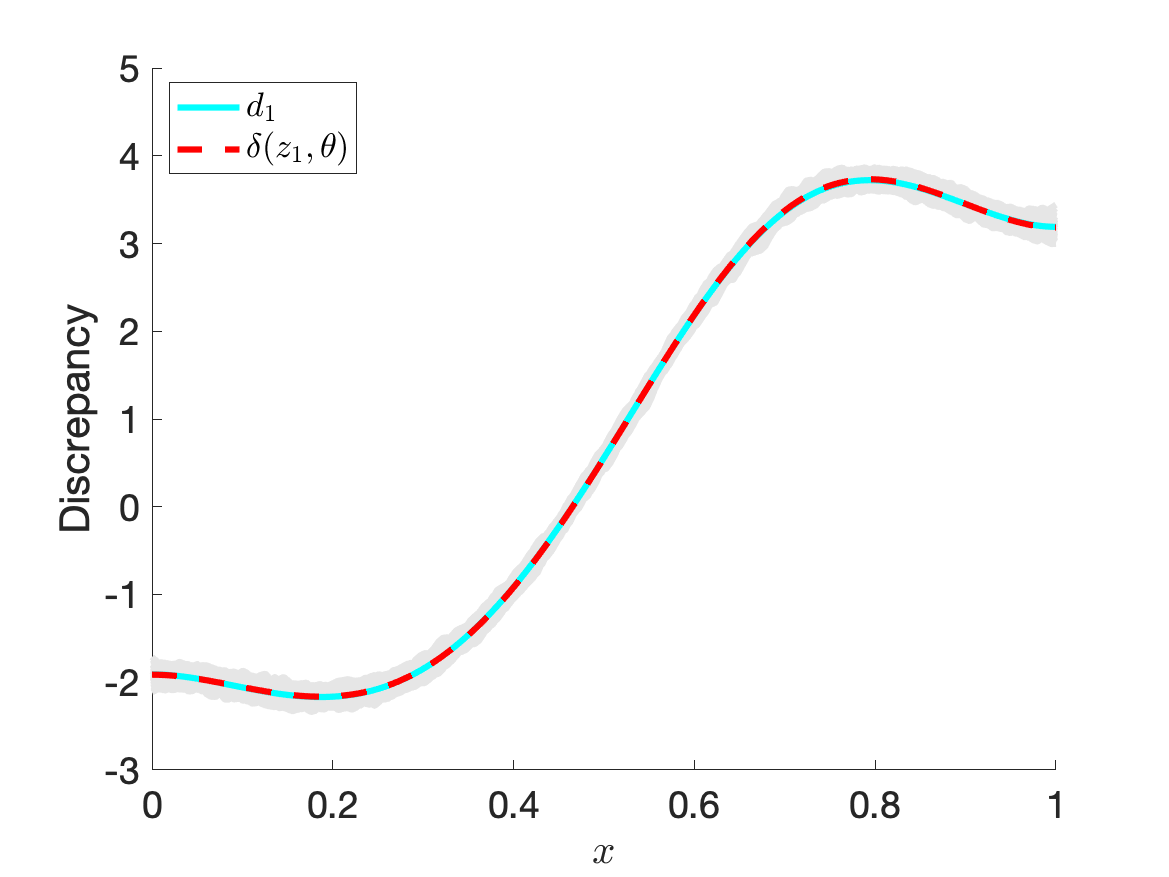}
    \includegraphics[width=0.49\textwidth]{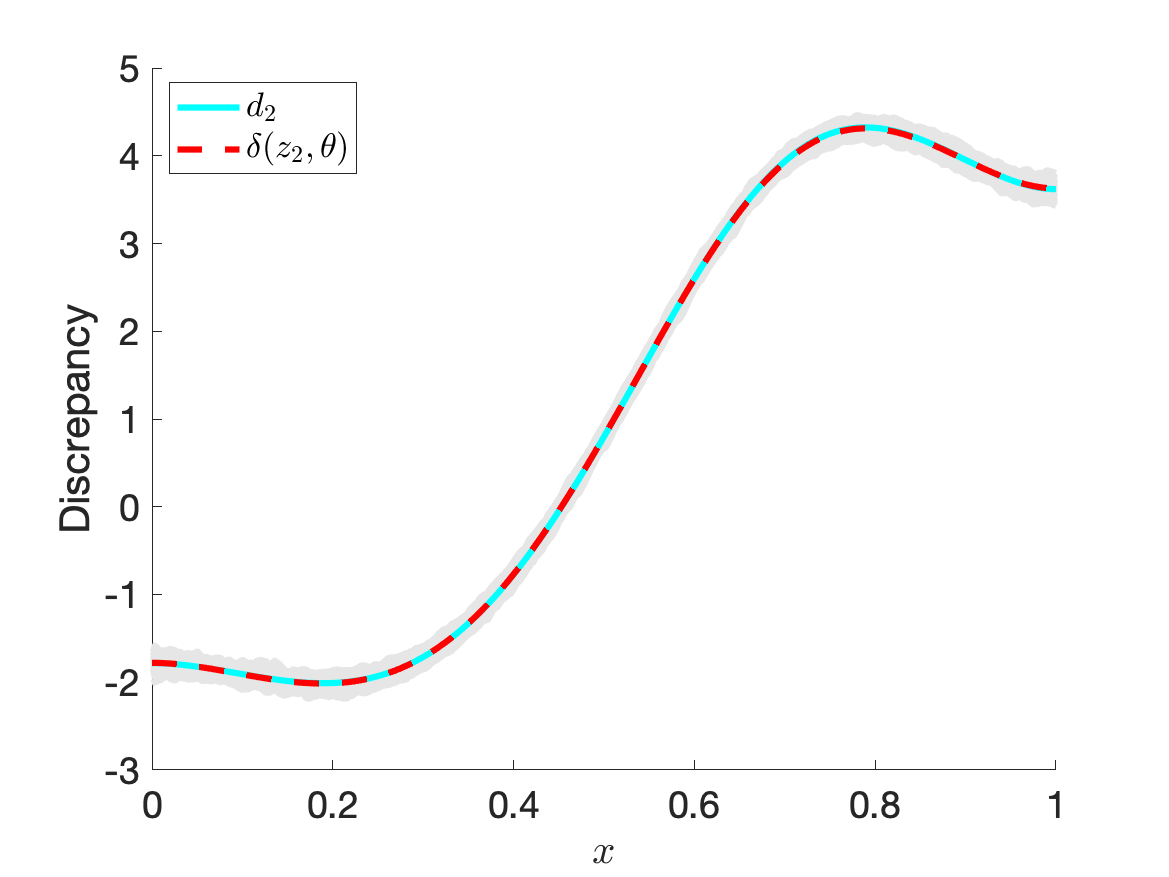}
  \includegraphics[width=0.49\textwidth]{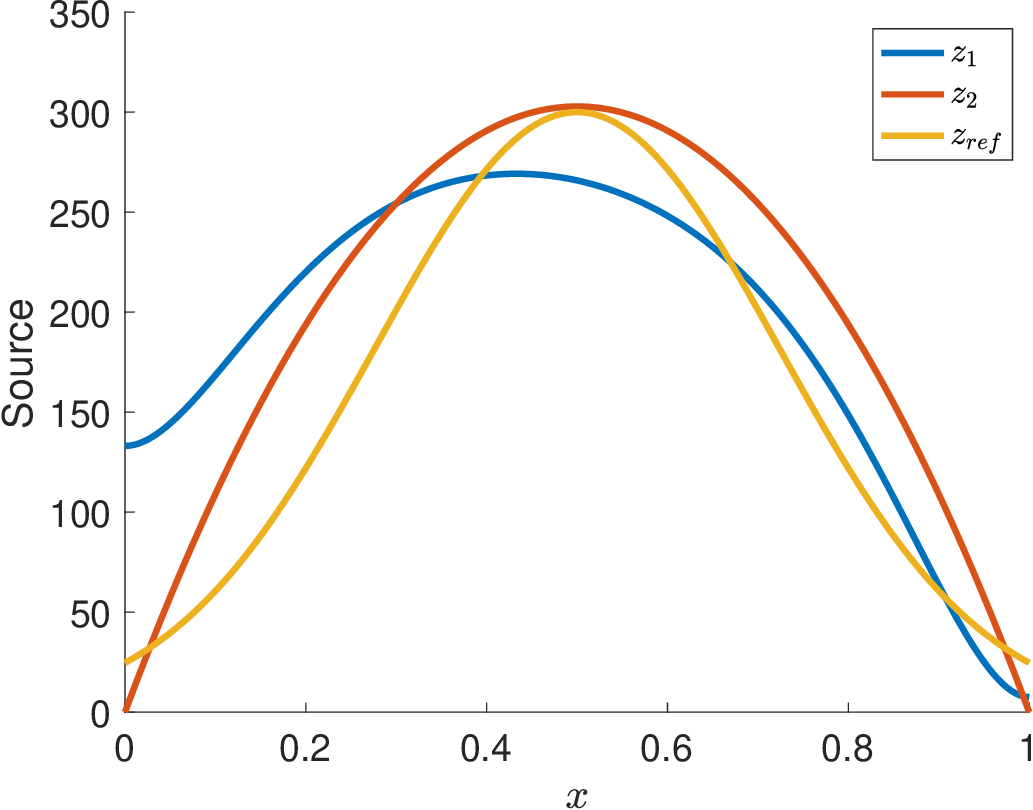}
        \includegraphics[width=0.49\textwidth]{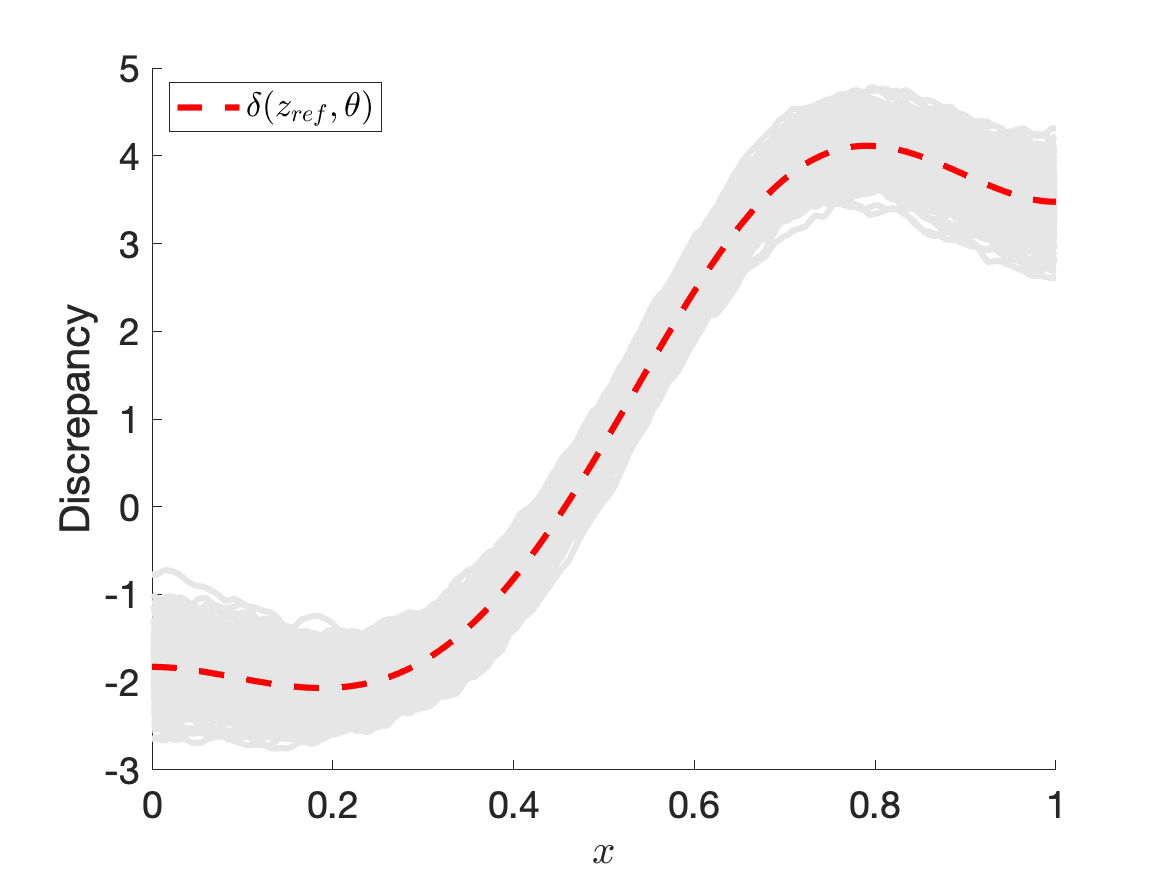} \\
    \caption{Top: discrepancy data $\vec{d}_\ell=S(\z_\ell)-\tilde{S}(\z_\ell)$, the posterior discrepancy mean evaluated at $\z_\ell$, and posterior discrepancy samples $\d(\z,\overline{\t}+\hat{\t}+\breve{\t})$ evaluated at $\z=\z_\ell$, $\ell=1,2$. Bottom left: optimization variable data $\z_1$, $\z_2$, and a reference source $\z_{ref}$; right: the posterior discrepancy $\d(\z,\overline{\t}+\hat{\t}+\breve{\t})$  evaluated at $\z=\z_{ref}$.}
  \label{fig:illustrative_example_discrepancy_posterior}
\end{figure}

It remains to specify the optimization variable variance $\alpha_\z$ and the noise variance $\alpha_{\vec{d}}$. To do this, we consider the posterior discrepancy distribution. The noise variance $\alpha_{\vec{d}}$ will influence how well the posterior discrepancy fits the data, while the optimization variable variance $\alpha_\z$ determines the uncertainty magnitude in directions not informed by the data. Figure~\ref{fig:illustrative_example_discrepancy_posterior} shows posterior discrepancy samples used to tune $\alpha_{\vec{d}}$ and $\alpha_\z$. Specifically, the top row of Figure~\ref{fig:illustrative_example_discrepancy_posterior} shows the discrepancy data $\vec{d}_\ell=S(\z_\ell)-\tilde{S}(\z_\ell)$, the posterior discrepancy mean evaluated at $\z_\ell$, and posterior discrepancy samples evaluated at $\z_\ell$, $\ell=1,2$. The noise variance was chosen to ensure a quality fit of the data with a relatively small posterior variance. We see a weak dependence of the discrepancy on $\z$ as $\vec{d}_1$ and $\vec{d}_2$ are similar, with slight differences in their maximum magnitudes. The bottom left panel of Figure~\ref{fig:illustrative_example_discrepancy_posterior} shows the optimization variable data $\z_1$ and $\z_2$, along with another source denoted as $\z_{ref}$, which is defined by the user as a ``reference". The role of $\z_{ref}$ is to assess the posterior uncertainty for a source term which is physically realistic but is not contained in the training data $\{\z_\ell\}_{\ell=1}^N$. The bottom right panel of Figure~\ref{fig:illustrative_example_discrepancy_posterior} displays the posterior discrepancy evaluated at $\z_{ref}$. The optimization variable prior variance $\alpha_\z$ was tuned to ensure that the posterior variance of $\d(\z_{ref},\t)$ has a reasonable magnitude of uncertainty. 

\bibliographystyle{amsplain}
\bibliography{dasco}

\providecommand{\bysame}{\leavevmode\hbox to3em{\hrulefill}\thinspace}
\providecommand{\MR}{\relax\ifhmode\unskip\space\fi MR }
\providecommand{\MRhref}[2]{%
  \href{http://www.ams.org/mathscinet-getitem?mr=#1}{#2}
}
\providecommand{\href}[2]{#2}
\begin{thebibliography}{10}

\bibitem{biegler_rom_opt}
Anshul Agarwal and Lorenz Biegler, \emph{A trust-region framework for
  constrained optimization using reduced order modeling}, Optimization and
  Engineering \textbf{14} (2013), 3--35.

\bibitem{trmm_lewis}
N.M Alexandrov, J.E.~Dennis {Jr.}, R.M. Lewis, and V.~Torczon, \emph{A
  trust-region framework for managing the use of approxima- tion models in
  optimization}, Structural Optimization \textbf{15} (1998), 16--23.

\bibitem{frontier_in_pdeco}
H.~Antil, D.~P. Kouri, M.~Lacasse, and D.~Ridzal (eds.), \emph{Frontiers in
  pde-constrained optimization}, Springer, 2018.

\bibitem{Arendt_2012}
Paul~D. Arendt, Daniel~W. Apley, and Wei Chen, \emph{Quantification of model
  uncertainty: Calibration, model discrepancy, and identifiability}, Journal of
  Mechanical Design \textbf{134} (2012).

\bibitem{Archer_01}
U.~M. Ascher and E.~Haber, \emph{Grid refinement and scaling for distributed
  parameter estimation problems}, Inverse Problems \textbf{17} (2001),
  571--590.

\bibitem{Biegler_11}
L.~Biegler, G.~Biros, O.~Ghattas, M.~Heinkenschloss, D.~Keyes, B.~Mallick,
  Y.~Marzouk, L.~Tenorio, B.~van Bloemen~Waanders, and K.~Willcox (eds.),
  \emph{Large-scale inverse problems and quantification of uncertainty}, John
  Wiley and Sons, 2011.

\bibitem{Biegler_07}
L.~T. Biegler, O.~Ghattas, M.~Heinkenschloss, D.~Keyes, and B.~van
  Bloemen~Waanders (eds.), \emph{Real-time pde-constrained optimization},
  vol.~3, SIAM Computational Science and Engineering, 2007.

\bibitem{Biros_05}
G.~Biros and O.~Ghattas, \emph{Parallel {L}agrange-{N}ewton-{K}rylov-{S}chur
  methods for {PDE}-constrained optimization. {Parts I-II}}, SIAM J. Sci.
  Comput. \textbf{27} (2005), 687-- 738.

\bibitem{shapiro_SIAM_review}
J.~Fr{\'e}d{\'e}ric Bonnans and Alexander Shapiro, \emph{Optimization problems
  with perturbations: A guided tour}, SIAM Review \textbf{40} (1998), no.~2,
  228--264.

\bibitem{Borzi_07}
A.~Borzi, \emph{High-order discretization and multigrid solution of elliptic
  nonlinear constrained optimal control problems}, J. Comp. Applied Math
  \textbf{200} (2007), 67--85.

\bibitem{griesse2}
Kerstin Brandes and Roland Griesse, \emph{Quantitative stability analysis of
  optimal solutions in {PDE}-constrained optimization}, Journal of
  Computational and Applied Mathematics \textbf{206} (2007), 908--926.

\bibitem{ohagan2014}
Jenn{\'y} Brynjarsd{\'o}ttir and Anthony O'Hagan, \emph{Learning about physical
  parameters: the importance of model discrepancy}, Inverse Problems
  \textbf{30} (2014), no.~11.

\bibitem{multifidelity_quasinewton_bryson}
Dean~{E.} Bryson and Markus~{P.} Rumpfkeil, \emph{Multifidelity quasi-newton
  method for design optimization}, AIAA JOURNAL \textbf{56} (2018), no.~10.

\bibitem{gentle_intro_bayes}
Tan Bui-Thanh, \emph{A gentle tutorial on statistical inversion using the
  bayesian paradigm}, Tech. report, The Institute for Computational Engineering
  and Sciences, 2012.

\bibitem{ghattas_infinite_dim_bayes_1}
Tan Bui-{T}hanh, Omar Ghattas, James Martin, and Georg Stadler, \emph{A
  computational framework for infinite-dimensional bayesian inverse problems.
  {P}art {I}: The linearized case, with applications to global seismic
  inversion}, SIAM Journal on Scientific Computing \textbf{35} (2013), no.~6,
  A2494--A2523.

\bibitem{stadler_cov_bc}
Y.~Daon and G.~Stadler, \emph{Mitigating the influence of the boundary on
  pde-based covariance operators}, Inverse Problems and Imaging \textbf{12}
  (2018), no.~5, 1083--1102.

\bibitem{fiacco}
Anthony~V. Fiacco and Abolfazl Ghaemi, \emph{Sensitivity analysis of a
  nonlinear structural design problem}, Computers and Operations Research
  \textbf{9} (1982), no.~1, 29--55.

\bibitem{Gardner_2021}
P.~Gardner, T.J. Rogers, C.~Lord, and R.J. Barthorpe, \emph{Learning model
  discrepancy: A gaussian process and sampling-based approach}, Mechanical
  Systems and Signal Processing \textbf{152} (2021), no.~107381.

\bibitem{multifidelity_gorodetsky}
Alex~{A.} Gorodetsky, Gianluca Geraci, Michael~S. Eldred, and John~{D.}
  Jakeman, \emph{A generalized approximate control variate framework for
  multifidelity uncertainty quantification}, Journal of Computational Physics
  \textbf{408} (2020), no.~1.

\bibitem{Griesse_part_2}
Roland Griesse, \emph{Parametric sensitivity analysis in optimal control of a
  reaction-diffusion system -- part {II}: practical methods and examples},
  Optimization Methods and Software \textbf{19} (2004), no.~2, 217--242.

\bibitem{Griesse_part_1}
\bysame, \emph{Parametric sensitivity analysis in optimal control of a reaction
  diffusion system. {I}. solution differentiability}, Numerical Functional
  Analysis and Optimization \textbf{25} (2004), no.~1-2, 93--117.

\bibitem{Haber_01}
E.~Haber and U.~M. Ascher, \emph{Preconditioned all-at-once methods for large,
  sparse parameter estimation problems}, Inverse Problems \textbf{17} (2001),
  1847--1864.

\bibitem{hart_2021_bayes}
Joseph Hart and Bart {van Bloemen Waanders}, \emph{Enabling hyper-differential
  sensitivity analysis for ill-posed inverse problems}, SIAM Journal on
  Scientific Computing \textbf{45} (2023), no.~4.

\bibitem{hart_bvw_cmame}
\bysame, \emph{Hyper-differential sensitivity analysis with respect to model
  discrepancy: Optimal solution updating}, Computer Methods in Applied
  Mechanics and Engineering \textbf{412} (2023).

\bibitem{HDSA}
Joseph Hart, Bart {van Bloemen Waanders}, and Roland Hertzog,
  \emph{Hyper-differential sensitivity analysis of uncertain parameters in
  {PDE}-constrained optimization}, International Journal for Uncertainty
  Quantification \textbf{10} (2020), no.~3, 225--248.

\bibitem{hart_hdsa_control_oed}
Joseph Hart, Bart {van Bloemen Waanders}, Lisa Hood, and Julie Parish,
  \emph{Sensitivity driven experimental design to facilitate control of
  dynamical systems}, Submitted. arXiv: https://arxiv.org/abs/2202.03312
  (2022).

\bibitem{Hazra_06}
S.~B. Hazra and V.~Schulz, \emph{Simultaneous pseudo-timestepping for
  aerodynamic shape optimization problems with state constraints}, SIAM J. Sci.
  Comput. \textbf{28} (2006), 1078--1099.

\bibitem{Higdon_2008}
Dave Higdon, James Gattiker, Brian Williams, and Maria Rightley, \emph{Computer
  model calibration using high-dimensional output}, Journal of the American
  Statistical Association \textbf{103} (2008), no.~482, 570--583.

\bibitem{Hintermuller_05}
M.~Hintermuller and L.~N. Vicente, \emph{Space mapping for optimal control of
  partial differential equations}, SIAM J. Opt. \textbf{15} (2005), 1002--1025.

\bibitem{Hinze_09}
M.~Hinze, R.~Pinnau, M.~Ulbrich, and S.~Ulbrich, \emph{Optimization with pde
  constraints}, Springer, 2009.

\bibitem{ohagan2001}
Marc~C. Kennedy and Anthony O'Hagan, \emph{Bayesian calibration of computer
  models}, Journal of the Royal Statistical Society \textbf{63} (2001), no.~3,
  425--464.

\bibitem{khristenko19}
Ustim Khristenko, Laura Scarabosio, Piotr Swierczynski, E.~Ullmann, and Barbara
  Wohlmuth, \emph{Analysis of boundary effects on {PDE}-based sampling of
  {W}hittle--{M}atérn random fields}, SIAM/ASA Journal on Uncertainty
  Quantification \textbf{7} (2019), 948--974.

\bibitem{Biegler_03}
M.~Heinkenschloss L.~T.~Biegler, O.~Ghattas and B.~van Bloemen~Waanders (eds.),
  \emph{Large-scale pde-constrained optimization}, vol.~30, Springer-Verlag
  Lecture Notes in Computational Science and Engineering, 2003.

\bibitem{Laird_05}
C.~D. Laird, L.~T. Biegler, B.~van Bloemen~Waanders, and R.~A. Bartlett,
  \emph{Time dependent contaminant source determination for municipal water
  networks using large scale optimization}, ASCE J. Water Res. Mgt. Plan.
  (2005), 125--134.

\bibitem{lindgren11}
Finn Lindgren, Håvard Rue, and Johan Lindström, \emph{An explicit link
  between {G}aussian fields and {G}aussian {M}arkov random fields: The
  stochastic partial differential equation approach}, Journal of the Royal
  Statistical Society Series B \textbf{73} (2011), 423--498.

\bibitem{Ling_2014}
You Ling, Joshua Mullins, and Sankaran Mahadevan, \emph{Selection of model
  discrepancy priors in bayesian calibration}, Journal of Computational Physics
  \textbf{276} (2014), 665--680.

\bibitem{willcox_multifi_opt_2012}
Andrew March and Karen Willcox, \emph{Provably convergent multifidelity
  optimization algorithm not requiring high-fidelity derivatives}, AIAA Journal
  \textbf{50} (2012), no.~5, 1079--1089.

\bibitem{Maupin}
Kathryn~A. Maupin and Laura~P. Swiler, \emph{Model discrepancy calibration
  across experimental settings}, Reliability Engineering \& System Safety
  \textbf{200} (2020).

\bibitem{multifidelity_wind_plant_geraci}
Friedrich Menhorn, Gianluca Geraci, Daniel~{T.} Seidl, Michael~S. Eldred, Ryan
  King, Hans-Joachim Bungartz, and Youssef Marzouk, \emph{Higher moment
  multilevel estimators for optimization under uncertainty applied to wind
  plant design}, AIAA SciTech Forum, 2020, pp.~1--16.

\bibitem{multi_fidelity_opt_willcox}
Leo~{W. T.} Ng and Karen~E. Willcox, \emph{Multifidelity approaches for
  optimization under uncertainty}, International Journal for Numerical Methods
  in Engineering (2014).

\bibitem{ohagan2002}
Jeremy Oakley and Anthony O'Hagan, \emph{Bayesian inference for the uncertainty
  distribution of computer model outputs}, Biometrika \textbf{89} (2002),
  no.~4, 769--784.

\bibitem{multifidelity_review_peherstorfer}
Benjamin Peherstorfer, Karen Willcox, and Max Gunzburger, \emph{Survey of
  multifidelity methods in uncertainty propagation, inference, and
  optimization}, SIAM Review \textbf{60} (2018), no.~3, 550--591.

\bibitem{ghattas_infinite_dim_bayes_2}
Noemi Petra, James Martin, Georg Stadler, and Omar Ghattas, \emph{A
  computational framework for infinite-dimensional bayesian inverse problems.
  {P}art {II}: Stochastic newton mcmc with application to ice sheet flow
  inverse problems}, SIAM Journal on Scientific Computing \textbf{36} (2014),
  no.~4, A1525--A1555.

\bibitem{saibaba_gsvd}
Arvind~{K.} Saibaba, Joseph Hart, and Bart {van Bloemen Waanders},
  \emph{Randomized algorithms for generalized singular value decomposition with
  application to sensitivity analysis}, Numerical Linear Algebra with
  Applications (2021).

\bibitem{stuart_inv_prob}
Andrew~M. Stuart, \emph{Inverse problems: a {B}ayesian perspective}, Acta
  Numerica (2010).

\bibitem{sunseri_2}
Isaac Sunseri, Alen Alexanderian, Joseph Hart, and Bart {van Bloemen Waanders},
  \emph{Hyper-differential sensitivity analysis for nonlinear bayesian inverse
  problems}, Submitted. arXiv: https://arxiv.org/abs/2202.02219 (2022).

\bibitem{sunseri_hdsa}
Isaac Sunseri, Joseph Hart, Bart {van Bloemen Waanders}, and Alen Alexanderian,
  \emph{Hyper-differential sensitivity analysis for inverse problems
  constrained by partial differential equations}, Inverse Problems \textbf{36}
  (2020), no.~12.

\bibitem{Vogel_99}
C.~R. Vogel, \emph{Sparse matrix computations arising in distributed parameter
  identification}, SIAM J. Matrix Anal. Appl. (1999), 1027--1037.

\bibitem{Vogel_02}
\bysame, \emph{Computational methods for inverse problems}, SIAM Frontiers in
  Applied Mathematics Series, 2002.

\bibitem{whittle54}
Peter Whittle, \emph{On stationary processes in the plane}, Biometrika
  \textbf{41} (1954).

\end{thebibliography}

\end{document}